\documentclass[9pt]{amsart}
\usepackage{amsmath, amssymb}
\usepackage{amsfonts,amscd}
\usepackage[a4paper, margin=1.25in, headsep=0.4in]{geometry}
\usepackage[parfill]{parskip}
\usepackage{xcolor}
\usepackage{comment}
\usepackage{mathtools}
\usepackage{commath}
\usepackage{appendix}
\usepackage{needspace}
\usepackage{braket}

\usepackage{xcolor}
\usepackage{comment}

\usepackage{bm}

\usepackage{mathrsfs}

\usepackage{enumitem}

\usepackage{booktabs, xltabular, array}

\setlist[enumerate]{labelindent=0pt, leftmargin=*}
\setlist[itemize]{labelindent=0pt, leftmargin=*}
\usepackage{hyperref}
\hypersetup{
    colorlinks,
    citecolor=red,
    filecolor=black,
    linkcolor=blue,
    urlcolor=black
}

\makeatletter
\def\@tocline#1#2#3#4#5#6#7{\relax
  \ifnum #1>\c@tocdepth
  \else
    \par\addpenalty\@secpenalty\addvspace{#2}%
    \begingroup
    \fontsize{9pt}{7pt}\selectfont
    \hyphenpenalty\@M
    \@ifempty{#4}{%
      \@tempdima\csname r@tocindent\number#1\endcsname\relax
    }{%
      \@tempdima#4\relax
    }%
    \parindent\z@
    \leftskip#3\relax
    \ifnum#1=2
      \advance\leftskip by 1.75em
    \fi
    \advance\leftskip\@tempdima\relax
    \rightskip\@pnumwidth plus4em
    \parfillskip-\@pnumwidth
    #5\leavevmode\hskip-\@tempdima
    \ifcase#1
      \or
      \or
      \or
      \else\hskip 3em
    \fi
    #6\nobreak\relax
    \ifnum#1=2
      \nobreak\leaders\hbox{.\kern.1em}\hfill
    \else
      \hfill
    \fi
    \hbox to\@pnumwidth{\@tocpagenum{#7}}\par
    \nobreak
    \endgroup
  \fi
}
\makeatother

\colorlet{eqcolor}{blue}
\makeatletter
\renewcommand{\tagform@}[1]{%
    \maketag@@@{\textcolor{eqcolor}{(\ignorespaces#1\unskip)}}%
}
\makeatother

\usepackage{setspace}
\usepackage{amsthm}
\allowdisplaybreaks

\newtheorem{theorem}{Theorem}

\newtheorem{prop}[theorem]{Proposition}
\newtheorem{lemma}[theorem]{Lemma}

\newtheorem{cor}[theorem]{Corollary}

\newtheorem{thm}[theorem]{Theorem}

\theoremstyle{definition}
\newtheorem{remark}[theorem]{Remark}
\newtheorem{setup}[theorem]{Setup}
\newtheorem{defn}[theorem]{Definition}

\newtheorem*{unmremark}{Remark}

\newcommand{\f}{\mathcal F}

\def\Z{\mathbb Z}
\def\N{\mathbb N}
\def\Q{\mathbb Q}
\def\R{\mathbb R}
\def\C{\mathbb C}
\def\F{\mathbb F}

\def\F{\mathbb F}

\let\hungarianumlaut\H 
\def\H{\mathbb H}
\def\E{\mathbb E}

\def\dim{\operatorname{dim}}

\def\O{\operatorname{O}}
\def\o{\operatorname{o}}
\def\log{\operatorname{log}}

\def\p{\operatorname{\preccurlyeq}}

\def\thtp{\operatorname{\theta_{\textit{f\,}}(\textit{p})}}

\def\thtpi{\operatorname{\theta_{\textit{f\,}}(\textit{p}_{\textit{i}})}}

\def\Bpm{S_M^{\pm}}

\def\chi{\mathbf 1}

 \newcommand{\tmod}{{\mathrm{mod} \, }}

\makeatletter
\renewenvironment{proof}[1][\proofname]{%
  \par\pushQED{\qed}%
  \normalfont            
  \topsep6pt \trivlist
  \item[\hskip\labelsep
        \bfseries #1\@addpunct{.}]
}{%
  \popQED\endtrivlist\@endpefalse
}
\makeatother

\makeatletter
\let\oldthebibliography\thebibliography
\renewcommand{\thebibliography}[1]{%
  \oldthebibliography{#1}%
  \small
  \setlength{\itemsep}{3pt}%
  \setlength{\parsep}{0pt}%
  \setlength{\labelwidth}{5.2em}%
  \setlength{\leftmargin}{\dimexpr\labelwidth+\labelsep\relax}%
}
\makeatother

\allowdisplaybreaks

\begin{document}

\title[Central limit theorems: Frobenius, Hecke and Kloosterman angles]{Central limit theorems for  Frobenius, Hecke and Kloosterman angles}

\author[Shreyas Pai]{Shreyas Pai}
\address{Shreyas Pai, Department of Mathematics, National University of Singapore, 10 Lower Kent Ridge Road, Singapore 119076}
\email{shreyas.pai@u.nus.edu}

\author[Kaneenika Sinha]{Kaneenika Sinha}
\address{Kaneenika Sinha, Department of Mathematics, IISER Pune, Dr Homi Bhabha Road, Pashan, Pune - 411008, Maharashtra, India}
\email{kaneenika@iiserpune.ac.in}
\author[Naman Thakkar]{Naman Thakkar}
\address{Naman Thakkar, Department of Mathematics, IISER Pune, Dr Homi Bhabha Road, Pashan, Pune - 411008, Maharashtra, India}
\email{thakkar.namanketan@students.iiserpune.ac.in}
\begin{abstract} 
  The Frobenius angles of non-CM elliptic curves and the Hecke angles of non-CM newforms are equidistributed with respect to the Sato--Tate measure, while for Kloosterman angles, this is known on average over natural families.  We study the fluctuations around these laws in large families, for smooth test functions.  We prove an unconditional central limit theorem for Frobenius angles of elliptic curves for a much larger class of test functions than previously considered; earlier results assumed the Riemann hypothesis for $L$-functions of Hecke eigenforms of full level.  This result follows from a general framework, developed here, which yields Gaussian moments for statistics of sequences in $[0,1]$ varying over families. In this framework, independence is replaced by two quantitative conditions: approximate independence at distinct indices, and convergence of the family averages to the Weyl limits at a suitable rate.  As further applications, we extend central limit theorems for Kloosterman angles and for Hecke angles of newforms to larger classes of test functions, the latter with the weight and the level varying simultaneously.  We also treat the dilated angles $\alpha\theta$ for integers $\alpha\geq2$, whose limiting distribution is uniform.  Finally, we discuss sharp statistics, obtaining a central limit theorem for Hecke angles and explaining the obstruction for the other two families. 
\end{abstract}
\bigskip



\maketitle

    \tableofcontents

\section{Introduction}

Many arithmetic quantities attached to primes, such as the traces of Frobenius of an elliptic curve, the Fourier coefficients of a modular form, and Kloosterman sums, behave in several aspects like random variables. The Sato–Tate law is a natural instance of this. For a non-CM Hecke newform $f$, write its normalized $p$-th Fourier coefficient as $a_f(p)=2\cos(\pi\theta_f(p))$ with $\theta_f(p)\in[0,1]$.   The angles $\theta_f(p)$ are equidistributed with respect to the measure $2\sin^2(\pi t)\,dt$ as $p\to\infty$.  This was conjectured by Sato and Tate in the 1960s for elliptic curves (see \cite{Tate}), extended to modular forms by Serre \cite{Serre1}, and is now a theorem \cite{Clozel-Harris-Taylor}, \cite{BGHT}.  Similarly, the normalized Kloosterman sums $\displaystyle K(p,a)/\sqrt p= 2\cos(\pi\theta_{p,a})$ define angles $\theta_{p,a}\in[0,1]$.  Katz \cite{Katz} showed that for a fixed prime $p$, these angles are Sato–Tate distributed as $a$ varies over $\mathbb F_p^\times$, in the limit $p\to\infty$, and conjectured the same for a fixed $a$ as $p$ varies; the latter remains open.  A natural question concerns the fluctuations around these laws: for a test function $h$, how does $\sum_{p\leq x}h(\theta_f(p))$ deviate from its expected value $\pi(x)\int_0^1h(t)\,2\sin^2(\pi t)\,dt$? For an individual form, this is a difficult question. Averaging over families makes such questions accessible. For example, in the years preceding the spectacular resolution of the Sato--Tate conjecture,   equidistribution questions for Hecke angles were approached by averaging over families of newforms.  The conjecture of Sato and Tate is viewed as a horizontal law (where $f$ is fixed and the primes vary).  In the {\it vertical} aspect, a prime $p$ is fixed and $f$ varies over increasing families (see \cite{Sarnak},  \cite{Serre}, \cite{CDF}, \cite{MS}, \cite{LW}).   In the {\it average} aspect, both $p$ and $f$ vary, with $f$ ranging over the set of newforms $\mathcal F_{N,k}$ of even weight $k$ with respect to $\Gamma_0(N)$ (see \cite{CDF}, \cite{Nagoshi}).  These results show that the Sato--Tate law holds on average over families with rapidly growing weights.  Moreover, Nagoshi \cite{Nagoshi} demonstrated an Erd\hungarianumlaut{o}s--Kac type theorem in this context: averaged over large families of Hecke newforms, suitably normalized sums of $a_f(p)$ over primes $p \leq x$ converge to the Gaussian distribution.  

Analogous questions have been studied for families of elliptic curves $E(a,b):\,y^2 = x^3 + ax +b$ with $a$ and $b$ ranging over integers in large boxes, where the relevant sequences are the Frobenius angles $\theta_E(p)$ (see \cite{FM}, \cite{DP}, \cite{JY}, \cite{Banks-Shpar}, \cite{Baier-Zhao}, \cite{Baier-Prabhu}, \cite{BPS}).   In particular, \cite{BPS} established central limit theorems for smooth statistics of these angles under the assumption of the Riemann Hypothesis for all $L$-functions associated to Hecke eigenforms with respect to the full modular group.  Our first main result removes the above assumption, and derives an unconditional central limit theorem for a much larger class of test functions than those considered in \cite{BPS}.  The key observation is that, for a fixed test function, the decay of its Fourier coefficients makes the following bound at each individual prime sufficient: for each prime $p$ and each $k \geq 1$, the average of $a_E(p^{2k})$ over all curves $E(a,b)$ with $a,b \in \mathbb F_p$ is $\O(k/\sqrt{p})$.  This follows from Birch's formulas (Proposition \ref{S(n)-elliptic}\ref{S(n)}-\ref{S0pm}) and the Ramanujan--Deligne bound (Proposition \ref{S(n)-elliptic}\ref{Spm-bounds}).  The cancellation over primes that \cite{BPS} obtained from the modular Riemann hypothesis is then not needed.

This is one among three applications of a general framework (Section \ref{sec: General Framework}). For Kloosterman angles, we extend the central limit theorem of Xi \cite[Theorem 2]{Xi} to test functions of limited smoothness, to sums of the normalized Kloosterman sums $K(p,a)/\sqrt p$ themselves, and to dilated angles, with less arithmetic input. For Hecke angles, we recover the known fixed-level results, extend them to test functions with only polynomial decay of the Fourier transform, and to dilated angles.  We also allow the weight and the level to vary simultaneously.  Detailed comparisons with earlier work are given in the remarks following Theorems \ref{Elliptic-curves-smooth}, \ref{Kloosterman-multiples-smooth} and \ref{Hecke-multiples-smooth}.

\subsection{A CLT for Frobenius angles}\label{sec: Main theorem Frobenius angles}
\begin{thm}\label{Elliptic-curves-smooth}
  Let $g$ be a continuous, real-valued and even function on $\R$ such that for some $\delta > 1$ and some $\kappa > 2$,
  $$g(t) \ll \frac{1}{(1 + |t|)^{\delta}} \quad  
  \text{ and }\quad 
  \widehat{g}(t) \ll \frac{1}{(1 + |t|)^{\kappa}} .$$
  Fix a real number $L \geq 1$, and define

$$h(t) := \sum_{n \in \Z}g(L(t+n)).$$
Assume that $h$ is not constant.  For a positive integer $\alpha$, define
$$\mu_{\alpha}(t) := \begin{cases}
    2\sin^2(\pi t) &\text{ if }\alpha = 1\\
    1 &\text{ if }\alpha \geq 2,
\end{cases}$$
and
$$V_h(\mu_{\alpha}) := \int_0^1 h(t)^2 \mu_{\alpha}(t) dt - \left(\int_0^1 h(t) \mu_{\alpha}(t) dt\right)^2. $$
For $(a,b) \in \Z^2$ such that $\Delta(a,b) := 4a^3 + 27b^2 \neq 0$, let $E = E(a,b)$ denote the elliptic curve defined by the equation
$Y^2 = X^3 + aX +b.$
For $A, B \geq 1$, define
$$\mathcal F_{A,B} := \{E = E (a,b): \,|a| \leq A,\,|b| \leq B,\,\Delta(a,b) \neq 0\},$$
which we regard as indexed by the pairs $(a,b)$, so that pairs giving isomorphic curves are counted separately.
 For a prime $p$, let $\widetilde{a}_E(p):= p+ 1 - N_E(p)$, where $N_E(p)$ is the number of points on the projective curve
$$Y^2Z = X^3 + aXZ^2 + bZ^3$$
over $\mathbb F_p$, the finite field of order $p$.  
Choose $\theta_E(p) \in [0,1]$ such that
$\displaystyle 2 \cos (\pi \theta_E(p)) = \frac{\widetilde{a}_E(p)}{\sqrt{p}}$. 

For $E \in \mathcal F_{A,B}$,
define
$$\mathcal S_x(E) := \frac{1}{\sqrt{{\pi}(x)V_{h}(\mu_{\alpha}})}\left(\sum_{p \leq x }h(\alpha\theta_E(p)) - {\pi}(x)\int_0^1 h(t)\mu_{\alpha}(t) dt\right).$$
If $A = A(x),\,B= B(x)$ satisfy $$\frac{\log A}{\log x},\,\frac{\log B}{\log x} \to \infty \text{ as }x \to \infty, $$
then, for all positive integers $r$,
$$\lim_{x \to \infty}\frac{1}{|\mathcal F_{A,B}|} \sum_{E \in \mathcal F_{A,B}}\mathcal S_x(E)^r  = \gamma_r :=\frac{1}{\sqrt{2\pi}}\int_{-\infty}^{\infty}  t^re^{-t^2/2}dt.$$

\end{thm}

\begin{unmremark}\leavevmode

\begin{enumerate}

\item In other words, Theorem \ref{Elliptic-curves-smooth} shows that viewing $\mathcal S_x$ as a random variable on $\mathcal F_{A,B}$ with the uniform probability measure, all its moments converge to the Gaussian moments.  Since the Gaussian distribution is determined by its moments, $\mathcal S_x$ converges in distribution to the standard Gaussian.
\item Theorem \ref{Elliptic-curves-smooth} strengthens previous results in multiple ways:
   \begin{enumerate}[label=(\alph*)]

   \item  The assumption of the Modular Riemann Hypothesis (MRH), that is, the Riemann Hypothesis for all $L$-functions associated to Hecke eigenforms $f$ with respect to the full modular group, is essential to the method of \cite{BPS}.  Theorem \ref{Elliptic-curves-smooth} bypasses this requirement, and provides an unconditional CLT.  Within the general framework developed in Section \ref{sec: Smooth statistics} (Theorems \ref{CLT general smooth} and \ref{CLT general smooth variant}), the decay of $\widehat{h}$ makes the bound of Proposition \ref{S(n)-elliptic}\ref{Spm-bounds} at individual primes sufficient. 
       \item The test functions $g$ considered in \cite[Theorem 1.5]{BPS} were (real-valued, even) Schwartz functions such that $\widehat{g}(t) \ll e^{-\lambda |t|^{\gamma}}$ for some $\lambda, \gamma >0$ as $|t| \to \infty$. Theorem \ref{Elliptic-curves-smooth} requires only polynomial decay for $g$ and $\widehat{g}$, namely  $\displaystyle g(t) \ll \frac{1}{(1 + |t|)^{\delta}}$ and $\displaystyle 
  \widehat{g}(t) \ll \frac{1}{(1 + |t|)^{\kappa}} $ for some $\delta >1$ and $\kappa >2$.  In particular, it holds for all real-valued, even Schwartz functions. It also covers test functions of limited smoothness that are not $C^{\infty}$.

     \item  Theorem \ref{Elliptic-curves-smooth}  treats the multiples $(\alpha \theta_E(p) )_{p \to \infty}$ for an integer $\alpha \geq 2$.  As a consequence of the Sato--Tate distribution of the sequence $(\theta_E(p))_{p \to \infty}$, the sequence $(\alpha \theta_E(p))_{p \to \infty}$ can be seen to be uniformly distributed modulo 1.  This becomes clear from the evaluation of the corresponding Weyl limits (see Remark \ref{dilated-Weyl}).  The general framework discussed in this article delivers a CLT for $\sum_{p \leq x} h(\alpha \theta_E(p))$ for any integer $\alpha \geq 1$.
     \end{enumerate}
       
       \item The sum over primes $\sum_{p \leq x} h(\alpha \theta_E(p))$ in the theorem includes the primes of bad reduction.  At these primes, the angles are not Frobenius angles, but their contribution to the sum is negligible on average (see Proposition \ref{bad-red-points}), since $h$ is bounded.

   \end{enumerate} 
\end{unmremark}

\subsection{A CLT for Kloosterman angles}\label{sec: Kloosterman angles main theorem}

The proof of Theorem \ref{Elliptic-curves-smooth} uses only two features of the family $\mathcal F_{A,B}$: averages over the family at a fixed prime are periodic modulo that prime, so that the Chinese remainder theorem makes the angles at distinct primes approximately independent; and at each prime, these averages are close to their limits, with an error of size $\O(|m|/\sqrt p)$ at frequency $m$, by Birch's formula and the Ramanujan--Deligne bound. Both features are shared by Kloosterman sums, where the role of the Ramanujan--Deligne bound is played by Katz’s bounds for symmetric powers of the Kloosterman sheaf. Katz \cite{Katz} showed that for a fixed prime $p$, the angles $\theta_{p,a}$, $a\in\mathbb F_p^\times$, become Sato–Tate distributed as $p\to\infty$; the corresponding horizontal statement for a fixed $a$ remains open. Averaging over $a$, we obtain the following central limit theorem.

\begin{thm}\label{Kloosterman-multiples-smooth}
  Let $g$, $h$, $\mu_{\alpha}$ and $V_h(\mu_{\alpha})$ be as defined in Theorem \ref{Elliptic-curves-smooth}. For a prime $p$, let $\F_p$ denote the finite field of order $p$, and let $\F_p^{\times}$ denote the set of nonzero elements of $\F_p$. For an integer $a$, the Kloosterman sum $K(p,a)$ is defined as follows:
\begin{equation*}
    K(p,a) = \sum_{x \in \F_p^{\times}} e\left(\frac{x + ax^{-1}}{p}\right), \, e(t) = e^{2\pi i t}.
\end{equation*}  
Choose the Kloosterman angle $\theta_{p,a} \in [0,1]$ such that $\displaystyle 2 \cos (\pi \theta_{p,a}) = \frac{K(p,a)}{\sqrt{p}}$.

For $N \geq 1$, let 
$\displaystyle \mathcal F_N = \{a:\,1 \leq a \leq N\}$. 
For a fixed positive integer $\alpha$, define $\mathcal S_x:\,\mathcal F_N \to \R$ by
$$\mathcal S_x(a) := \frac{1}{\sqrt{{\pi}(x)V_{h}(\mu_{\alpha}})}\left(\sum_{p \leq x }h(\alpha\theta_{p,a}) - {\pi}(x)\int_0^1 h(t)\mu_{\alpha}(t) dt\right).$$
If $N = N(x)$ satisfies  $\displaystyle \frac{\log N}{\log x} \to \infty \text{ as }x \to \infty $,
then for all positive integers $r$,
$$\lim_{x \to \infty}\frac{1}{N} \sum_{a=1}^N\mathcal S_x(a)^r  = \gamma_r.$$

\end{thm}

\begin{unmremark}\label{Kloosterman-CLT-remarks}\leavevmode

\begin{enumerate}

\item  The above theorem extends the horizontal central limit theorem of Xi \cite[Theorem 2]{Xi} in several directions:

\begin{enumerate}[label=(\alph*)]

\item \cite{Xi} considers smooth test functions supported in a closed subinterval of $(0,1)$. Theorem \ref{Kloosterman-multiples-smooth} applies to symmetric test functions, $h(1-t)=h(t)$, which need only be of limited smoothness: it suffices that $\widehat g(t)\ll(1+|t|)^{-\kappa}$ for some $\kappa>2$.  Combined with the extension to 2-periodic test functions (Theorem~\ref{2-periodic}), our method recovers \cite[Theorem 2]{Xi} and extends it to non-symmetric test functions of limited smoothness (see Corollary \ref{Nagoshi-type-Kloosterman}).  The growth condition $\log N/\log x\to\infty$ is the same as in \cite{Xi}.

\item Under the same growth conditions on $N$ as above, our framework also shows that as $a$ ranges from 1 to $N$, all moments of the sums $\pi(x)^{-1/2}\sum_{p\leq x}K(p,a)/\sqrt p$ converge to the  Gaussian moments, and so do the moments of the corresponding normalized sums for $U_j(K(p,a)/2\sqrt p)$, where $U_j(\cos (\pi t))$ denotes the $m$-th Chebyshev polynomial of the second kind.

\item  Theorem \ref{Kloosterman-multiples-smooth} gives a central limit theorem for the dilated sequence $(\alpha \theta_{p,a} \,\tmod 1)_p$ for each positive integer $\alpha$.  We note that the dilated angles are, up to sign, the Kloosterman angles for the field extensions $\mathbb{F}_{p^\alpha}$.

\end{enumerate}

\item Besides the Chinese remainder theorem, the proof of Theorem \ref{Kloosterman-multiples-smooth} uses only Katz’s bound for complete sums of the Chebyshev polynomial at a single prime (Proposition \ref{Katz}).  In \cite{Xi}, the sums over $a\leq N$ are completed with additive characters, which requires bounds for Kloosterman sums twisted by additive characters \cite[Lemmas 3 and 4]{Xi}.  In our framework, the growth condition on $N$ allows the incomplete sums over $a$ to be replaced by complete sums with a trivial error (Proposition \ref{appl of crt}).

\end{enumerate}
\end{unmremark}

\subsection{A CLT for Hecke angles}\label{sec: Hecke angles: main theorems}

In our third application, the family averages are no longer governed by periodicity modulo $p$, but by the Eichler–Selberg trace formula. This has two consequences. First, the trace formula estimates are uniform in both the weight and the level, which allows both to vary with $x$. Second, the convergence to the Weyl limits is now geometric in the frequency, while the error terms grow exponentially in it. As a result, the admissible growth of the family depends on the decay of the test function, and we obtain a graded family of results, one for each decay class.

\begin{thm}\label{Hecke-multiples-smooth}
 Let $g$, $h$, $\mu_{\alpha}$ and $V_h(\mu_{\alpha})$ be as defined in Theorem \ref{Elliptic-curves-smooth}, except that we only assume $\widehat{g}(t) \ll 1/(1+|t|)^{\kappa}$ for some $\kappa > 1$ (instead of $\kappa >2$).  For a positive integer $N$, and an even, positive integer $k$, let $\mathcal F_{N,k}$ denote the set of normalized Hecke newforms of weight $k$ with respect to the congruence subgroup $\Gamma_0(N)$ of the full modular group, where each $f \in \f_{N,k}$ has a Fourier expansion
$$f(z) = \sum_{n=1}^{\infty} n^{(k-1)/2}a_f(n)e^{2\pi i nz},\quad \quad a_f(1) = 1.$$
For a prime $p$, choose the Hecke angle $\thtp \in [0,1]$ such that
$\displaystyle 2 \cos (\pi \thtp) = a_f(p)$.

Let $N = N(x)$ and $k = k(x)$ \textup{(}with $k$ even\textup{)} be such that
$\omega(N) = \o(\sqrt{\pi(x)})$.  

For a fixed positive integer $\alpha$, define $\mathcal S_x\colon \mathcal F_{N,k} \to \R$ by
$$\mathcal S_x(f) := \frac{1}{\sqrt{\pi(x) V_h(\mu_{\alpha})}}\left(\sum_{p\leq x}h(\alpha \thtp)-\pi(x)\int_0^1h(t)\mu_{\alpha}(t)dt\right).$$

Suppose there exists $M = M(x)$ such that
\begin{align*}
\sum_{m > M} \abs{\widehat{h}(m)} = \o\left(\frac{1}{\sqrt{\pi(x)}}\right) \text{ and }\frac{\log (kN/8^{\omega(N)})}{M\log x} \to \infty.
\end{align*}

Then, for all positive integers $r$,
$$\lim_{x \to \infty}\frac{1}{|\mathcal F_{N,k}|} \sum_{f \in \mathcal F_{N,k}}\mathcal S_x(f)^r  = \gamma_r.$$  In particular, this conclusion holds for each of the following cases, in which the decay condition on $\widehat{g}$ is successively strengthened and the growth condition on $kN$ correspondingly relaxed.

    \begin{enumerate}[label=\normalfont\bfseries(\alph*), align=left, labelsep=0em] 
    \item\label{polynomial decay} $\displaystyle \widehat{g}(t) \ll \frac{1}{(1 +|t|)^{\kappa}}$ for some $\kappa > 1$ and $\displaystyle \log (kN/8^{\omega(N)}) \geq x^{1/(\kappa-1)}$.
    \item\label{Schwartz}$\displaystyle \widehat{g}(t) \ll \frac{1}{(1 +|t|)^B}$ for every $B > 0$ and $\displaystyle \log (kN/8^{\omega(N)}) \geq x^{\eta}$ for some $\eta > 0$.

    \item\label{exp-decay}
    $\displaystyle \widehat{g}(t) \ll e^{-\lambda|t|^{\gamma}}$ for some $\lambda,\gamma >0$ and $\displaystyle \frac{\log (kN/8^{\omega(N)})}{(\log x)^{1 + \frac{1}{\gamma}}} \to \infty.$
    
    \item\label{compact-support} $\widehat{g}$ is compactly supported and $\displaystyle \frac{\log (kN/8^{\omega(N)})}{\log x} \to \infty.$

\end{enumerate}

\end{thm}

\begin{unmremark}\leavevmode

Theorem \ref{Hecke-multiples-smooth}  extends previous central limit theorems in several directions.
\begin{enumerate}
    \item  In \cite{Nagoshi} and \cite[Theorem 1.6]{BPS}, the level $N$ is fixed, and only the weight $k$ tends to infinity. In Theorem \ref{Hecke-multiples-smooth}, the weight and the level may vary simultaneously, subject to the growth conditions on $kN$ in Theorem \ref{Hecke-multiples-smooth}, provided $\omega(N)=\o\big(\sqrt{\pi(x)}\big)$.  The condition on $\omega(N)$ ensures that the primes dividing $N$, at which the Hecke relations \eqref{prod of coeff prime power} are not available, contribute negligibly.  The conclusions of Theorem \ref{Hecke-multiples-smooth} extend to 2-periodic test functions (Theorem \ref{2-periodic}).  This allows us to recover Nagoshi's theorem for $\sum_{p\leq x}a_f(p)$ \cite{Nagoshi}, and its generalization to $\sum_{p\leq x}a_f(p^j)$ for Hecke eigenforms in \cite[Theorem 6.1]{murty-prabhu}, and to extend them to varying level (Corollary \ref{Nagoshi-revisited}).  In the level aspect, \cite[Section 4]{Cho-Kim}  considers Nagoshi's statistic for fixed $k\geq 2$ and $\log N/\log x\to\infty$, with $f$ ranging over a Hecke eigenbasis of $S_k(\Gamma_0(N))$, that is, a basis consisting of eigenforms of the Hecke operators $T_n$, $(n,N)=1$; our framework applies to this family as well (Remark \ref{natural families S(N,k)}).   
    \item Cases \ref{polynomial decay} and \ref{Schwartz} require only $\widehat{g}(t) \ll 1/(1+|t|)^{\kappa}$ for some $\kappa > 1$, and rapid decay of $\widehat{g}$, respectively, in contrast with the exponential-type decay $\displaystyle \widehat{g}(t) \ll e^{-\lambda|t|^{\gamma}}$ ($\lambda,\gamma >0$) assumed in \cite{BPS}. For fixed $N$ and $\alpha=1$, cases \ref{exp-decay} and \ref{compact-support} recover the corresponding results of \cite{BPS}.
    \item For $\alpha\geq 2$, the sequences $(\alpha\thtp)_{p\to\infty}$ have the uniform distribution as their limiting measure, and Theorem \ref{Hecke-multiples-smooth} gives a central limit theorem for them as well.
    
\end{enumerate}

\end{unmremark}

\subsection{General framework}\label{sec: General Framework}\leavevmode

The above theorems are built upon a general structure, whose prototype is the Erd\hungarianumlaut{o}s–Kac theorem on $\omega(n)$, the number of distinct prime factors of $n$: for any $a \in \R$,
$$\lim_{x \to \infty}\frac{1}{x}\#\left\{n \leq x:\,\frac{\omega(n) - \log \log x}{\sqrt{\log \log x}} \leq a\right\} = \frac{1}{\sqrt{2\pi}}\int_{-\infty}^{a} e^{-t^2/2} dt.$$
Here, $\omega(n)$ is the sum $\sum_{p\leq x}g_p(n)$, where $g_p(n) = 1$ if $p \mid n$ and 0 otherwise, and $n$ is chosen uniformly from $\{n \leq x\}$.  The random variables $g_p$ are not independent, but approximately independent: $\displaystyle\E[g_pg_q] = \E[g_p]\E[g_q] + \O(1/x)$ for distinct primes $p$ and $q$.  Recently, Murty, Murty and Pujahari \cite{MMP2} proved an all-purpose generalization of the Erd\hungarianumlaut{o}s--Kac theorem for $\omega(a_n)$ for integer-valued sequences $(a_n)_{n \geq 1}$, and presented interesting applications.  This setup rests on three principles: writing an arithmetic quantity as a sum of random variables on growing probability spaces, approximate rather than exact independence of these random variables, and the method of moments as a tool for proving a Gaussian law.  Similar principles underlie central limit theorems in other arithmetic settings.  In this article, we are concerned with one such setting, namely statistics of sequences in $[0,1]$ that vary over large families, with the angles attached to elliptic curves, modular forms and Kloosterman sums serving as our main examples.   

In our setting, these principles take the following form.  The probability space is a large family of arithmetic objects: elliptic curves with coefficients in a box, integers $a\leq N$, or newforms of a given level and weight.  To each member $X$ of the family is associated a sequence $(X_n)_{n\in\mathcal A}$ in $[0,1]$, typically indexed by the primes, such as the Frobenius, Kloosterman or Hecke angles.  For a test function $h$, the random variables are $h_n(X)=h(X_n)$, and the smooth statistics are the sums $\sum_{n\in\mathcal A\cap[1,x]}h_n(X)$, which have $A(x):=|\mathcal A\cap[1,x]|$ terms.  In our examples, $\mathcal A$ consists of primes, and $A(x)$ is essentially $\pi(x)$.  The family grows as $x\to\infty$ (see Setup \ref{spaces-rvs} for the precise setup).  Gaussian laws for such statistics, averaged over the family, have been established in a wide variety of contexts, such as  modular cusp forms \cite{Prabhu-Sinha}, \cite{BPS}, Hilbert cusp forms \cite{Das-Prabhu},  automorphic forms for $GL(N),\,N \geq 3$ \cite{LNW1}, holomorphic Siegel cusp forms \cite{KWY1}, and more generally cusp forms on split simple groups over $\Q$ \cite{KWY2},  families of elliptic curves \cite{Baier-Prabhu}, \cite{BPS}, Kloosterman sums \cite{Xi}, and eigenvalues of random regular graphs \cite{Johnson}.

The goal of this article is to formulate general principles behind such results.  More precisely, we lay down explicit sufficient conditions under which all moments of these sums, suitably normalized, converge to the Gaussian moments as $x \to \infty$.  The  random variables $\{h_n:\, n \in \mathcal A\}$  form a triangular array, and in general, they are not mutually independent, so that the classical central limit theorem does not apply.  Our main general theorem replaces independence by quantitative conditions on the underlying sequences: approximate independence of the $X_n$ at distinct $n$, measured through the functions $\cos(2\pi m X_n)$, together with convergence of the family averages to the Weyl limits at a suitable rate, tested against the truncations of $h$ (Theorem \ref{CLT general smooth}).  This parallels \cite{MS}, where a variant of the Erd\hungarianumlaut{o}s–Tur\'{a}n inequality expresses the discrepancy of a $\mu$-equidistributed sequence in terms of the rate at which its Weyl sums converge to the Weyl limits, with applications to the vertical distribution of Hecke eigenvalues.  Here, the same comparison between family averages and Weyl limits, combined with approximate independence at distinct indices, shows that the smooth statistics satisfy Gaussian laws.

In each of our applications, the arithmetic information is available only away from a small exceptional set: the curves with bad reduction at $p$, the integers $a$ divisible by $p$, and the primes dividing the level. Our framework therefore allows such sets to be discarded at each $n$, provided they are thin on average over the family (Theorem \ref{CLT general smooth variant}).

Our framework is formulated for general families of sequences, and we expect it to apply in other settings, such as those mentioned above, whenever the corresponding arithmetic input is available (see also Remark \ref{natural families S(N,k)}). In this article, we restrict ourselves to three families, which illustrate different mechanisms for verifying the hypotheses (see Section \ref{sec: smooth comparison}), and in each of which the framework leads to new results. We do not attempt to recover or extend the central limit theorems in the other settings listed above, and hope to return to this question in future work.

Central limit theorems of a different nature arise for the zeros of zeta functions of curves over a fixed finite field, as the genus tends to infinity (see, for example, \cite{FR} for hyperelliptic curves and \cite{BDFLS} for Artin--Schreier curves). Here, the statistics concern the $2g$ zeros of a single zeta function, whose angles behave like eigenvalues of random matrices, and the variance grows only logarithmically. In contrast, in our setting the angles at distinct primes are approximately independent, and the variance of the statistics grows linearly with the number of primes.

The smooth statistics mentioned above naturally lead to a study of the {\it sharp} statistics,$$N_{I,x}(X):=\#\{n\in\mathcal A\cap[1,x]:\ X_n\in I\}$$
for a fixed interval $I \subset (0,1)$.  The quantities $N_{I,x}-A(x)\mu(I)$ are the discrepancies of the sequences at the interval $I$.  As in the discrepancy bounds of \cite{MS}, we approximate the indicator function of the interval $I$ with the Beurling--Selberg polynomials (details in Section \ref{sec:linear-statistics}).  The comparison with smooth statistics requires a truncation degree $M \gg \sqrt{A(x)}$, and the resulting error terms impose much stronger growth conditions on the families than those required for smooth statistics.  Theorem \ref{CLT gen linear} gives general conditions under which these discrepancies satisfy Gaussian laws.  For Hecke angles, these conditions can be met, and we obtain a central limit theorem for $\#\{p\leq x:\ \{\alpha\theta_f(p)\}\in I\}$, which extends \cite{Prabhu-Sinha} to varying level and to dilated angles (Theorem \ref{Prabhu-Sinha type CLT, diff}).  For Frobenius and Kloosterman angles, the method would require local discrepancy bounds of essentially square-root strength, beyond the best known bounds $\O(p^{-1/4})$ (we explain this obstruction in Remark \ref{limitations of Theorem 52}\ref{obstruction-linear}).

\subsection{Organization of the article}

This article is organized as follows.  In Section \ref{sec: Smooth statistics}, after introducing Weyl limits for families of sequences in $[0,1]$, we develop a general framework that converts quantitative information about these limits into Gaussian laws for smooth statistics (Theorem \ref{CLT general smooth}).  The hypotheses of this theorem are three conditions: a truncation condition, approximate independence at distinct indices, and convergence of the family averages to the Weyl limits, tested against the truncated test function; the method of moments combines them into a Gaussian law, with convergence of all moments. Theorem \ref{CLT general smooth variant} refines this, allowing thin exceptional subsets of the family to be discarded, and Theorem \ref{2-periodic} extends both results to 2-periodic test functions. The verification of these hypotheses in the three applications follows a common pattern, summarized in Table \ref{tab: organization of sections 3,4 and 5}.  In Section \ref{sec: some useful tools}, we collect some elementary tools that are used repeatedly in the forthcoming sections: properties of the periodized test functions $h$, Chebyshev polynomials of the second kind, which express the functions $\cos(\pi m\theta)$ in terms of the arithmetic data $a_E(p^k)$, $a_f(p^k)$ and $U_k(K(p,a)/2\sqrt{p})$, and two elementary estimates for products and sums of products.

Section \ref{sec: Elliptic curves} proves Theorem \ref{Elliptic-curves-smooth}. The key arithmetic input is the bound at individual primes in Proposition \ref{S(n)-elliptic}\ref{Spm-bounds}, which follows from Birch’s formula and the Ramanujan--Deligne bound. Approximate independence follows from periodicity modulo primes and the Chinese remainder theorem (Proposition \ref{Condition 2 elliptic}), and the verification of the convergence to the Weyl limits, where the decay of the test function makes cancellation over primes unnecessary, is carried out in Proposition \ref{Weyl limits elliptic} and Section \ref{proof of theorem 1} (see also Remark \ref{Comments MRH free}).  Section \ref{sec: Kloosterman} proves Theorem \ref{Kloosterman-multiples-smooth} along the same lines, with Katz’s bounds for symmetric powers of the Kloosterman sheaf in place of Birch’s formula and the Ramanujan--Deligne bound. Section \ref{sec: modular forms} proves Theorem \ref{Hecke-multiples-smooth} using the Eichler--Selberg trace formula, which is uniform in the weight and the level. Each of these sections ends with the corresponding results for Nagoshi-type statistics, that is,  statistics for test functions defined by the Chebyshev polynomials of the second kind. In Section \ref{sec: smooth comparison}, we compare the three families of Sections \ref{sec: Elliptic curves}-\ref{sec: modular forms}, and explain how the arithmetic input determines the growth conditions and the admissible test functions in each family.

In Section \ref{sec:linear-statistics}, we turn to sharp statistics.  We approximate indicator functions of intervals by Beurling--Selberg trigonometric polynomials, and extend the moment computations of Section \ref{Central limit theorems for smooth statistics} to such test functions (Proposition \ref{rem:uniform}).  Using this, we prove a central limit theorem for sharp statistics (Theorem \ref{CLT gen linear} and Corollary \ref{CLT gen linear union}).  We apply it to Hecke angles (Theorem \ref{Prabhu-Sinha type CLT, diff}), and explain why the method does not extend to Frobenius and Kloosterman angles (Remark \ref{limitations of Theorem 52}).  As the reader navigates the general framework for smooth and sharp statistics along with applications, they will be aided by a list in Section \ref{sec: Notation} which recalls and collects the notation for the ingredients of the framework.

\subsection{Acknowledgments and AI declaration}

Some of the results of this paper are contained in the first-named author’s Master’s thesis.  We thank the Indian Institute of Science Education and Research (IISER) Pune for excellent working conditions.  Part of this work was carried out during the second-named author’s visit to the Lodha Mathematical Sciences Institute, Mumbai in December 2025, as part of its semester programme on Arithmetic Statistics.  She thanks the institute for its hospitality. We thank Stephan Baier for helpful discussions.

The first-named author was supported by a DST INSPIRE Scholarship for Higher Education during his studies at IISER Pune, when part of this work was carried out, and the third-named author is supported by the same scholarship.  The research of the second-named author is supported by ANRF Core Research Grant CRG/2023/001743.

{\bf AI Declaration.} The authors used Claude Opus 5.5 (Anthropic) to improve the exposition and to proofread the manuscript for mathematical correctness and flow. It also  assisted in drafting Proposition \ref{rem:uniform} and in streamlining the proof of Theorem \ref{Prabhu-Sinha type CLT, diff}. All mathematical content was verified by the authors, who take full responsibility for it.

\section{From Weyl limits to Gaussian laws}\label{sec: Smooth statistics}

In this section, we develop the general framework. In Section \ref{sec: ed-Weyl}, we recall Weyl limits of sequences in $[0,1]$, and in Section \ref{sec:2.2}, we extend them to families of sequences.  In Section \ref{Central limit theorems for smooth statistics}, we prove the central limit theorem for smooth statistics (Theorem \ref{CLT general smooth}), its refinement allowing thin exceptional sets (Theorem \ref{CLT general smooth variant}), and its extension to 2-periodic test functions (Theorem \ref{2-periodic}).  Section \ref{sec: table} summarizes how these results are applied in Sections \ref{sec: Elliptic curves}-\ref{sec: modular forms} (Table \ref{tab: organization of sections 3,4 and 5}), and Section \ref{sec: some useful tools} collects some tools used throughout the article.

\enlargethispage{-3\baselineskip}

\subsection{Equidistribution and Weyl limits}\label{sec: ed-Weyl}

\begin{defn}\label{mu-equi-seq}
A sequence $X = (X_n)_{n \geq 1} \subseteq [0,1]$ is said to be equidistributed with respect to a continuous probability
measure $\mu(t)\,dt$ (or $\mu$-equidistributed)  if for  any $f \in \mathcal C[0,1]$,
\begin{align}\label{sequence-ed}
    \lim_{N \to \infty} \frac{1}{N}\sum_{n=1}^N f(X_n) = \int_0^ 1 f(t) \mu(t) dt.
\end{align}

\end{defn}
By a well-known generalization of Weyl's equidistribution criterion (see \cite[Theorem 7.3, Chapter 1]{KN}), the sequence $X$ is $\mu$-equidistributed if and only if 
$$\lim_{N \to \infty} \frac{1}{N} \sum_{n=1}^N e(mX_n) = \int_0^1 e(mt)\mu(t) dt, \quad  e(t) = e^{2\pi i t}$$
for all integers $m$.  The quantity
$$c_m(X) = c_m := \lim_{N \to \infty} \frac{1}{N} \sum_{n=1}^N e(mX_n)$$
is referred to as the $m$-th Weyl limit of the sequence $X$.  Clearly, $c_0 = 1$. 

In this article, we will encounter $\mu$-equidistributed sequences for the following functions $\mu$.
\begin{enumerate}
\item $\mu(t) = 1$ (uniform distribution), for which  $c_m = 0$ for all $m \neq 0$.  
\item 
    $\mu(t) = 2\sin^2(\pi t)$ (the Sato--Tate measure) for which $c_{\pm 1} = -1/2$ and $c_m = 0$ for $|m| \geq 2$.
\item For a prime $p$, the $p$-adic Plancherel measure (see \eqref{p-adic plancherel measure}) for which $\displaystyle c_m = 
    \frac{1-p}{2p^{|m|}}$ for all $m \neq 0$.
\end{enumerate}
In all these examples, $c_m = c_{-m}$ for all integers $m$ so that $\mu(t) = \mu(1-t)$, and 
$$c_m = \lim_{N \to \infty} \frac{1}{N}\sum_{n=1}^N\cos(2\pi mX_n).$$
\begin{remark}\label{dilated-Weyl}
If a sequence $(X_n)_{n\geq 1} $ is equidistributed with respect to the Sato--Tate measure, then for any integer $\alpha \geq 2$, the reduced dilated sequence $(\alpha X_n\, \tmod 1)$ is uniformly distributed.  This is because
$$\lim_{N \to \infty}\frac{1}{N} \sum_{n=1}^N e(m\alpha X_n) = c_{m\alpha}(X) = 0\text{ for all }m \neq 0.$$
\end{remark}

\subsection{Smooth statistics for a family of sequences}\label{sec:2.2}

In this subsection, we set up the framework and define Weyl limits for families of sequences.  We have the following setup.
\begin{setup}\label{spaces-rvs}
    Let $\mathcal A$ be an infinite subset of $\N$, and define $A(x) := |\mathcal A \cap [1,x]|$.  For each integer $k \geq 1$, let $\f_k$ be a finite probability space with the uniform measure, such that $|\mathcal{F}_k| \to \infty$ as $k \to \infty$. To each $X \in \mathcal F_k$, we associate a sequence $(X_n)_{n\in \mathcal A}$ with $X_n \in [0,1]$ for all $n$.  For a function $\Phi:\,\mathcal F_k \to \C$, we write
    $$ \langle \Phi \rangle_k :=\frac{1}{|\mathcal{F}_k|} \sum_{X \in \mathcal F_k} \Phi(X).$$
Let $h\colon\R \to \R$ be an even, continuous, one-periodic function with the Fourier expansion
\begin{equation*}
\sum_{m \in \Z}\widehat{h}(m) e(mt), 
\end{equation*}
where $\widehat{h}(m) = \displaystyle\int_0^1 h(t) e(-mt) dt.$  We further assume that  
$\displaystyle\sum_{m \in \Z} |\widehat{h}(m)| < \infty$.  Thus, the Fourier series converges absolutely and uniformly in $t$, and 
\begin{equation}\label{h-FE}
h(t) = \sum_{m \in \Z}\widehat{h}(m) e(mt) = \sum_{m \in \Z}\widehat{h}(m)\cos (2\pi mt).
\end{equation}

For $n \in \mathcal A$, define the random variable $h_n:\,\mathcal F_k \to \R$ such that 
\begin{equation*}
    h_n(X) := h(X_n) \text{ for }X \in \mathcal F_k.
\end{equation*}
Then,
\begin{align*}
    \langle h_n \rangle_k  = \frac{1}{|\mathcal{F}_k|} \sum_{X \in \mathcal F_k} h(X_n) =\sum_{m \in \Z}\widehat{h}(m)\left(\frac{1}{|\mathcal{F}_k|} \sum_{X \in \mathcal F_k} \cos(2\pi mX_n)\right).
\end{align*}
Since the test function $h$ is even, only cosine sums appear in $\langle h_n \rangle_k$.

\end{setup}
  Inspired by the Weyl limits defined in the previous section, we define Weyl-type limits for the families $(\mathcal F_k)_{k \geq 1}$ as $k \to \infty$, which provide a framework to explore the distributional properties of the random variables $h_n$, as well as their sums.

\begin{defn}\label{defn-c-mn}
  Assume the notation of Setup \ref{spaces-rvs}.  For $n \in \mathcal A$, $m \in \Z$ and $k \in \N$, define 
  \begin{equation}\label{Weyl sum of X}
    c_{m,n,k} :=  \langle \cos(2\pi mX_n) \rangle_k = \frac{1}{|\mathcal{F}_k|}\sum_{X \in \mathcal F_k}\cos(2\pi mX_n).
\end{equation}
 
Let $k=k(x) \to\infty$ as $x\to\infty$.  For each $m \in \Z$, define
    \begin{align}\label{Weyl limit d_l}
        d_m := \lim_{x \to \infty}\frac{1}{A(x)}  \sum_{n \in \mathcal A \cap [1,x]}\langle \cos(2\pi mX_n)\rangle_{k(x)}
        &= \lim_{x \to \infty}\frac{1}{A(x)}\sum_{n \in \mathcal A \cap [1,x]}c_{m,n,k(x)},
        \end{align}
    provided the limit exists.  We refer to $d_m$ as the $m$-th Weyl limit of $(\mathcal F_k)_{k \geq 1}$ along $\mathcal A$ (the name is justified by the next lemma).

    \end{defn}
    \begin{lemma}\label{mu-pdf}
    Suppose that $d_m$ exists for every $m \in \Z$ and $\sum_{m\in \Z}|d_m|<\infty$.  Define
    \begin{align}\label{mu gen}
        \mu(t) := \sum_{m \in \Z} d_m e(-mt).
    \end{align}
    Then, $\mu$ is a continuous, real-valued, nonnegative function on $[0,1]$ with $\mu(1-t)= \mu(t)$, $$\int_0^1 e(mt)\mu(t) dt = d_m \text{ for every }m \in \Z.$$  In particular, $\displaystyle\int_0^1 \mu(t) dt=d_0 =1$.  Thus, $\mu(t)\,dt$ is a symmetric probability measure on $[0,1]$ whose $m$-th Fourier coefficient is the Weyl limit $d_m$.
        
    \end{lemma}
    \begin{proof}

Since $\sum_{m\in \Z}|d_m|<\infty$, the series $\sum_{m \in \Z} d_m e(-mt)$ converges uniformly, and therefore $\mu$ is continuous.  Since $d_m \in \R$ and $d_m = d_{-m}$, $\mu(t)$ is real-valued with $\mu(1-t) = \mu(t)$.  Integrating term by term, which is justified by uniform convergence, and using the orthogonality relation
\begin{equation}\label{orth-exp}
\int_0^1 e((m - \ell)t) dt = \delta_{m,\ell}
\end{equation}
gives  
$\displaystyle \int_0^1 e(mt)\mu(t) dt = d_m$, 
and in particular,
$\displaystyle \int_0^1 \mu(t) dt = d_0 = 1$. 

To show nonnegativity, let $\displaystyle F_N(y) = \sum_{|m| \leq N}\left(1 - \frac{|m|}{N+1}\right)e(my) \geq 0$ be the Fej\'{e}r kernel.  Then,
$$\frac{1}{2}(F_N(y-t) + F_N(-y-t)) = \sum_{|m| \leq N}\left(1 - \frac{|m|}{N+1}\right)\cos(2\pi my)e(-mt) \geq 0.$$
Taking $y = X_n$, averaging over $X \in \mathcal F_k$ and $n\in \mathcal A \cap [1,x]$, and letting $x \to \infty$, we obtain, by \eqref{Weyl limit d_l},
$$\sigma_{N}(t) :=\sum_{m = -N}^{N }\left(1 - \frac{|m|}{N+1}\right) d_m e(-mt) \geq 0.$$
Since $\sum_m|d_m|<\infty$ and $1-\frac{|m|}{N+1}\to1$ for each $m$, dominated convergence gives $\sigma_N(t)\to\mu(t)$ as $N\to\infty$; hence $\mu(t) \geq 0$.
\end{proof}

\begin{remark}\label{discrepancy-polynomial}
    We make some observations and collect notation which will be repeatedly used in Lemma \ref{properties of powers of Z_M},  Theorem \ref{CLT general smooth} and Theorem \ref{CLT general smooth variant}.  Let $h(t)$ be as defined in Setup \ref{spaces-rvs} and $\mu(t)$ be as defined in \eqref{mu gen}.
    \begin{enumerate}[label=(\roman*), align=left, labelsep=0em]

    \item \label{Variance: Z_M gen}   Define the function
    \begin{equation}\label{Z_M poly gen}
         Z(t) :=  h(t) - \int_0^1 h(y) \mu(y) dy.   
    \end{equation} 
    We immediately see that $$\int_0^1 Z(t) \mu(t) dt = 0.$$
    Since $\int_0^1 \mu(t) dt = 1$, a direct computation also shows that
\begin{equation*}
         \int_0^1 Z(t)^2 \mu(t) dt =\int_0^1 h(t)^2 \mu(t)dt - \left(\int_0^1 h(t) \, \mu(t) dt\right)^2.
    \end{equation*}

\item  Henceforth, we denote
    \begin{equation}\label{Vh-defn}
        V_{h}(\mu) := \int_0^1 Z(t)^2 \mu(t) dt.
    \end{equation}
Note that $V_{h}(\mu) \geq 0$, and $V_{h}(\mu) > 0$ unless $h$ is constant on the support of $\mu$.  In the next theorem, $A(x)V_{h}(\mu)$ will appear as the  variance of $\sum_{n \in \mathcal A \cap [1,x]} h_n$ (asymptotically).

        \item \label{Z_M poly gen Fourier series}

       Since $$\int_0^1 h(y) \mu(y) dy = \sum_{m \in \Z}\widehat{h}(m)d_m,$$
       $Z(t)$ has the Fourier expansion
\begin{equation*}
        Z(t) = \sum_{m \in \Z} \widehat{Z}(m) e(mt),
\end{equation*}
    where
\begin{equation*}
         \widehat{Z}(m)=\begin{cases}
                            \widehat{h}(m) - \sum_{l \in \Z}d_l \widehat{h}(l) &\text{ if }m = 0\\
                            \widehat{h}(m) &\text{ if }|m| \geq 1. 
                        \end{cases}
\end{equation*}

\item\label{hMZM} For $M  \geq 1$, denote
    $$h_M(t) := \sum_{|m| \leq M} \widehat{h}(m) e(mt) \hspace{0.5cm}\text{and} \hspace{0.5cm} Z_M(t) := h_M(t) - \int_0^1h_M(y) \mu(y) dy. $$

    \end{enumerate}
\end{remark}

\subsection{Central limit theorems for smooth statistics}\label{Central limit theorems for smooth statistics}
In this subsection, we prove a central limit theorem for the sums $ \sum_{n \in \mathcal A \cap [1,x]} h_n$, where the random variables $h_n:\,\mathcal F_k \to \R$ are as defined in Setup \ref{spaces-rvs}.  Theorems \ref{Elliptic-curves-smooth}-\ref{Hecke-multiples-smooth} follow by verifying the hypotheses of this general theorem for each family (Sections \ref{sec: Elliptic curves}-\ref{sec: modular forms}).  The hypotheses consist of three conditions:  a truncation condition \ref{mean square convergence}, which allows us to replace $h$ by a trigonometric polynomial $h_M$,    approximate independence of $X_n$ at distinct $n$ \ref{condn 1, smooth}, and the convergence of the family averages  $c_{m,n,k}$ to the Weyl limits $d_m$ tested against the powers of $h_M$ \ref{condn 2, smooth}.  The method of moments then combines these into a Gaussian law.

We explain this framework in two phases.  Theorem \ref{CLT general smooth} contains the core argument, connecting the family averages $c_{m,n,k}$ and the Weyl limits $d_m$ with the moments of the sums under consideration.  Theorem \ref{CLT general smooth variant} presents a further technical refinement, which allows us to ignore thin subsets of $\mathcal F_k$ while executing the moment calculations.  In all three families of our interest, this refinement is needed because of the form in which the arithmetic information about these families is available.  However, we prove Theorem \ref{CLT general smooth} first, since it isolates the primary mechanism of the CLT without the additional notation of the thin subsets.  Theorem \ref{2-periodic} then extends both results to 2-periodic test functions.

At the heart of our study is the centered sum
$$\sum_{n \in \mathcal A \cap [1,x]} Z(X_n) = \sum_{n \in \mathcal A \cap [1,x]}h_n(X) - A(x)\int_0^1 h(y) \mu(y) dy,$$
where $Z$ is the centered test function of Remark \ref{discrepancy-polynomial}.  The proof of the CLT consists in computing the moments  $ \left \langle \left(\sum_nZ(X_n)\right)^r\right \rangle_k$.  In these computations, the families $\mathcal F_k$ enter through the averages $c_{m,n,k}$ and the limits $d_m$, and the test function enters through the Fourier coefficients of the powers of its truncations.
We record the relevant properties of the latter in a lemma below.

\begin{lemma}\label{properties of powers of Z_M}
Let $h$ and $\mu$ be as defined in Setup \ref{spaces-rvs} and Lemma \ref{mu-pdf}.  Assume the notation of Remark \ref{discrepancy-polynomial}.    
    For each integer $r \geq 1$,
\begin{enumerate}[label=\normalfont(\alph*), align=left, labelsep=0em] 
\item\label{cosine-expansion} $Z_M(t)^r$ is a trigonometric polynomial of degree at most $Mr$, with real and even Fourier coefficients $\widehat{Z_M^{r}}(m)$. Hence,
\begin{equation*}
Z_M(t)^r = \sum_{|m| \leq Mr}\widehat{Z_M^{r}}(m) \cos (2\pi mt).
\end{equation*}
\item\label{ZM-power-F} 
$
\widehat{Z_M^{r}}(m) \ll 1
$, with the implied constant depending only on $h$ and $r$.
\item\label{DMr-all}
$\displaystyle D(M,r) :=\int_0^1 Z_M(t)^r\mu(t) dt= \sum_{|m| \leq Mr} \widehat{Z_M^{r}}(m) d_m $.
In particular,
\begin{equation}\label{DM1}
D(M,1)  = 0,
\end{equation}
and
\begin{equation}\label{DMr}
|D(M,r)| \leq \left(2\sum_{m \in \Z}\abs{\widehat{h}(m)}\right)^r.
\end{equation}
\item\label{DM2}  If $M = M(x)$ satisfies $\displaystyle \sum_{|m| > M(x)}\abs{\widehat{h}(m)} \to 0$, then
\begin{equation*}
D(M,2) = \int_0^1 Z_M(t)^2\mu(t) dt = V_h(\mu) + \o(1) \text{ as }x \to \infty.
\end{equation*}

\end{enumerate}

\end{lemma}

\begin{proof}
The Fourier coefficients 
$$\widehat{Z_M^{r}}(m) = \int_0^1 Z_M(t)^{r}e(-mt) \,dt$$
are real and even, since the coefficients $\widehat{h}(m)$ are real and even (as $h$ is real and even). 
This proves \ref{cosine-expansion}.

Observe that $\displaystyle |Z_M(t)| \leq 2\sum_{m \in \Z}|\widehat{h}(m)|$ and 
$ \displaystyle |Z_M(t)|^{r} \leq \left(2\sum_{m \in \Z}|\widehat{h}(m)|\right)^{r} \ll 1$ uniformly in $M$ and $t$.  Therefore, $\widehat{Z_M^{r}}(m) \ll 1$, with the implied constant depending only on $h$ and $r$.  This proves \ref{ZM-power-F}.

 \ref{DMr-all} follows by the orthogonality of the set $\{e(mt):\,m \in \Z\}$.  This immediately implies that
$$D(M,1) = \int_0^1 Z_M(t) \mu(t) dt  = \int_0^1 h_M(t) \mu(t) dt - \int_0^1h_M(y) \mu(y) \left(\int_0^1 \mu(t) dt\right)dy = 0,$$
proving \eqref{DM1}.  Equation \eqref{DMr} follows by the bound for $\abs{Z_M(t)}^r$ mentioned in the proof of \ref{ZM-power-F}.

 To prove \ref{DM2}, we observe that $$ \sup_{t \in [0,1]}\abs{Z(t) - Z_M(t)} \leq 2\sum_{|m| > M}\abs{\widehat{h}(m)} = \o(1).$$ Since
\begin{align*}
        D(M,2)&=\int_0^1 (Z_M(t))^2\mu(t)\,dt\\
        &=\int_0^1 (Z_M(t)-Z(t))^2\mu(t)\,dt+\int_0^1 (Z(t))^2\mu(t)\,dt+2\int_0^1 (Z_M(t)-Z(t))Z(t)\mu(t)\,dt,
    \end{align*}
    we derive \ref{DM2}.
\end{proof}
\begin{thm}\label{CLT general smooth}
    
    Let $(\mathcal F_k)_{k \geq 1}$ and $h$ be as in Setup \ref{spaces-rvs}.  Let $\mathcal A$ be an infinite subset of $\N$, and $A(x) = \lvert \mathcal A \cap [1,x]\rvert$.  Let $k = k(x)\to\infty$ be such that the Weyl limits $d_m$ of Definition \ref{defn-c-mn} exist for each $m \in \Z$, and satisfy $\sum_{m \in \Z} |d_m| < \infty$.  Let $\displaystyle\mu(t) = \sum_{m \in \Z} d_m e(-mt)$ be the associated probability density (Lemma \ref{mu-pdf}).   
    Put
    \[Z(t) := h(t) - \int_0^1 h(y)\mu(y)\, dy,\quad V_h(\mu) := \int_0^1 Z(t)^2 \mu(t)\, dt,\]
    and assume that $V_h(\mu) > 0$.  For $M  \geq 1$, denote
    \[h_M(t) := \sum_{|m| \leq M} \widehat{h}(m) e(mt), \quad Z_M(t) := h_M(t) - \int_0^1h_M(y) \mu(y)\, dy.\]
Let $M = M(x) \geq 1$ be such that 
\begin{align}\label{ZM to Z uniform}
 \sum_{|m| > M} \abs{\widehat{h}(m)} \to 0\text{ as }x \to \infty.   
\end{align}

For an integer $u \geq 1$ and positive integers $r_1,r_2,\dots,r_u$, put $r = r_1 + r_2 + \dots +r_u$.  Let $\sum_{\underline{n}}$ denote the sum over all $u$-tuples $\underline{n} = (n_1,n_2,\dots,n_u)$ of distinct elements of $\mathcal A \cap [1,x]$, and  let $\sum_{\underline{m}}$ denote the sum over all $u$-tuples $\underline{m} = (m_1,m_2,\dots,m_u) \in \Z^u$ with $|m_i| \leq Mr_i$ for each $i$.

    Suppose that the following asymptotic estimates hold as $x \to \infty$. 
    \begin{enumerate}[label=\normalfont\bfseries (C\arabic*)] 
        \item \label{mean square convergence}
        {\bf Truncation.} For every integer $l \geq 1$, 
        $$\displaystyle \Bigg\langle \abs{\sum_{n \in \mathcal{A} \cap [1,x]}(Z(X_n) - Z_M(X_n) ) }^l\Bigg\rangle_k= \o\left(A(x)^{l/2}\right).$$

        \item \label{condn 1, smooth} {\bf Approximate independence.} For all positive integers $u, r_1,r_2,\dots, r_u$, 
           \begin{align*} 
\sum_{\underline{n}}\sum_{\underline{m}} \left\lvert \prod_{i=1}^u \left \langle\cos (2\pi m_i X_{n_i}) \right \rangle_k -\left \langle \prod_{i=1}^u \cos (2\pi m_i X_{n_i}) \right \rangle_k\right \rvert
= \o(A(x)^{r/2}).
            \end{align*}
 \textup{(}Recall that $\displaystyle \left \langle\cos (2\pi m_i X_{n_i})\right \rangle_k = c_{m_i,n_i,k}.$\textup{)}

        \item \label{condn 2, smooth} {\bf Convergence to the Weyl limits.}
For all positive integers $u,r_1,r_2,\dots,r_u$,
        \begin{align*}
\sum_{\underline{n}}\prod_{i=1}^u\sum_{|m_i| \leq Mr_i} \abs{\widehat{Z_M^{r_i}}(m_i)(c_{m_i,n_i,k}-d_{m_i})}=\o\left(A(x)^{r/2}\right),\quad \widehat{Z_M^{r}}(m) = \int_0^1 Z_M(t)^r e(-mt) dt.
        \end{align*}
        
    \end{enumerate}
    
For $X \in \mathcal F_{k(x)}$, define
$$\mathcal S_x(X) := \frac{1}{\sqrt{A(x) V_h(\mu)}}\left(\sum_{n \in \mathcal A \cap [1,x]} h_n(X) - A(x)\int_0^1 h(t) \, \mu(t) \, dt\right).$$
Then, all the moments of $\mathcal S_x$ converge to the Gaussian moments as $x \to \infty$.  That is, for all positive integers $r$,
$$\lim_{x \to \infty}\left \langle \mathcal S_x(X)^r \right \rangle_{k(x)} = \gamma_r :=\frac{1}{\sqrt{2\pi}}\int_{-\infty}^{\infty}  t^re^{-t^2/2}dt.$$
In particular, since the Gaussian distribution is determined by its moments, $\mathcal S_x$ converges in distribution to the standard Gaussian.

\end{thm}

\begin{proof} The proof of this theorem is structured into a few key steps: truncation,  moment expansion of the truncated polynomial in terms of the family averages $c_{m,n,k}$ and Fourier coefficients of the higher powers of the test function, and replacing $c_{m,n,k}$ by the Weyl limits $d_m$ in order to isolate the main terms in the moment expansion, which match the Gaussian moments.

\noindent
\hypertarget{step1}{\textbf{Step 1 (Reduction to trigonometric polynomials).}}

We have
\begin{equation}\label{T_x -A(x).Avg '} 
   \sum_{n \in \mathcal A \cap [1,x]} h_n(X)-A(x)\int_0^1 h(t) \, \mu(t)\,dt=\sum_{n \in \mathcal{A} \cap [1,x]}Z(X_n).
\end{equation}

By \eqref{T_x -A(x).Avg '}, $\mathcal S_x = S_x + T_x,$ where
$$S_x := \frac{ \sum_{n \in \mathcal A \cap [1,x]}Z_M(X_n)}{\sqrt{A(x) V_h(\mu)}}\text{ and }
 T_x := \frac{ \sum_{n \in \mathcal A \cap [1,x]}(Z(X_n) - Z_M(X_n))}{\sqrt{A(x) V_h(\mu)}}$$
By Condition \ref{mean square convergence}, for each $l \geq 1$,
\begin{align*}
\|T_x\|_l := \left(\left \langle |T_x|^l \right \rangle \right)^{1/l} = \frac{\o(A(x)^{1/2})}{(A(x) V_h(\mu))^{1/2}}   \to 0.
\end{align*}

In Steps \hyperlink{step2}{2}, \hyperlink{step3}{3} and \hyperlink{step4}{4} below, we show that for every $r \geq 1$,
\begin{equation}\label{moment-identity}
\lim_{x \to \infty} \left \langle (S_x)^r \right \rangle_k = \begin{cases}
                        \frac{r!}{2^{r/2} (r/2)!} &\text{ if }r \text{ is even}\\
                        0 &\text{ if }r \text{ is odd}
                    \end{cases} = \frac{1}{\sqrt{2\pi}}\int_{-\infty}^{\infty}t^r e^{-t^2/2} dt = \gamma_r.
                    \end{equation}
In Step \hyperlink{step5}{5}, we deduce the same for $\mathcal S_x$.

\par\smallskip\noindent
\hypertarget{step2}{\textbf{Step 2 (Moment expansion and decoupling).}}

We have 
\begin{equation}\label{multinomial gen'} 
(A(x) V_h(\mu))^{r/2}S_x^r = \left(\sum_{n \in \mathcal A \cap [1,x]}Z_M(X_n)\right)^r\\
=\sum_{u=1}^r \sum_{\underline{{\bf r}}}\frac{r!}{r_1!r_2!\dots r_u!u!}\sum_{\underline{n}}\prod_{i=1}^u Z_M(X_{n_i})^{r_i},
\end{equation}
where 
 $\sum_{\underline{{\bf r}}}$ denotes the sum over all $u$-tuples $\underline{\bf r} = (r_1,r_2,\dots,r_u)$ of positive integers such that $r_1 + r_2 +\dots +r_u = r.$  

 By Lemma \ref{properties of powers of Z_M}\ref{cosine-expansion}, we have
\begin{equation}\label{higher-moments-MT-ET-gen}
    \begin{split}
   &\left\langle \left(\sum_{n \in \mathcal A \cap [1,x]}Z_M(X_n)\right)^r\right \rangle\\
   &=\sum_{u=1}^r \sum_{\underline{{\bf r}}}\frac{r!}{r_1!r_2!\dots r_u!u!}\sum_{\underline{n}}\left \langle\prod_{i=1}^u \sum_{|m_i| \leq Mr_i} \widehat{Z_M^{r_i}}(m_i) \cos (2\pi m_i X_{n_i})\right\rangle\\
   &=\sum_{u=1}^r \sum_{\underline{{\bf r}}}\frac{r!}{r_1!r_2!\dots r_u!u!}\sum_{\underline{n}}\sum_{\underline{m}} \prod_{i=1}^u \widehat{Z_M^{r_i}}(m_i)\left\langle \prod_{i=1}^u\cos (2\pi m_i X_{n_i}) \right\rangle\\
   &= \sum_{u=1}^r \sum_{\underline{{\bf r}}}\frac{r!}{r_1!r_2!\dots r_u!u!}\sum_{\underline{n}} \prod_{i=1}^u \left(\sum_{|m_i| \leq Mr_i} \widehat{Z_M^{r_i}}(m_i)c_{m_i,n_i,k}\right)+ \o\left(A(x)^{r/2}\right),
    \end{split}
\end{equation}
where the last equality follows by Lemma \ref{properties of powers of Z_M}\ref{ZM-power-F} and Condition \ref{condn 1, smooth}.

\par\smallskip\noindent
\hypertarget{step3}{\textbf{Step 3 (Replacing averages over $\mathcal F_k$ by the Weyl limits $d_m$).}}

Write
\begin{align*}
    \sum_{\underline{n}}\prod_{i=1}^u \left(\sum_{|m_i| \leq Mr_i} \widehat{Z_M^{r_i}}(m_i)c_{m_i,n_i,k}\right)
\end{align*}
as 
\begin{align*}
    \sum_{\underline{n}}\prod_{i=1}^u \left( E_k(M,n_i,r_i)+ D(M,r_i)\right),
\end{align*}
where
\begin{equation}\label{r-th moment error term of Z_M at index n}
    E_k(M,n_i,r_i):=\sum_{|m_i| \leq Mr_i} \widehat{Z_M^{r_i}}(m_i)(c_{m_i,n_i,k}-d_{m_i}) 
\end{equation}
and $D(M,r_i)$ is as in Lemma \ref{properties of powers of Z_M}\ref{DMr-all}.

Let $\underline{{\bf r}} = (r_1,r_2,\dots,r_u)$ be fixed.
For a tuple $\underline{n}$, 
\begin{equation}\label{dcd}
\begin{split} &\prod_{i=1}^u \left(E_k(M,n_i,r_i)+ D(M,r_i)\right)\\
={}&\prod_{i=1}^u D(M,r_i)+\O\Bigg(\sum_{\substack{I \subseteq \{1,2,\dots,u\} \\I \neq \emptyset}}\prod_{i \in I}\left |E_k(M,n_{i},r_{i})\right| \prod_{i \notin I} \left|D(M,r_i)\right|\Bigg).
\end{split}
\end{equation}
    Since $D(M,1) = 0$ (see \eqref{DM1}), a set $I$ contributes to the $\O$-term in equation \eqref{dcd} only if $I$ contains every $i$ with $r_i = 1$.
Fix such a set $I$.  Bounding $D(M,r_i) \ll 1$ (see \eqref{DMr}) for $i \notin I$ and enlarging the range of $n_i$ for $i \notin I$ to $\mathcal A \cap [1,x]$ (which contributes at most $A(x)^{u - |I|}$), Condition \ref{condn 2, smooth} applied to the tuple $(r_i)_{i \in I}$ gives
\begin{equation}\label{Cond 3 error term}
\begin{split}
&\sum_{\underline{n}}\prod_{i \in I}\left|E_k(M,n_{i},r_{i})\right| \prod_{i \notin I} \left|D(M,r_i)\right|\\    ={}& \o\left(A(x)^{\sum_{i \in I} r_i/2}A(x)^{u-|I|}\right).\\
\end{split}
\end{equation}
Observing that  $\displaystyle\sum_{i \in I}r_i + 2(u-|I|) \leq r$, and summing \eqref{dcd} over $\underline{n}$ yields
$$\sum_{\underline{n}}\prod_{i=1}^u \left(E_k(M,n_i,r_i)+ D(M,r_i)\right) = \sum_{\underline{n}}\prod_{i=1}^u D(M,r_i) + \o(A(x)^{r/2}).$$

\par\smallskip\noindent
\hypertarget{step4}{\textbf{Step 4 (Evaluating the main terms).}}

By Steps \hyperlink{step2}{2} and \hyperlink{step3}{3}, for a fixed positive integer $r$, we have
\begin{equation}\label{higher-moments-MT-ET-gen-1}
    \left\langle \left(\sum_{n \in \mathcal A \cap [1,x]}Z_M(X_n)\right)^r\right \rangle_k
   =\sum_{u=1}^r \sum_{\underline{{\bf r}}}\frac{r!}{r_1!r_2!\dots r_u!u!}\sum_{\underline{n}}\prod_{i=1}^u D(M,r_i) + \o(A(x)^{r/2}).
   \end{equation}

We now take different cases for the choice of $(r_1,r_2,\dots,r_u)$, and evaluate the main term on the right hand side of \eqref{higher-moments-MT-ET-gen-1}.

\textbf{Case 1:} $r$ is even and $\underline{{\bf r}}=(2,2,\dots,2)$. Hence, $r=2u$.

In this case, 
\begin{align*}
    \sum_{\underline{n}}\prod_{i=1}^u D(M,2) = \sum_{\underline{n}}\left(\int_0^1 Z_M(t)^2\mu(t)\,dt\right)^u.
\end{align*}
The number of distinct $u$-tuples of elements in $\mathcal A \cap [1,x]$ is 
$$A(x)(A(x) - 1) \dots (A(x) - (u-1)) = A(x)^u + \O(A(x)^{u-1}).$$
Thus, by Lemma \ref{properties of powers of Z_M}\ref{DM2}, we have 
\begin{align*}
    \sum_{\underline{n}} D(M,2)^{u}= \left(A(x)^{u} + \O\left(A(x)^{u-1}\right)\right)\left(V_h(\mu)+\o(1)\right)^{u} = A(x)^{r/2}V_h(\mu)^{r/2} + \o(A(x)^{r/2}),
\end{align*}
and the contribution of this case to \eqref{higher-moments-MT-ET-gen-1} is
$$\frac{r!}{(r/2)!2^{r/2}}A(x)^{r/2}V_h(\mu)^{r/2} + \o(A(x)^{r/2}).$$

\textbf{Case 2:}  $\underline{{\bf r}}$ is a $u$-tuple such that $r_i=1$ for some $1\leq i\leq u$.  By equation \eqref{DM1},
$$\prod_{i=1}^u D(M,r_i) = 0.$$

\textbf{Case 3:} $\underline{{\bf r}}$ is a $u$-tuple such that $r_i\geq 2$ for all $1\leq i\leq u$ and $r_j\geq 3$ for some $1\leq j\leq u$.

In this case, $u \leq (r-1)/2$.  We use the bound in \eqref{DMr}  to obtain 
\begin{align*}
    \sum_{\underline{n}}\prod_{i=1}^u D(M,r_i)\ll \sum_{\underline{n}} 1 \ll A(x)^{(r-1)/2} =\o(A(x)^{r/2}).
\end{align*}
For even $r$, exactly one composition ${\bf r} = (2,2,\dots,2)$ with $u = r/2$ falls in Case 1; for $r$ odd, none does.  All other compositions fall in Case 2 or Case 3.
Dividing \eqref{higher-moments-MT-ET-gen-1} by $(A(x)V_h(\mu))^{r/2}$,
we get
$$\lim_{x \to \infty} \left\langle \left(\frac{\sum_{n \in \mathcal A \cap [1,x]}Z_M(X_n)}{\sqrt{A(x) V_h(\mu)}}\right)^r\right \rangle = \begin{cases}
                        \frac{r!}{2^{r/2} (r/2)!} &\text{ if }r \text{ is even},\\
                        0 &\text{ if }r \text{ is odd}
                    \end{cases}.$$
This proves the moment identity \eqref{moment-identity} required in Step \hyperlink{step1}{1}.

\par\smallskip\noindent
\hypertarget{step5}{\textbf{Step 5 (From $S_x$ to $\mathcal S_x$).}}
Since $h$ is real-valued, so are $Z_M$ and $S_x$. For any even integer $v \geq 2$, $|S_x|^v = S_x^v$, and so, $\displaystyle \left \langle \,\abs{S_x}^v \right \rangle_k = \left \langle S_x^v \right \rangle_k \to \gamma_v$ as $x \to \infty$ by \eqref{moment-identity}. By the binomial theorem and the Cauchy--Schwarz inequality, for any $r \geq 1$,
\begin{align*}
    \abs{\langle \mathcal S_x^r \rangle_k - \langle  S_x^r \rangle_k} \leq \sum_{j=1}^r \binom{r}{j} \big\langle|S_x|^{2(r-j)}\big\rangle_k^{1/2}\,\|T_x\|_{2j}^{\,j}\longrightarrow0,
\end{align*}
since $\|T_x\|_{2j} \to 0$ for each $j$ by Step \hyperlink{step1}{1}, and $\left \langle|S_x|^{2(r-j)}\right \rangle = \left \langle S_x^{2(r-j)}\right \rangle,$ which converges to $\gamma_{2r-2j}$ (by \eqref{moment-identity}) and is therefore bounded for each $j$.  This gives us $\langle \mathcal S_x ^r\rangle_k \to \gamma_r $ for every $r \geq 1$.  This proves the theorem.  In particular, since the Gaussian distribution is determined by its moments, the method of moments (\cite[Theorem 30.2]{Billingsley}) shows that $\mathcal S_x$ converges in distribution to the standard Gaussian.   
\end{proof}

\begin{remark}\label{remark after general CLT}\leavevmode

\begin{enumerate}[label=(\alph*)]
\item Condition \eqref{ZM to Z uniform} holds whenever $M(x) \to \infty$, since $\displaystyle \sum_m|\widehat h(m)|<\infty$, and $\displaystyle \sum_{|m|>M}|\widehat h(m)| \to 0$.  It also holds with $Z_M = Z$ if $h$ is a trigonometric polynomial of degree at most $M$.  In the latter case, $M$ may be taken to be fixed.

\item Condition \ref{mean square convergence} holds in particular if
$$\sup _{X \in \mathcal F_k} \abs{\sum_{n \in \mathcal A \cap [1,x]}(Z - Z_M)(X_n)} = \o\left(\sqrt{A(x)}\right),$$
which is how it is verified in all our applications.
\item  By Lemma \ref{properties of powers of Z_M}, \ref{condn 2, smooth} implies
$$\sum_{\underline{n}}\prod_{i=1}^u \left|\left \langle Z_M(X_{n_i})^{r_i}\right \rangle_k - \int_0^1Z_M(t)^{r_i}\mu(t) dt \right| = \o(A(x)^{r/2}).$$
That is, the distribution of $X_n$ over $\mathcal F_k$ is close to $\mu$ when tested against the powers of $Z_M$, with the discrepancies small {\it on average over $n$.}
        
 \item Since $\abs{\widehat{Z_M^r}(m)} \ll_r 1$ uniformly in $M$ and $m$, \ref{condn 2, smooth} follows from the ``unweighted'' condition    \begin{align*}
\sum_{\underline{n}}\prod_{i=1}^u\sum_{|m_i| \leq Mr_i} \abs{(c_{m_i,n_i,k}-d_{m_i})}=\o\left(A(x)^{r/2}\right),
\end{align*}
which involves only the family averages and the Weyl limits.  This condition is often convenient to verify; we use it for Hecke angles (see Proposition \ref{Towards GC 3}).  For Frobenius and Kloosterman angles, however, verifying \ref{condn 2, smooth} directly allows a much larger truncation level $M$, which is what leads to convergence of moments of $\mathcal S_x$ to the Gaussian moments.

\item The idea of working with the powers of the approximating trigonometric polynomials directly, instead of expanding them, together with the uniform boundedness of their Fourier coefficients, appears in an unpublished note of Sun, Wen and Zhang \cite{SWZ}, who used it to give a simpler proof of the main theorem of \cite{Prabhu-Sinha}. This idea underlies the weighted condition \ref{condn 2, smooth}, and later, \ref{Cond GC 3}, in the general setting of this article.  In Section \ref{sec:linear-statistics}, Proposition \ref{rem:uniform} and Theorem \ref{Prabhu-Sinha type CLT, diff} adapt it further to test functions that depend on $x$.

\item Condition \ref{mean square convergence} for a single $l \geq 1$, together with \ref{condn 1, smooth} and \ref{condn 2, smooth}, still yields convergence of $\mathcal S_x$ in distribution to the standard Gaussian.  Indeed, Steps 1–4 of the proof use only \ref{condn 1, smooth} and \ref{condn 2, smooth}, and show that the moments of $S_x$ converge to the Gaussian moments; hence $S_x$ converges in distribution to the standard Gaussian, by the method of moments \cite[Theorem 30.2]{Billingsley}.  Moreover, by \ref{mean square convergence} for a single $l$ and Markov’s inequality, $T_x$ converges to $0$ in probability. By Slutsky’s theorem \cite[Chapter 6, Theorem 15]{Rohatgi-Saleh}, it follows that $\mathcal S_x=S_x+T_x$ converges in distribution to the standard Gaussian.  This version may be useful when the tail $\sum_n(Z-Z_M)(X_n)$ can be controlled only on average.

    \end{enumerate}
\end{remark}

In some applications of the above theorem, the arithmetic information needed to compute the family averages of $\cos (2\pi mX_n)$ is available only outside a small subset $\mathcal G_{n,k}$ of the family $\mathcal F_k$.  For example, for the families $\mathcal F_{A,B}$ of elliptic curves, the identities relating the Frobenius angle $\theta_E(p)$ to $a_E(p^m)$ hold only for curves $E(a,b)$ such that $p \nmid ab\Delta(a,b)$.  The following variant of Theorem \ref{CLT general smooth} allows such subsets to be discarded, provided they are sufficiently thin in the sense of \eqref{thin-sets}.  We use this variant in the proof of Theorems \ref{Elliptic-curves-smooth}, \ref{Kloosterman-multiples-smooth} and \ref{Hecke-multiples-smooth}.

\begin{thm}\label{CLT general smooth variant}
We use the setup and notations in Theorem \ref{CLT general smooth}.  For each $n \in \mathcal A$, let $\mathcal G_{n,k} \subseteq \mathcal F_k$ be such that for every $u \geq 1$,
\begin{align}\label{thin-sets}
\left \langle\left(\sum_{ n \in \mathcal A \cap [1,x] }\chi_{\mathcal G_{n,k}}\right)^u \right \rangle_k= \o\left(A(x)^{u/2}\right).
\end{align}
Put
$$\mathcal H_{n,k} := \mathcal F_k \setminus \mathcal G_{n,k}. $$
For $n \in \mathcal A$ and $X \in \mathcal F_k$, define 
$$W(X_n) := \chi_{\mathcal H_{n,k}}(X)Z(X_n).$$
For $m \in \Z$ and $n \in \mathcal A$, define
$$\widetilde{c}_{m,n,k} := \left \langle \chi_{\mathcal H_{n,k}}(X) \cos (2\pi m X_n) \right \rangle_k = \frac{1}{|\mathcal{F}_k|}\sum_{X \in \mathcal H_{n,k}} \cos (2\pi m X_n).$$
Suppose that the following asymptotic estimates hold as $x \to \infty$.
\begin{enumerate}[label=\normalfont\bfseries (GC\arabic*)]
\item\label{Cond GC 1}{\bf Truncation.}
Denote 
\begin{align*}
W_M(X_n) = \chi_{\mathcal H_{n,k}}(X)Z_M(X_n).
\end{align*}
For every $l \geq 1$, 
$$ \Bigg\langle \left |\sum_{n \in \mathcal A \cap [1,x]} (W - W_M)(X_n)
\right| ^l\Bigg\rangle_k= \o\left(A(x)^{l/2}\right).$$

\item\label{Cond GC 2}{\bf Approximate independence.}

For all positive integers $u, r_1,r_2,\dots, r_u$, 
           \begin{align*} 
\sum_{\underline{n}}\sum_{\underline{m}} \abs{\prod_{i=1}^u \left \langle \chi_{\mathcal H_{n_i,k}}(X) \cos (2\pi m_i X_{n_i})\right \rangle_k-\left \langle \prod_{i=1}^u\chi_{\mathcal H_{n_i,k}}(X) \cos (2\pi m_i X_{n_i}) \right \rangle_k }
&= \o(A(x)^{r/2}).
            \end{align*}
\item\label{Cond GC 3}{\bf Convergence to the Weyl limits.}
For all positive integers $u,r_1,r_2,\dots,r_u$,
        \begin{align*}
\sum_{\underline{n}}\prod_{i=1}^u\sum_{|m_i| \leq Mr_i} \abs{\widehat{Z_M^{r_i}}(m_i)(\widetilde{c}_{m_i,n_i,k}-d_{m_i})}=\o\left(A(x)^{r/2}\right).
\end{align*}
\textup{(}We recall that $d_m$ is our usual Weyl limit of $(\mathcal F_k)_{k\geq 1}$ along $\mathcal A$ defined in \eqref{Weyl limit d_l}.\textup{)}
\end{enumerate}
As in Theorem \ref{CLT general smooth}, for $X \in \mathcal F_{k(x)}$, define
$$\mathcal S_x(X) := \frac{1}{\sqrt{A(x) V_h(\mu)}}\left(\sum_{n \in \mathcal A \cap [1,x]} h_n(X) - A(x)\int_0^1 h(t) \, \mu(t) \, dt\right).$$
Then, all the moments of $\mathcal S_x$ converge to the Gaussian moments as $x \to \infty$. 
\end{thm}
\begin{unmremark}
 Conditions \ref{Cond GC 1}-\ref{Cond GC 3} are obtained from  \ref{mean square convergence}-\ref{condn 2, smooth} by replacing, for each $n \in \mathcal A$ and $X \in \mathcal F_k$, the quantities $Z(X_n)$, $Z_M(X_n)$ and $\cos(2\pi mX_n)$ by $\chi_{\mathcal H_{n,k}}(X)Z(X_n)$, $\chi_{\mathcal H_{n,k}}(X)Z_M(X_n)$ and $\chi_{\mathcal H_{n,k}}(X)\cos(2\pi mX_n)$ respectively, and correspondingly, $c_{m,n,k}$ by $\widetilde{c}_{m,n,k}$.  Taking $\mathcal G_{n,k} = \emptyset$ for each $n\in \mathcal A$ recovers Theorem \ref{CLT general smooth}.  The sets $\mathcal G_{n,k}$ may also be all of $\mathcal F_k$ for some $n \in \mathcal A$ (for example, in Section \ref{sec: modular forms},  Setup~\ref{spaces-rvs-modular}).
\end{unmremark}
\begin{proof}
We denote, for a random variable $Y$ on $\mathcal F_k$ and a positive integer $\ell$, $\|Y\|_\ell := \left \langle |Y|^{\ell}\right \rangle_k^{1/\ell}$.  Also denote $\sigma_x := \sqrt{A(x)V_{h}(\mu)}$ and $C_h := 2 \sum_{m \in \Z} \abs{\widehat{h}(m)}$.

\hypertarget{step1'}{\smallskip\noindent\textbf{Step 1 (Decomposition of $\mathcal S_x$).}}

 Define the random variables $\mathcal U_x,\,\mathcal T_x$ and $\mathcal V_x$ on $\mathcal F_k$ as follows:
 \begin{align*}
     \mathcal U_x(X) &= \frac{1}{\sigma_x}\sum_{n \in \mathcal A \cap [1,x]}W_M(X_n),\\
     \mathcal T_x(X) &= \frac{1}{\sigma_x}\sum_{n \in \mathcal A \cap [1,x]}(W - W_M)(X_n)\\
    \intertext{and}
    \mathcal V_x(X) &= \frac{1}{\sigma_x}\sum_{n \in \mathcal A \cap [1,x]}\chi_{\mathcal G_{n,k}}(X)(Z(X_n)).
 \end{align*}
We denote $\mathcal E_x = \mathcal T_x + \mathcal V_x$.  Then, $\mathcal S_x = \mathcal U_x + \mathcal E_x.$

By \ref{Cond GC 1},  $\|\mathcal T_x\|_r \to 0$ for each $r \geq 1$ (as $x \to \infty$).  Since $|Z(t)| \leq C_h$, we have
$$\abs{\mathcal V_x(X)} \leq \frac{C_h}{\sigma_x}\sum_{n \in \mathcal A \cap [1,x]}\chi_{\mathcal G_{n,k}}(X).$$
Thus, by \eqref{thin-sets}, for each $r \geq 1$,
$$\|\mathcal V_x \|_r \leq \frac{C_h}{\sigma_x}\left \langle \left(\sum_{n \in \mathcal A \cap [1,x]}\chi_{\mathcal G_{n,k}}(X)\right)^r \right \rangle_k^{1/r} \to 0.$$
By Minkowski's inequality,
$\mathcal E_x = \mathcal T_x + \mathcal V_x$ satisfies $\|\mathcal E_x\|_r \to 0$ for every $r \geq 1$.

\hypertarget{step2'}{\smallskip\noindent\textbf{Step 2 (Moments of the truncated sums).}}

We show that $\displaystyle \left \langle \mathcal U_x^r \right \rangle \to \gamma_r$ for every $r \geq 1$ by following Steps \hyperlink{step2}{2}, \hyperlink{step3}{3} and \hyperlink{step4}{4} in the proof of Theorem \ref{CLT general smooth}, with $c_{m,n,k}$ replaced by $\widetilde{c}_{m,n,k}$.  Conditions \ref{condn 1, smooth} and \ref{condn 2, smooth} get replaced by Conditions \ref{Cond GC 2} and \ref{Cond GC 3} respectively.

Expanding $\mathcal U_x^r$ as in \eqref{higher-moments-MT-ET-gen} produces the averages 
$$\left \langle \prod_i \chi_{\mathcal H_{n_i,k}}(X) \cos (2\pi m_i X_{n_i})\right \rangle_k,$$
which are compared with $\prod_i \widetilde{c}_{m_i,n_i,k}$ by \ref{Cond GC 2}.  Note that unlike in \ref{condn 2, smooth}, the term $m_i = 0$ in \ref{Cond GC 3} need not vanish.  Since $d_0 = 1$, 
 $$\abs{\widetilde{c}_{0,n,k} - d_0 } = \frac{\abs{\mathcal G_{n,k}}}{\abs{\mathcal F_k}},$$
 and therefore, by equation \eqref{thin-sets},
 $\displaystyle \sum_{n \in \mathcal A \cap [1,x]}\abs{\widetilde{c}_{0,n,k} - d_0 } = \o\left(\sqrt{A(x)}\right)$.  (In applications, it is usually absorbed into the other terms.) 

\smallskip\noindent\textbf{Step 3 (Conclusion).}

Exactly as in Step \hyperlink{step5}{5} of the proof of Theorem \ref{CLT general smooth}, since $\mathcal S_x = \mathcal U_x + \mathcal E_x$, the binomial theorem and the Cauchy--Schwarz inequality give 
\begin{align*}
    \abs{\left \langle \mathcal S_x ^r \right \rangle_k - \left \langle \mathcal U_x ^r \right \rangle}_k \leq \sum_{j=1}^ r \left \langle \mathcal U_x ^{2(r-j)} \right \rangle_k^{1/2}\|\mathcal E_x\|_{2j}^j \to 0,
\end{align*}
since the even moments $\left \langle \mathcal U_x ^{2(r-j)} \right \rangle_k$ are bounded by Step \hyperlink{step2'}{2}, and $\|\mathcal E_x\|_{2j} \to 0$ for each $j$ by Step \hyperlink{step1'}{1}.  Thus, $\left \langle \mathcal S_x^r \right \rangle_k \to \gamma_r$ for each $r \geq 1$.
\end{proof}

Theorem \ref{CLT general smooth variant} can also be adapted to 2-periodic test functions, $h(t)=\sum_{m\in\Z} \widehat h(m)e(mt/2)$, with $\widehat h(m)$ real, even, and absolutely summable.  This leads to CLTs for a larger class of functions than those addressed by Theorem \ref{CLT general smooth variant}, such as polynomials in $\cos(\pi t)$.  
Another advantage of extending the scope of Theorem \ref{CLT general smooth variant} to 2-periodic test functions is that they can approximate the indicator function of an arbitrary interval $[a,b]\subset (0,1)$. In Theorem \ref{CLT gen linear}, we work with such a test function, constructed from the Beurling–Selberg polynomials.  We state this adaptation below.

\begin{thm}\label{2-periodic}

Assume the setup of Theorem \ref{CLT general smooth variant}, but let $H\in C(\mathbb R)$ be a real-valued, even, 2-periodic function,
$$H(t)=\sum_{m\in\mathbb Z}\widehat H(m)\cos(\pi mt),\qquad\sum_{m\in\Z}|\widehat H(m)|<\infty,$$
in place of $h$. For $m\in\mathbb Z$, put
$$C_{m,n,k}:=\Big\langle\mathbf 1_{\mathcal H_{n,k}}(X)\cos(\pi mX_n)\Big\rangle_k,\qquad D_m:=\lim_{x\to\infty}\frac{1}{A(x)}\sum_{\substack{n\leq x\\ n\in\mathcal A}}\big\langle\cos(\pi mX_n)\big\rangle_k,$$
and assume that the limits $D_m$ exist and satisfy $\sum_m|D_m|<\infty$. Suppose \ref{Cond GC 1}-\ref{Cond GC 3} hold with $\cos(2\pi mX_n)$, $\widetilde{c}_{m,n,k}$, $d_m$ and $\widehat h$ replaced by $\cos(\pi mX_n)$, $C_{m,n,k}$, $D_m$ and $\widehat H$. Then the conclusion of Theorem \ref{CLT general smooth variant} holds for the sums $\sum_{n\leq x,\,n\in\mathcal A}H(X_n)$, with $\mu$ replaced by the probability density
$$\nu(t) := \sum_{m \in \Z} D_m \cos (\pi mt),\quad t \in [0,1].$$

\end{thm}

\begin{proof}
Apply Theorem \ref{CLT general smooth variant} to the sequences $Y_n:=X_n/2\in[0,\frac{1}{2}]$ and the even, 1-periodic test function $h(t):=H(2t)$, noting that $h(Y_n)=H(X_n)$ and $\cos(2\pi mY_n)=\cos(\pi mX_n)$.
\end{proof}

\begin{remark}\label{rmk-2-periodic}
In the applications, the Weyl limits at $m=2j$ are those already computed, $D_{2j}=d_j$.  For the angles themselves (that is, for $\alpha = 1$), the odd frequencies $m$ give new Weyl sums, $\cos(\pi m\theta)=\frac{1}{2}\left(U_m(\cos\pi\theta)-U_{m-2}(\cos\pi\theta)\right)$, which in all three families involve only odd powers or odd-degree symmetric powers.  These are controlled by the same arithmetic input as that for even frequencies, and $D_m = 0$ for all odd frequencies $m$.  In particular, in all the three applications considered in this article, for $\alpha = 1$,
$$\nu(t) = \sum_{j \in \Z} d_j \cos(2\pi jt) = \mu(t).$$That is, the mean and variance are those of $H$ with respect to the original limiting measure $\mu(t)$.

\end{remark}

\subsection{A table for key steps}\label{sec: table}
In the following table, we indicate how the proofs of Theorems \ref{Elliptic-curves-smooth}-\ref{Hecke-multiples-smooth} are organized, pointing to the propositions and calculations for the verification of the hypotheses of Theorem \ref{CLT general smooth variant}.  The reader can refer to Section~\ref{sec: Notation} for notations.

\begin{table}[h!]
    \centering
    \begin{tabularx}{\textwidth}{|>{\raggedright\arraybackslash}p{3cm}|>{\centering\arraybackslash}X|>{\centering\arraybackslash}X|>{\centering\arraybackslash}X|} \hline
    
     & Frobenius angles\par (Section~\ref{sec: Elliptic curves}, Setup \ref{spaces-rvs-elliptic}) & Kloosterman angles\par (Section~\ref{sec: Kloosterman}, Setup \ref{Spaces-rvs-Kloosterman}) & Hecke angles\par (Section~\ref{sec: modular forms}, Setup \ref{spaces-rvs-modular}) \\ \hline
                                    
     $\f_k$  & $E \in \mathcal F_{A,B}$ & $a \in \mathcal F_N$  &  $f \in \mathcal F_{N,k}$  \\ \hline
     
     $\mathcal{A}$ & $\mathcal{P}_{\geq 5}$ &          $\mathcal{P}$ & $\mathcal{P}$ \\ \hline

    $(X_n)_{n \in \mathcal{A}}$          & $(\alpha\theta_E(p)\bmod\,1)_{p \geq 5}$  & $(\alpha\theta_{p,a}\bmod\,1)_p$ & $(\alpha\theta_f(p) \, \bmod 1)_p$ \\ \hline
    Thinness of $\mathcal G_{n,k}$              & $\mathcal G_{p,A,B}$ (Prop.~\ref{bad-red-points}\ref{omega_x-higher moments-ell}) & $\mathcal G_{p,N}$ (Lem.~\ref{thin-set-sum-Kl}) &  $\mathcal G_{p,N,k}$ (Eq.~\eqref{thin-sets-modular})\\ \hline 
    $\widetilde{c}_{m,n,k}$ (\& $d_m$)      & $\widetilde{c}_{m,p,A,B}$ (Prop.~\ref{Weyl limits elliptic}) & $\widetilde{c}_{m,p,N}$ (Prop.~\ref{Weyl limits Kloosterman})  & $\widetilde{c}_{m,p,N,k}$ (Prop.~\ref{Weyl-limits-modular}) \\ \hline

     Arithmetic input & Prop.~\ref{S(n)-elliptic}\ref{S0pm}--\ref{Spm-bounds} & Prop.~\ref{Katz} & Prop.~\ref{trace formula} \\ \hline

    \ref{Cond GC 1} & Prop.~\ref{Condition 1 elliptic} & Prop.~\ref{Condition 1 Kloosterman} & Prop.~\ref{Condition 1 modular} \\ \hline    
    \ref{Cond GC 2} & Prop.~\ref{Condition 2 elliptic}\ref{independence-rv-ell}, Eq.~\eqref{GC-2-estimate-ell} & Prop.~\ref{appl of crt} & Prop.~\ref{Towards GC 2}, Eq.~\eqref{GC-2-mod}\\ \hline
    \ref{Cond GC 3} & Eqs.~\eqref{GC3-inner}--\eqref{GC3-estimate-ell}     &  Eqs.~\eqref{GC3-inner-Kl}--\eqref{GC3-estimate-Kl} & Prop.~\ref{Towards GC 3}, Eq.~\eqref{GC-3-mod} \\ \hline

    Nagoshi's statistics & Cor.~\ref{Nagoshi-revisited-elliptic} & Cor.~\ref{Nagoshi-type-Kloosterman} & Cor.~\ref{Nagoshi-revisited} \\ \hline
   
    \end{tabularx}
    \vspace{0.1cm}
    \caption{Organization of Sections \ref{sec: Elliptic curves}, \ref{sec: Kloosterman} and \ref{sec: modular forms}.}
    \label{tab: organization of sections 3,4 and 5}
\end{table}

\subsection{Some useful  tools}\label{sec: some useful tools}

We now collect some tools and observations which will be used for all the applications of Theorem~\ref{CLT general smooth variant}.

\subsubsection{Test functions and Poisson summation}\label{sec: test-fns-PSF}
The following lemma shows that our general CLTs (Theorem \ref{CLT general smooth}-\ref{2-periodic}) are applicable to the test functions considered in Theorems \ref{Elliptic-curves-smooth}-\ref{Hecke-multiples-smooth}.

\begin{lemma}\label{h in setup5}
   Let $g \in \mathcal C(\R)$ be real-valued and even such that $\displaystyle g(t)\ll\frac{1}{(1+|t|)^\delta}$ for some $\delta>1$ and let $L \geq 1$.  
   \begin{enumerate}[label=\normalfont(\alph*)]
   \item\label{PSF for h}Then, $h(t) := \sum_{n \in \Z}g(L(t+n))$
   converges absolutely and uniformly, $h$ is continuous, even and 1-periodic, and $\displaystyle\widehat h(m)=\frac{1}{L}\,\widehat g\left(\frac {m}{L}\right)$ for every $m \in \Z$. 
\item\label{hm}  Suppose $\displaystyle\widehat g(t)\ll(1+|t|)^{-\kappa}$ for some $\kappa >1$.  Then, $\displaystyle \sum_m|\widehat h(m)|<\infty$, so $h(t)=\sum_m\widehat h(m)e(mt)$, with absolute and uniform convergence.
\item\label{m>M} Under the hypothesis of \ref{hm}, for any $M \geq 1$, we have $\displaystyle \sum_{|m| > M}\abs{\widehat{h}(m)} \ll_{L,\kappa} \frac{1}{M^{\kappa - 1}}$.
\item\label{mhm} If $\kappa>2$, then $\sum_m(1+|m|)\,|\widehat h(m)|<\infty$.
\end{enumerate}
\end{lemma}
\begin{proof}
For $t\in[0,1]$ and $|n|\ge2$, we have $|g(L(t+n))|\ll\left(1+L(|n|-1)\right)^{-\delta}$, which is summable since $\delta>1$. The Weierstrass M-test gives uniform convergence, and by periodicity this extends to all of $\mathbb R$. By a process of unfolding, we get
$$\widehat h(m)=\int_0^1\sum_ng\left(L(t+n)\right)e(-mt)\,dt=\int_{-\infty}^{\infty}g(Lt)\,e(-mt)\,dt=\frac{1}{L}\,\widehat g\left(\frac mL\right).$$
This proves \ref{PSF for h}.

Parts \ref{hm} and \ref{mhm} follow from the bound $|\widehat h(m)|\ll_{L,\kappa}(1+|m|)^{-\kappa}$.

Part \ref{m>M} follows from the inequality $\displaystyle \sum_{|m| > M}\abs{\widehat{h}(m)} \ll_L \int_{M}^{\infty} \frac{1}{(1 + t)^{\kappa}} dt\ll_{L,\kappa} \frac{1}{M^{\kappa - 1}}$.

\end{proof}

\begin{unmremark}
 Part \ref{mhm} will be required only in the CLTs for Frobenius and Kloosterman angles.  
\end{unmremark}

\subsubsection{Chebyshev Polynomials}\label{sec: Chebyshev}

Since all three families considered in the article are Sato--Tate distributed on average, it is natural to expand test functions in an orthonormal basis adapted to the Sato--Tate measure.  For $n \geq 0$, the Chebyshev polynomial of the second kind $U_n$ is the polynomial of degree $n$ determined as follows:
\begin{equation}\label{eq:chebyshev}
U_n(\cos (\pi t)) := \sum_{j = 0}^{n} e^{i(n-2j)\pi t} = \frac{\sin((n+1) \pi t)}{\sin \pi t}, \qquad t \in [0,1].    
\end{equation}

Let $\mu(t) = 2\sin^2(\pi t)$ 
and equip $L^2([0,1], \R)$ with the inner product
\begin{equation*}\label{inner prod wrt Sato--Tate measure}
  \langle f, g \rangle := \int_0^1 f(t)\, g(t)\, \mu(t)\, dt.
\end{equation*}
The functions $U_n(\cos (\pi t))$ are orthonormal with respect to this inner
product:
\begin{equation}\label{eq:mean zero}
  \int_0^1 U_n(\cos (\pi t))\, U_m(\cos (\pi t))\, \mu(t)\, dt = \delta_{m,n}.
\end{equation}

These polynomials satisfy the following properties:
\begin{lemma}
For all $n,m \geq 0$,
\begin{equation}\label{formula of U_n U_m}
    U_n(\cos (\pi t))U_m(\cos (\pi t))=\sum_{r=0}^{\min(m,n)}U_{n+m-2r}(\cos (\pi t)).
\end{equation}
In particular, for all $m \geq 1$, we have
\begin{equation}\label{recurrence of Chebyshev poly}
    U_{m+1}(\cos (\pi t))=(2\cos (\pi t))U_{m}(\cos (\pi t))-U_{m-1}(\cos (\pi t)).
\end{equation}
\end{lemma}

\subsubsection{Elementary product estimates}\label{sec: product-estimates}

We now record two elementary estimates for products, which will be used repeatedly to verify conditions \ref{Cond GC 2} and \ref{Cond GC 3}, in Sections \ref{sec: Elliptic curves}, \ref{sec: Kloosterman} and \ref{sec: modular forms}.

\begin{lemma}\label{product-estimates} \leavevmode

\begin{enumerate}[label=\normalfont(\roman*), align=left, labelsep=0em]

\item\label{product-sum} If $x_1,\dots,x_u,y_1,\dots,y_u\in\mathbb C$ satisfy $|x_i|,|y_i|\leq 1$, then
$$\left|\prod_{i=1}^ux_i-\prod_{i=1}^uy_i\right|\le\sum_{i=1}^u|x_i-y_i|.$$

\item\label{sum-product} If $F_1,\dots,F_u$ are nonnegative functions defined on $\mathcal A$, then
$$\sum_{\underline n}\prod_{i=1}^uF_i(n_i)\leq\prod_{i=1}^u\sum_{n \in \mathcal A \cap [1,x]}F_i(n),$$
where $\underline n$ runs over $u$-tuples of distinct elements $n_i\le x$ of $\mathcal A$.
\end{enumerate}
    
\end{lemma}
\begin{proof}
    Part \ref{product-sum} follows from the telescoping identity

$$\prod_{i=1}^ux_i-\prod_{i=1}^uy_i=\sum_{j=1}^u\left(x_1\cdots x_{j-1}\right)(x_j-y_j)\left(y_{j+1}\cdots y_u\right),$$
since each bracketed product has absolute value at most $1$. Part \ref{sum-product} holds because the right-hand side is the sum of the same nonnegative terms over all $u$-tuples of elements in $\mathcal A \cap [1,x]$, distinct or not.
\end{proof}

\section{Sequences arising from Elliptic curves}\label{sec: Elliptic curves}

In this section, we prove Theorem \ref{Elliptic-curves-smooth} by verifying the hypotheses of Theorem \ref{CLT general smooth variant} for the families $\mathcal F_{A,B}$. We first recall the arithmetic input (Propositions \ref{ell-curves-basic} and \ref{S(n)-elliptic}). After the setup (Setup \ref{spaces-rvs-elliptic}), we prove a bound for the exceptional sets (Proposition \ref{bad-red-points}), compute the Weyl limits (Proposition \ref{Weyl limits elliptic}), bound the tail (Proposition \ref{Condition 1 elliptic}), and establish approximate independence (Proposition \ref{Condition 2 elliptic}).  The proof of Theorem \ref{Elliptic-curves-smooth} is completed in Section \ref{proof of theorem 1}. We conclude with the corresponding result for Nagoshi-type statistics (Corollary \ref{Nagoshi-revisited-elliptic}).

For $(a,b)\in\Z^2$ with $\Delta(a,b)=4a^3+27b^2\neq 0$, Let $E = E(a,b)$ be the elliptic curve given by the equation
\begin{align*}
    y^2=x^3+ax+b.
\end{align*}
For a prime $p$, let $\tilde a_E(p):=p+1-N_E(p)$, where $N_E(p)$ is the
number of points on the projective curve
\[
Y^2Z=X^3+aXZ^2+bZ^3
\]
over $\mathbb F_p$. The model $E(a,b)$ has good reduction at
$p$ if and only if $p\nmid2\Delta(a,b)$. We consider the $L$-series
attached to the model $E(a,b)$,
\[
L(E,s)=\sum_{n=1}^\infty\frac{\tilde a_E(n)}{n^s}
:=\prod_{p\mid2\Delta(a,b)}\Big(1-\frac{\tilde a_E(p)}{p^s}\Big)^{-1}
\prod_{p\nmid2\Delta(a,b)}\Big(1-\frac{\tilde a_E(p)}{p^s}+\frac{1}{p^{2s-1}}\Big)^{-1},
\qquad \operatorname{Re}(s)>\tfrac32,
\]
whose local factors agree with those of the Hasse--Weil $L$-function of $E$
at every prime $p\nmid2\Delta(a,b)$.

By $a_E(n)$, we denote the normalized $n$-th coefficient, given by 
\begin{align*}
    a_E(n)=\frac{\widetilde{a}_E(n)}{\sqrt{n}}.
\end{align*}
 In order to prove Theorem \ref{Elliptic-curves-smooth}, we state a few properties of $a_E(n)$.  We refer the reader to \cite{Baier-Zhao} and \cite{Baier-Prabhu} for a curated discussion of these properties.

\begin{prop}\label{ell-curves-basic} Let $E = E(a,b)$ as defined above. 
    \begin{enumerate}[label=\normalfont(\roman*), align=left, labelsep=0em]
        \item  For a prime $p$,
        $\displaystyle
            \widetilde{a}_E(p)=\begin{cases}
                p+1-N_E(p)&\text{if  }p \nmid 2\Delta(a,b)\\
                -1,\,0 \text{ or }1&\text{if  }p \mid 2\Delta(a,b).
            \end{cases}
        $
        
        In the second case, $\widetilde{a}_E(2) = 0$, and $\widetilde{a}_E(3) = 0$ if $3 \mid \Delta(a,b)$.
        \item For any $n\in\N$, we have 
        \begin{align*}
            |a_E(n)|\leq d(n),
        \end{align*}
        where $d(n)$ is the number of divisors of $n$. In particular, if $p$ is a prime, then 
        \begin{align*}
            \abs{a_E(p)}\leq 2.
        \end{align*}
        \item \label{aE-prime-powers} 
        
        For a prime $p \geq 5$ and an integer $m \geq 1$ ,  
        \begin{align*}
            a_E(p^m) = \begin{cases}
                U_m\left(\frac{a_E(p)}{2}\right) &\text{ if }p \nmid \Delta(a,b)\\
                (a_E(p))^m &\text{ if }p \mid \Delta(a,b),\\
            \end{cases}
        \end{align*}
        where $U_m(x)$ denotes the $m$-th Chebyshev polynomial of the second kind defined in Section \ref{sec: Chebyshev}.
    \end{enumerate}
\end{prop}
\noindent By the above proposition, $|a_E(p)|\leq 2$ for every prime $p$, and hence, we may write 
$$a_E(p)=2\cos (\pi\theta_E(p)) \text{ for a unique }\theta_E(p)\in[0,1].$$
If $E$ has good reduction at a prime $p$, that is, if $p \nmid 2\Delta(a,b)$, we have 
\begin{align}\label{aE as Um}
a_E(p^m)=U_m(\cos(\pi\theta_E(p))).
\end{align}

 We are interested in the distribution of $a_E(p)$ as we vary $(a,b)$ over a suitable rectangle $|a|\leq A$ and $|b|\leq B$ for positive integers $A,B$. We recall the following proposition from \cite{Baier-Prabhu}. 
\begin{prop}\label{S(n)-elliptic}
    For $A,B \geq 1$, let
 $\displaystyle\mathcal F_{A,B} := \{E = E(a,b):\,|a| \leq A,\,|b| \leq B,\,4a^3 + 27b^2 \neq 0\}.$
    \begin{enumerate}[label=\normalfont(\roman*), align=left, labelsep=0em] 
    \item \label{S(n)}
      For an odd $n\in \N$, and for $\epsilon > 0$,
     \begin{align*}
        \frac{1}{\abs{\mathcal F_{A,B}}}\sum_{\substack {E(a,b) \in \mathcal F_{A,B} \\ (ab\Delta(a,b),n) = 1}}a_{E(a,b)}(n)=S(n)+ \O_{\epsilon,\omega(n)}\left(d(n)s(n)^{\frac{1}{2}+ \epsilon}\left(\frac{1}{A}+\frac{1}{B}\right)\right)
    \end{align*}
    where $d(n)$ denotes the number of positive divisors of $n$, $s(n)$ is the largest square free divisor of $n$, and $S(n)$ is defined as 
    \begin{align*}
        S(n) :=\frac{1}{(s(n))^2}\underset{(ab\Delta(a,b),n) = 1}{\sum_{a=1}^{s(n)}\sum_{b=1}^{s(n)}}a_{E(a,b)}(n).
    \end{align*}
    
    \item \label{Sn-multiplicative} $S(n)$ is multiplicative.

    \item \label{Spm} For a prime $p \geq 5$, and $k \in \N$,

\begin{align*}
    S(p^k)&=S_0(p^k)-S_1(p^k)-S_2(p^k)\\
\intertext{where}
    S_0(p^k)&=\frac{1}{p^2}\underset{p\nmid \Delta(a,b)}{\sum_{a=1}^{p}\sum_{b=1}^{p}}a_{E(a,b)}(p^k),\\
    S_1(p^k)&=\frac{1}{p^2}\sum_{a=1}^{p-1}a_{E(a,0)}(p^k)  =\O\left(\frac{k}{p}\right),\\
\intertext{and}
    S_2(p^k)&=\frac{1}{p^2}\sum_{b=1}^{p-1}a_{E(0,b)}(p^k) =\O\left(\frac{k}{p}\right).
\end{align*}

\item \label{S0pm}   For a prime $p\geq 5$, and $k\in\N$, we have 
    \begin{align*}
        S_0(p^k)=\begin{cases}\displaystyle-\left(1-\frac{1}{p}\right)\frac{1}{\sqrt{p}}\left(\frac{1}{p^{(k+1)/2}} +\sum_{f \in \mathcal F_{1,k+2}} a_f(p)\right)&\text{ if }k\text{ is even }\\
        0&\text{ otherwise,}
        \end{cases}
    \end{align*}
    where $\mathcal F_{1,k+2}$ denotes the set of normalized cuspidal Hecke eigenforms of weight $k+2$ with respect to $\mathrm{SL}_2(\Z)$, and for each such Hecke eigenform, $a_f(p)$ denotes its normalized $p$-th Fourier coefficient \textup{(}more details in Section \ref{sec: modular forms}\textup{)}. 
    
    \item \label{Spm-bounds}
  For a prime $p\geq 5$, and $k\in\N$, we have 
  $$S(p^{2k}) = \O\left( \frac{k}{\sqrt{p}}\right).$$
  \end{enumerate}
\end{prop}

\begin{proof}
For a proof of \ref{S(n)} and \ref{Sn-multiplicative},  we refer the reader to \cite[Section 4]{Baier-Prabhu}.  The proofs use two  facts: $a_E(n)$ is a multiplicative function of $n$, and for $(ab\Delta(a,b),n) = 1$, $a_E(n)$ depends only on $a$ and $b$ modulo $s(n)$.  An average over a large box therefore reduces to complete sums modulo $s(n)$, up to a boundary error.  Averaging over boxes of this kind has been used extensively to study Frobenius traces on average (see \cite{FM}, \cite{DP}, \cite{JY}, \cite{Baier-Zhao}). 

Part \ref{Spm} is \cite[Lemma 9.1]{Baier-Prabhu}.  Their origin is the work of Birch \cite{Birch}, who combined Deuring’s description of the number of elliptic curves over $\mathbb F_p$ with a given trace in terms of Hurwitz class numbers \cite{Deuring}, with the Eichler–Selberg trace formula.  In this way he expressed the moments of $a_E(p)$ over all curves modulo $p$ through traces of Hecke operators on $S_{k+2}(\mathrm{SL}_2(\mathbb Z))$, and showed that they converge to the Sato–Tate moments as $p\to\infty$.  

Combining \cite[(62)]{Baier-Prabhu} and \cite[Lemma 8.1]{Baier-Prabhu}, one obtains part \ref{S0pm}.  (The statement of \cite[Lemma 8.2]{Baier-Prabhu} omits the sign and the term $\frac{1}{p^{(k+1)/2}}$, but this does not affect the bounds derived there.)
    
Finally, by \ref{S0pm}, the Ramanujan--Deligne bound $|a_f(p)|\leq 2$, and $\dim S_{2k+2}(\mathrm{SL}_2(\mathbb Z))\ll k$, we get that for any integer $k \geq 1$, $S_0(p^{2k})  = \O\left({k}/{\sqrt{p}}\right)$.  By \ref{Spm}, we deduce \ref{Spm-bounds}.  Since Deligne’s bound is unconditional, so is \ref{Spm-bounds}; this is the only arithmetic input at a single prime that our proof of Theorem \ref{Elliptic-curves-smooth} requires.
\end{proof}

In order to derive Theorem \ref{Elliptic-curves-smooth}, we define the ingredients of Theorem \ref{CLT general smooth variant} in the context of families of elliptic curves.

\begin{setup}\label{spaces-rvs-elliptic}
For $A,B \geq 1$, let
 $\displaystyle\mathcal F_{A,B}$ be as in Theorem \ref{Elliptic-curves-smooth}, regarded as a finite probability space with the uniform measure.  Let $\mathcal A$ denote the set of primes greater than or equal to 5, so that $A(x) := |\mathcal A \cap [1,x]| = \pi(x) - 2$ for $x \geq 3$.
     For each such prime $p$, define 
    $$\mathcal G_{p,A,B} = \{E = E(a,b)\in \mathcal F_{A,B}:\,p \mid ab\Delta(a,b)\},$$
    and
    $$\mathcal H_{p,A,B} = \{E \in \mathcal F_{A,B}:\,p \nmid ab\Delta(a,b)\}.$$
To each positive integer $\alpha$ and $E \in \mathcal F_{A,B}$, we associate the sequence $\{\alpha\theta_E(p) \,\tmod\,1\}_{p \to \infty}.$
Define
\begin{align*}
\widetilde{c}_{m,p,A,B}(\alpha) &:= \frac{1}{\abs{\mathcal F_{A,B}}} \sum_{E \in \mathcal H_{p,A,B}}\cos (2\pi m\alpha\theta_E(p)),
\end{align*}
and
\begin{align*}
d_m(\alpha) &:= \lim_{x \to \infty}\frac{1}{{\pi}(x)-2}\sum_{5 \leq p \leq x}\frac{1}{\abs{\mathcal F_{A,B}}} \sum_{E \in \mathcal F_{A,B}}\cos (2\pi m\alpha\theta_E(p)).
\end{align*}
Define $\mu_{\alpha}(t) := \sum_{m \in \Z} d_m(\alpha) e(-mt)$.  We now adapt the setup of Theorems \ref{CLT general smooth} and \ref{CLT general smooth variant} to $\mathcal F_{A,B}$, and  consider the following functions:
 \begin{align*}
     h(\alpha \theta_E(p)) &= \textstyle\sum_{m \in \Z} \widehat{h}(m) e(m\alpha \theta_E(p)),\\
     h_M(\alpha\theta_E(p)) &= \textstyle\sum_{|m| \leq M}\widehat{h}(m) e(m\alpha \theta_E(p)),\\
     Z(\alpha \theta_E(p)) &= h(\alpha \theta_E(p)) - \textstyle\int_0^1 h(t) \mu_{\alpha}(t)\, dt  = h(\alpha \theta_E(p)) - \sum\nolimits_{m \in \Z}\widehat{h}(m) d_m(\alpha),\\
     Z_M(\alpha \theta_E(p)) &= h_M(\alpha \theta_E(p)) - \textstyle\int_0^1 h_M(t) \mu_{\alpha}(t)\, dt  = h_M(\alpha \theta_E(p)) - \sum\nolimits_{|m| \leq M}\widehat{h}(m) d_m(\alpha),\\
     W(\alpha \theta_E(p))  &\coloneqq \chi_{\mathcal H_{p,A,B}}(E)Z(\alpha \theta_E(p)),\\
     W_M(\alpha \theta_E(p)) &\coloneqq \chi_{\mathcal H_{p,A,B}}(E)Z_M(\alpha \theta_E(p)).
 \end{align*}

\end{setup}
In Propositions \ref{bad-red-points}-\ref{Condition 2 elliptic}, we verify the hypotheses and the conditions enumerated in Theorem \ref{CLT general smooth variant}.

In the following proposition, we show that the thin sets $\mathcal G_{p,A,B}$ have a negligible contribution when we execute moment computations for Theorem \ref{Elliptic-curves-smooth}.  
\begin{prop}\label{bad-red-points}
Let $A, B \geq 1$.

\begin{enumerate}[label=\normalfont(\roman*), align=left, labelsep=0em] 
\item\label{p-bad-prime-fixed} For a fixed prime $p$,
$$\frac{1}{\abs{\mathcal F_{A,B}}}\sum_{E \in \mathcal G_{p,A,B}} 1 \ll \frac{1}{p} + \frac{p}{A} + \frac{p}{B}.$$

\item\label{omega_x-higher moments-ell} For $E = E(a,b)\in \mathcal F_{A,B}$, define
$$ \omega_x(E) := \#\{p \leq x:\,p \mid ab\Delta(a,b)\}.$$
For any $u \geq 1$,
\begin{align*}
\left \langle\left(\sum_{ p \leq x }\chi_{\mathcal G_{p,A,B}}\right)^u \right \rangle = &\frac{1}{\abs{\mathcal F_{A,B}}}\sum_{E \in \mathcal F_{A,B}} \omega_x(E)^u\\
\ll_u& (\log \log x)^u + (x\pi(x))^u\left(\frac{1}{A} + \frac{1}{B}\right).
\end{align*}

\end{enumerate}
    
\end{prop}
\begin{proof}
 Note that 
    $$\frac{1}{4AB} \#\{|a| \leq A,|b| \leq B:\, p \mid 4a^3 + 27b^2\} = \O\left(\frac{1}{p}\right) + \O\left(\frac{p}{A} + \frac{p}{B}\right).$$
    This is because 
    $$\#\{(a,b) \in \mathbb F_p \times \mathbb F_p:\,4a^3 + 27b^2 \equiv 0 \,(\tmod p)\} = \O(p),$$
    and the total number of $\mathbb F_p \times \mathbb F_p$ grids inside the box $[-A,A] \times [-B,B]$ is $\leq (2A+1)(2B+1)/p^2$.  The number of integer points in the box $[-A,A] \times [-B,B]$ which are not covered by the above grids is at most $p(2A + 2B +2)$.  Moreover,
    $$\#\{|a| \leq A,|b| \leq B:\, p \mid ab\} \ll \frac{AB}{p} + p(A + B).$$
    
Thus,
\begin{equation*}
\frac{1}{4AB} \#\{|a| \leq A,|b| \leq B:\, p \mid ab(4a^3 + 27b^2)\} = \O\left(\frac{1}{p}\right) + \O\left(\frac{p}{A} + \frac{p}{B}\right).
\end{equation*} 
This proves \ref{p-bad-prime-fixed}.

To prove \ref{omega_x-higher moments-ell}, expand $\omega_x(E)^u$ as a sum over $u$-tuples of primes $p_1,p_2,\dots, p_u \leq x$, not necessarily distinct, with the condition that $p_i \mid ab\Delta(a,b)$ for each $i$. Fix one such tuple $(p_1,p_2,\dots,p_u)$, and let $p_1,p_2,\dots, p_j$ be distinct primes for some $1\leq j\leq u$. Put $q=p_1p_2\cdots p_j$.

Define
$\psi(a,b) := \chi_{q \mid ab\Delta(a,b)}$.  Under the bijection $(a,b)\,( \tmod q) \leftrightarrow ((a,b)\,( \tmod p_i))_i$, the function $\psi$ factors as a product of the indicator functions $\chi_{p_i \mid ab\Delta(a,b)}$.  Hence,
$$\frac{1}{q^2} \sum_{a,b \,(\tmod q)}\psi(a,b) = \prod_{i=1}^j \frac{\rho(p_i)}{p_i^2},\,\rho(p_i):= \#\{(a,b) \,(\tmod p_i):\,p_i \mid ab\Delta(a,b)\}.$$
As above, covering the box $[-A,A] \times [-B,B]$ by $q \times q$ blocks and accounting for the number of integer points $(a,b)$ which are not covered by the above grids, which is $\ll q(A+B)$, we get
\begin{align*}
\frac{1}{\abs{\mathcal F_{A,B}}}\sum_{E \in \mathcal F_{A,B}} \psi(a,b) = &\prod_{i=1}^j \frac{\rho(p_i)}{p_i^2} + \O\left(q\left(\frac{1}{A} + \frac{1}{B}\right)\right)\\
\ll&\frac{1}{q} + q\left(\frac{1}{A} + \frac{1}{B}\right).
\end{align*}
Expanding $\omega_x(E)^u$ as a sum of $u$-tuples of primes $p \leq x$ and grouping the tuples by their set $\{p_1,p_2,\dots,p_j\}$ of distinct primes ($1 \leq j \leq u$), we get
\begin{align*}
\left \langle \omega_x(E)^u \right \rangle \ll_u&\sum_{j=1}^u \sum_{p_1 < p_2 < \dots <p_j \leq x}\left(\frac{1}{p_1p_2\dots p_j} + p_1p_2\dots p_j\left(\frac{1}{A} + \frac{1}{B}\right)\right)\\
\ll_u&(\log \log x)^u + (x\pi(x))^u\left(\frac{1}{A} + \frac{1}{B}\right).
\end{align*}
This proves \ref{omega_x-higher moments-ell}. 
\end{proof}

We now derive the Weyl limits associated with this setup.

\begin{prop}\label{Weyl limits elliptic}
Let $A,B \geq 1$, and let $\alpha$ be an integer greater than or equal to 1. 
Then, the following properties hold.

\begin{enumerate}[label=\normalfont(\roman*), align=left, labelsep=0em]
\item\label{cmp elliptic} For a prime $p \geq 5$, an integer $m$, and a positive integer $\alpha$,
\begin{align*}
 \widetilde{c}_{m,p,A,B}(\alpha) = \begin{cases}
     1 + \O\left(\frac{1}{p} + \frac{p}{A} + \frac{p}{B}\right) &\text{ if }m=0\\
     \frac{1}{2}\left(S(p^2) - 1\right) + \O\left(\frac{1}{p} + \frac{p}{A} + \frac{p}{B}\right) &\text{ if }|m| = 1 = \alpha\\
      \frac{1}{2}\left(S(p^{2|m|\alpha}) - S(p^{2|m|\alpha-2})\right) + \O\left(|m|p^{\frac{1}{2} + \epsilon}\left(\frac{1}{A} + \frac{1}{B}\right)\right) &\text{ if }|m| \geq 2\text{ or }\alpha \geq 2.
 \end{cases}
 \end{align*}

\item\label{dm elliptic} Let $A = A(x)$ and $B = B(x)$ such that 
$$ \frac{x}{A} \text{ and }\frac{x}{B}\to0 \text{ as }x \to \infty.$$
Then, for an integer $m$,
\begin{align*}
&d_m(\alpha) = \begin{cases}
\phantom{-} 1&\text{ if }m=0,\\
-\frac{1}{2} &\text{ if } \alpha = 1,\,|m| = 1,\\
\phantom{-} 0&\text{ if } |m| \geq 2\text{ or }\alpha \geq 2.\\
\end{cases} 
\end{align*}
Therefore,
$$\mu_{\alpha}(t) = \sum_{m \in \Z} d_m(\alpha) e(-mt) = \begin{cases}
    2\sin^2(\pi t) &\text{ if }\alpha = 1\\
    1 &\text{ if }\alpha \geq 2.
\end{cases}$$
\end{enumerate}

\end{prop}

\begin{proof}

Applying Propositions \ref{ell-curves-basic}\ref{aE-prime-powers}, \ref{S(n)-elliptic}\ref{S(n)} and \ref{S(n)-elliptic}\ref{Spm}, we have, for a prime $p \geq 5$, and for $|m| \geq 1$,
\begin{align*}
&\frac{1}{\abs{\mathcal F_{A,B}}} \sum_{E \in \mathcal H_{p,A,B}} \cos (2\pi m\alpha\theta_{E(a,b)}(p))\\
=&    \frac{1}{\abs{\mathcal F_{A,B}}}\sum_{\substack{E \in \mathcal F_{A,B}\\(p,ab(4a^3 + 27b^2)) = 1}}\frac{1}{2}\left( a_{E(a,b)}(p^{2m\alpha})  -  a_{E(a,b)}(p^{2m\alpha-2})\right)\\
=&  
 \frac{1}{2}\left(S(p^{2|m|\alpha}) - S(p^{2|m|\alpha-2})\right)  
  + \O_{\epsilon}\left(|m| p^{\frac{1}{2}+ \epsilon}\left(\frac{1}{A} + \frac{1}{B}\right)\right).
\end{align*}
If $|m| = 1 = \alpha$, then
\begin{align*}
&\frac{1}{\abs{\mathcal F_{A,B}}} \sum_{E \in \mathcal H_{p,A,B}} \cos (2\pi m\alpha\theta_{E(a,b)}(p))\\ 
&=\frac{1}{2\abs{\mathcal F_{A,B}}}\sum_{E \in \mathcal H_{p,A,B}}a_E(p^2) - \frac{1}{2}\frac{\abs{\mathcal H_{p,A,B}}}{\abs{\mathcal F_{A,B}}}\\
&=\frac{1}{2}\left(S(p^2) + \O_{\epsilon}\left(p^{1/2+\epsilon}\left(\frac{1}{A} + \frac{1}{B}\right)\right)\right) -\frac{1}{2}\left(1 - \frac{\abs{\mathcal G_{p,A,B}}}{\abs{\mathcal F_{A,B}}}\right)\\
&= \frac{1}{2}\left(S(p^2) - 1\right) +\O\left(\frac{1}{p} + \frac{p}{A} + \frac{p}{B}\right),
\end{align*}
where the last estimate follows by Proposition \ref{bad-red-points}.
Finally, for $m = 0$,
\begin{align*}
\frac{1}{\abs{\mathcal F_{A,B}}} \sum_{E \in \mathcal H_{p,A,B}} \cos (2\pi m\alpha\theta_{E(a,b)}(p)) &= 1 - \frac{\abs{\mathcal G_{p,A,B}}}{\abs{\mathcal F_{A,B}}}\\
&= 1 +  \O\left(\frac{1}{p} + \frac{p}{A} + \frac{p}{B}\right).
\end{align*}
This proves \ref{cmp elliptic}.  Combining the above estimate with Proposition \ref{S(n)-elliptic}\ref{Spm-bounds}, and taking into account the primes of bad reduction with the help of Proposition \ref{bad-red-points} \ref{p-bad-prime-fixed}, we get
\begin{align*}
    &\frac{1}{{\pi}(x)-2}\sum_{p \leq x}\frac{1}{\abs{\mathcal F_{A,B}}} \sum_{E \in \mathcal F_{A,B}}\cos (2\pi m\alpha\theta_{E}(p))\\
    &= \frac{1}{{\pi}(x)}\sum_{5 \leq p \leq x} 
        \begin{cases}
         1&\text{ if }m=0\\
         \frac{1}{2}\left(S(p^{2|m|\alpha}) - S(p^{2|m|\alpha - 2})\right) &\text{ if }|m| \geq 1
        \end{cases}
    \\
    &\hspace{0.5cm}+ \O\left(\frac{1}{{\pi}(x)}\sum_{5 \leq p \leq x}\left((|m| +1)p\left(\frac{1}{A} + \frac{1}{B}\right)\right)\right)+ \O\left(\frac{\log \log x}{\pi(x)}\right).
\end{align*}
By Proposition \ref{S(n)-elliptic}\ref{Spm-bounds}, we get that for any integer $k \geq 1$,
$\displaystyle 
S(p^{2k})  = \O(k/\sqrt{p})$.  

We have the estimate $\displaystyle \sum_{p \leq x} \frac{1}{\sqrt{p}} \ll \sqrt{\frac{\pi(x)}{\log x}}$.  Thus, for a fixed $\alpha$ and $m \neq 0$,
\begin{align*}
    &\frac{1}{{\pi}(x)}\sum_{5 \leq p \leq x}\frac{1}{\abs{\mathcal F_{A,B}}} \sum_{E \in \mathcal F_{A,B}}\cos (2\pi m\alpha\theta_{E}(p))\\
    &= -\frac{1}{2} \delta_{|m|\alpha, 1} + \O\left(\frac{1}{{\pi}(x)}\sum_{5 \leq p \leq x} \frac{|m|\alpha}{\sqrt{p}}\right) + \O\left(\frac{1}{{\pi}(x)}\sum_{5 \leq p \leq x} |m|p\left(\frac{1}{A} + \frac{1}{B}\right)\right)\\
    &= -\frac{1}{2} \delta_{|m|\alpha,1} + \O\left(\frac{|m|\sqrt{\pi(x)}}{{\pi}(x)\sqrt{\log x}}\right) + \O\left(\frac{|m|x}{A} + \frac{|m|x}{B}\right) + \O\left(\frac{\log \log x}{\pi(x)}\right).
\end{align*}

 This proves \ref{dm elliptic}.
\end{proof}

The following proposition verifies \ref{Cond GC 1}. 

\begin{prop}\label{Condition 1 elliptic}
Let $h$ be as in Theorem \ref{Elliptic-curves-smooth} and let $\displaystyle M = M(x) \asymp x^{\frac{1}{2(\kappa - 1)}}$. Then
$$\sup_{E \in \mathcal F_{A,B}}\abs{\sum_{p \leq x} (W - W_M)(\alpha \theta_E(p))} = \o\left(\sqrt{\pi(x)}\right).$$

\end{prop}
\begin{proof} 

As $\widehat{Z}(m) = \widehat{h}(m)$ for $m \neq 0$, by Lemma \ref{h in setup5}\ref{m>M}, we obtain
$$ \abs{\sum_{p \leq x} (W - W_M)(\alpha \theta_E(p))} \leq \sum_{p \leq x} \sum_{|m| > M}\abs{\widehat{Z}(m)} \ll \frac{\pi(x)}{M^{\kappa -1}}\text{ for each }E \in \mathcal F_{A,B}.$$
\end{proof}

The next proposition provides the approximate independence required in condition \ref{Cond GC 2}. As $(a,b)$ ranges over a large box, the reductions of $E(a,b)$ modulo distinct primes behave independently, by the Chinese remainder theorem. Thus, joint averages of the angles at distinct primes factor into products of individual averages, up to an error that is small when the box is large compared with the primes involved.

\begin{prop}\label{Condition 2 elliptic}
     Let $p_1,p_2,\dots,p_u$ be $u$ distinct primes greater than or equal to 5, and let $q=p_1p_2\cdots p_u$.  For a prime $p$ and an  integer $k \geq 0$, put
     $$S_p(p^k) := \frac{1}{p^2}\sum_{\substack {a,b\,\tmod p \\ p \nmid ab\Delta(a,b)}}a_{E(a,b)}(p^k).$$
\begin{enumerate}[label=\normalfont(\roman*), align=left, labelsep=0em] 

\item \label{traces-elliptic} Let $n_1,n_2,\cdots,n_u$ be nonnegative integers.  Then,
 \begin{equation*}
  \frac{1}{\abs{\mathcal F_{A,B}}} \sum_{\substack{E \in \mathcal F_{A,B}\\(ab\Delta(a,b),q)=1}}\prod_{i=1}^u a_E(p_i^{n_i}) = \prod_{i=1}^u S_{p_i}(p_i^{n_i}) + \O_{u}\left(\prod_{i=1}^u (n_i +1)p_i\right)\left(\frac{1}{A} + \frac{1}{B}\right).
 \end{equation*}
 \item \label{independence-rv-ell}  Let $m_1,m_2,\dots,m_u$ be $u$  integers, and let $\alpha$ be a positive integer. Then,
    \begin{equation*}
    \begin{aligned}
         &\frac{1}{\abs{\mathcal F_{A,B}}} \sum_{\substack{E \in \mathcal F_{A,B}\\(ab\Delta(a,b),q)=1}}\prod_{i=1}^u \cos 2\pi m_i\alpha\theta_E(p_i)
         =\prod_{i=1}^u \widetilde{c}_{m_i,p_i,A,B}(\alpha) +\O_{u}\left(\prod_{i=1}^u (|m_i| \alpha +1)p_i\right)\left(\frac{1}{A} + \frac{1}{B}\right).
    \end{aligned}
    \end{equation*}
    \end{enumerate}
\end{prop}
\begin{proof}
Define
$$\Phi(a,b) := \begin{cases}
    \prod_{i=1}^ u a_{E(a,b)}(p_i^{n_i}) &\text{ if }(ab\Delta(a,b),q) = 1\\
    0 &\text{ otherwise.}
\end{cases}$$
$\Phi(a,b)$ is periodic modulo $q$ in both $a$ and $b$, and $|\Phi(a,b)| \leq \prod_{i=1}^u (n_i + 1)$.  The Chinese Remainder Theorem identifies $(a,b)\, \tmod q$ with the tuple $((a,b)\, \tmod p_i)_i$.  Thus, 
$$\frac{1}{q^2}\sum_{a,b\,\tmod q}\Phi(a,b) = \prod_{i=1}^u \frac{1}{p_i^2}\sum_{\substack {a,b\,\tmod p_i \\ p_i \nmid ab \Delta(a,b)}} a_{E(a,b)}(p_i^{n_i}) = \prod_{i=1}^u S_{p_i}(p_i^{n_i}).$$
We cover the box $[-A,A] \times [-B,B]$ by complete $q \times q$ grids.  The number of integer points $(a,b)$ which are not covered by the above grids is $\ll q(A+B)$. (We note that pairs with $\Delta(a,b) = 0$ are automatically excluded, since every prime divides $\Delta$.)  Thus,
\begin{equation*}
  \frac{1}{\abs{\mathcal F_{A,B}}} \sum_{\substack{E \in \mathcal F_{A,B}\\(ab\Delta(a,b),q)=1}}\prod_{i=1}^u a_E(p_i^{n_i}) = \prod_{i=1}^u S_{p_i}(p_i^{n_i}) + \O_{u}\left(q\left(\frac{1}{A} + \frac{1}{B}\right)\left(\prod_{i=1}^u (n_i +1)\right)\right).
 \end{equation*}
This proves part \ref{traces-elliptic}.

Denote, for each $1 \leq i \leq u$,
$$T_i = \begin{cases}\frac{1}{2}\left(S_{p_i}(p_i^{2|m_i|\alpha})-S_{p_i}(p_i^{2|m_i|\alpha-2})\right) &\text{ if }|m_i| \geq 1,\\
S_{p_i}(1) &\text{ if }m_i = 0.
\end{cases}$$
   Applying part \ref{traces-elliptic} with $u = 1$ and $q = p_i$,
   \begin{align}\label{Part i single prime}
   \frac{1}{\abs{\mathcal F_{A,B}}} \sum_{E \in \mathcal H_{p_i,A,B}}\cos(2\pi m_i\alpha\theta_E(p_i)) = T_i + \O\left((|m_i|\alpha+1)p_i\left(\frac{1}{A}+\frac{1}{B}\right)\right).
   \end{align}
Writing each cosine as a combination of at most two coefficients of the form $a_E(p^{2n})$, expanding the product, and applying part \ref{traces-elliptic} to each of the (at most $2^u$) terms, the main terms combine to form $\prod_i T_i$.  We have
    \begin{align*}
         \frac{1}{\abs{\mathcal F_{A,B}}} \sum_{\substack{E \in \mathcal F_{A,B}\\(ab\Delta(a,b),q)=1}}\prod_{i=1}^u \cos(2 \pi m_i \alpha\theta_{E(a,b)}(p_i))
        = \prod_{i=1}^u T_i +\O_u\left(\prod_{i=1}^u\left((|m_i|\alpha+1)p_i\right)\left(\frac{1}{A}+\frac{1}{B}\right)\right).
         \end{align*}
Note that $|T_i| \leq 1$, since $T_i$ is the average of $\cos(2\pi m_i\alpha\theta)$ over the classes $(a,b)\,\tmod p_i$ with $p_i\nmid ab\Delta(a,b)$, normalized by $p_i^2$.  By \eqref{Part i single prime} and Lemma \ref{product-estimates}\ref{product-sum} applied to $\prod_{i=1}^u T_i - \prod_{i=1}^u \widetilde{c}_{m_i,p_i,A,B}(\alpha)$ and the triangle inequality,
\begin{align*}
        & \biggl|\frac{1}{|\mathcal F_{A,B}|} \sum_{\substack{E \in \mathcal F_{A,B}\\(ab\Delta(a,b),q)=1}}\prod_{i=1}^u \cos(2 \pi m_i \alpha\theta_{E(a,b)}(p_i)) - \prod_{i=1}^u \widetilde{c}_{m_i,p_i,A,B}(\alpha)\biggr|\\
        \leq&\biggl|\frac{1}{|\mathcal F_{A,B}|} \sum_{\substack{E \in \mathcal F_{A,B}\\(ab\Delta(a,b),q)=1}}\prod_{i=1}^u \cos(2 \pi m_i \alpha\theta_{E(a,b)}(p_i)) - \prod_{i=1}^u T_i\biggr| 
        + \Bigl|\prod_{i=1}^u T_i - \prod_{i=1}^u \widetilde{c}_{m_i,p_i,A,B}(\alpha)\Bigr|\\  
    \leq& \sum_{i=1}^u \bigl|T_i - \widetilde{c}_{m_i,p_i,A,B}\bigr| +\O_u\left(\prod_{i=1}^u\left((|m_i|\alpha+1)p_i\right)\left(\frac{1}{A}+\frac{1}{B}\right)\right)
=\O_u\left(\prod_{i=1}^u\left((|m_i|\alpha+1)p_i\right)\left(\frac{1}{A}+\frac{1}{B}\right)\right).
\end{align*}
\end{proof}

\subsection{Proof of Theorem \ref{Elliptic-curves-smooth}}\label{proof of theorem 1}
 We now prove Theorem \ref{Elliptic-curves-smooth} as an application of Theorem \ref{CLT general smooth variant}.  

 Since $A = A(x),\,B= B(x)$ satisfy $$\frac{\log A}{\log x},\,\frac{\log B}{\log x} \to \infty,$$
 we get that
 \begin{equation}\label{growth of A and B}
 \lim_{x \to \infty} \frac{x^{\beta}}{A} = 0 = \lim_{x \to \infty} \frac{x^{\beta}}{B} \text{ for any }\beta >0.
 \end{equation}
Recall the notation of Setup \ref{spaces-rvs-elliptic}.  Recall that $h$  satisfies the properties mentioned in Setup \ref{spaces-rvs} (see Lemma \ref{h in setup5}).

For $E = E(a,b) \in \mathcal F_{A,B}$, $\sum_{p \leq x }\chi_{\mathcal G_{p,A,B}}(E) \leq \omega_x(E)$.  By Proposition \ref{bad-red-points}\ref{omega_x-higher moments-ell} and the above growth condition,
$$\left \langle\left(\sum_{ p \leq x }\chi_{\mathcal G_{p,A,B}}\right)^u \right \rangle_k= \o\left(\pi(x)^{u/2}\right).$$
Hence, the hypothesis in equation \eqref{thin-sets} is satisfied.  We now show that the conditions \ref{Cond GC 1}-\ref{Cond GC 3} are satisfied for $M(x) \asymp x^{1/(2(\kappa-1))}$.

By Proposition \ref{Condition 1 elliptic}, Condition \ref{Cond GC 1} holds for every $l \geq 1$.

Recall that $\sum_{\underline{p}}$ denotes the sum over $u$-tuples of distinct primes lying between 5 and $x$.  By Proposition \ref{Condition 2 elliptic}\ref{independence-rv-ell} applied to each prime $p_i$ such that $5 \leq p_i \leq x$,
\begin{align}\label{GC-2-estimate-ell}
&\sum_{\underline{p}}\sum_{\underline{m}} \left\lvert \Bigg \langle \prod_{i=1}^u \chi_{\mathcal H_{p_i,A,B}}(E)\cos(2\pi m_i \alpha \theta_E(p_i)) \Bigg \rangle_{\mathcal F_{A,B}}
   - \prod_{i=1}^u \widetilde{c}_{m_i,p_i,A,B}(\alpha)\right \rvert\\
\ll& \sum_{\underline{p}}\sum_{\underline{m}}\left(\prod_{i=1}^u (|m_i| \alpha + 1)p_i\right)\left(\frac{1}{A} + \frac{1}{B}\right)\notag\\
 \ll_{ \alpha}&M^{2u}\pi(x)^u x^u   \left(\frac{1}{A} + \frac{1}{B}\right)\notag\\
=& \o(\pi(x)^{r/2}),\notag
\end{align}
where the last asymptotic follows by the growth conditions in \eqref{growth of A and B}.  Condition \ref{Cond GC 2} is therefore verified.
Before verifying \ref{Cond GC 3}, we observe that by Lemma \ref{h in setup5}\ref{mhm}, we have
\begin{align}\label{mZm}
\sum_{|m| \leq M} (|m| + 1)\abs{\widehat{Z_M}(m)} \ll 
1,\end{align}
which is where we require  $\kappa > 2$.

By Proposition \ref{S(n)-elliptic}\ref{Spm-bounds} and Proposition \ref{Weyl limits elliptic}, we have, for each $1 \leq i \leq u$,
\begin{align*}
\abs{\widetilde{c}_{m_i,p_i,A,B}(\alpha)-d_{m_i}(\alpha)} \ll \frac{1}{p_i}\delta_{m_i,0} + \frac{1}{p_i}\delta_{|m_i|\alpha,1}  + \frac{|m_i|\alpha + 1}{\sqrt{p_i}} + (|m_i|\alpha + 1)p\left(\frac{1}{A} + \frac{1}{B}\right).
\end{align*} 
    Thus,
\begin{align}\label{GC3-inner}
&\sum_{|m_i| \leq Mr_i} \abs{ \widehat{Z_M^{r_i}}(m_i) (\widetilde{c}_{m_i,p_i,A,B}(\alpha)-d_{m_i}(\alpha))}\notag\\
&\ll_{\alpha} \left(\sum_{|m_i| \leq Mr_i}(|m_i| + 1)\abs{\widehat{Z_M^{r_i}}(m_i)}\right)\left(\frac{1}{\sqrt{p_i}} + p_i\left(\frac{1}{A} + \frac{1}{B}\right)\right).
\end{align}
As the function $w(m) = |m| + 1$ is submultiplicative,
\begin{align}\label{mZmpowers}
\sum_{|m_i| \leq Mr_i}(|m_i| + 1)\abs{\widehat{Z_M^{r_i}}(m_i)} \leq\left ( \sum_{|m_i|\leq M}(|m_i| + 1)\abs{\widehat{Z_M}(m_i)} \right)^{r_i} \ll 1, 
\end{align}
by \eqref{mZm}.
Combining the above with \eqref{GC3-inner} and Lemma \ref{product-estimates}\ref{sum-product}, and under the growth conditions \eqref{growth of A and B},
\begin{align}\label{GC3-estimate-ell}
&\sum_{\underline{p}}\prod_{i=1}^u \sum_{|m_i| \leq Mr_i} \abs{ \widehat{Z_M^{r_i}}(m_i) (\widetilde{c}_{m_i,p_i,A,B}(\alpha)-d_{m_i}(\alpha))}\\
    &\ll\sum_{\underline{p}}\prod_{i=1}^u \left(\frac{1}{\sqrt{p_i}} + p_i\left(\frac{1}{A} + \frac{1}{B}\right)\right)\notag\\
    &\ll_u\left(\frac{\pi(x)}{\log x}\right)^{u/2} + x^u\pi(x)^u \left(\frac{1}{A} + \frac{1}{B}\right)\notag\\
&=\o\left(\pi(x)^{r/2}\right).\notag
    \end{align}    
This proves Condition \ref{Cond GC 3}.  Thus, by Theorem \ref{CLT general smooth variant}, the moments of the normalized random variable
$$\widetilde{\mathcal S}_x(E) := \frac{1}{\sqrt{(\pi(x) - 2)V_h(\mu_{\alpha})}}\left(\sum_{5 \leq p \leq x} h(\alpha \theta_E(p)) - (\pi(x)-2)\int_0^1 h(t) \mu_{\alpha}(t) \, dt\right)$$
converge to the Gaussian moments as $x \to \infty$.  

We now pass from $\widetilde{\mathcal S}_x$ to $\mathcal S_x$.  We have
\begin{align*}
&\sum_{p \leq x} h(\alpha \theta_E(p)) - \pi(x)\int_0^1 h(t)\mu_{\alpha}(t) dt \\
=&\sum_{5 \leq p \leq x} h(\alpha \theta_E(p)) - (\pi(x) - 2)\int_0^1 h(t)\mu_{\alpha}(t) dt + \sum_{p \in \{2,3\}}Z(\alpha\theta_E(p)).
\end{align*}
That is, $\mathcal S_x = c_x(\widetilde{\mathcal S}_x + \mathcal R_x),$ where
$$c_x := \sqrt{\frac{\pi(x) -2}{\pi(x)}},\,R_x(E) = \frac{\sum_{p \in \{2,3\}}Z(\alpha\theta_E(p))}{\sqrt{(\pi(x)-2)V_h(\mu_{\alpha})}}.$$
As $\displaystyle \abs{\sum_{p \in \{2,3\}}Z(\alpha\theta_E(p))} \leq 2\sum_{m \in \Z}\abs{\widehat{Z}(m)} \ll 1$ uniformly for $E \in \mathcal F_{A,B}$, we get that $\displaystyle \abs{R_x(E)} \ll \frac{1}{\sqrt{\pi(x)}}$ and $c_x \to 1$ as $x \to \infty$.  Thus, by the binomial theorem and the Cauchy--Schwarz inequality, we have, for every $r \geq 1$,
\begin{align*}
\abs{\left \langle \left(\widetilde{\mathcal S}_x + \mathcal R_x \right)^r \right \rangle - \left \langle \left(\widetilde{\mathcal S}_x\right)^r\right\rangle } 
\leq&  \sum_{j=1}^r \binom{r}{j}\|R_x \|_{\infty}^j\left \langle \left(\widetilde{\mathcal S}_x\right)^{2(r-j)}\right\rangle^{1/2} \to 0,\\
\end{align*}
as the moments $\left \langle \left(\widetilde{\mathcal S}_x\right)^{2(r-j)}\right\rangle$ converge, and are therefore, bounded.  Thus,
$$\left \langle \mathcal S_x^r\right \rangle = c_x^r\left \langle \left(\widetilde{\mathcal S}_x + \mathcal R_x \right)^r \right \rangle \to \gamma_r \text{ for all }r \geq 1.$$
This proves Theorem \ref{Elliptic-curves-smooth}. 

\begin{remark}\label{Comments MRH free}

The only arithmetic input at a single prime is the bound $S(p^{2m})\ll m/\sqrt p$ (Proposition \ref{S(n)-elliptic}\ref{S0pm} and \ref{Spm-bounds}), which follows from Birch’s formula and Deligne’s bound. In the verification of \ref{Cond GC 3}, the weights $\widehat{Z_M^r}(m)$ make this bound summable over the frequencies (see \eqref{mZm}–\eqref{mZmpowers}), so no cancellation over primes in sums $\sum_{p\leq x}a_f(p)$ is needed.  Such cancellation, uniformly in the weight, is where the modular Riemann hypothesis was used in \cite[Section 9]{Baier-Prabhu} and \cite[Lemma 3.5]{BPS}.

\end{remark}
The general CLT for 2-periodic functions, Theorem \ref{2-periodic}, leads to the following corollary.
\begin{cor}\label{Nagoshi-revisited-elliptic}
    Fix an integer $j\geq 1$, and let $A=A(x)$ and $B=B(x)$ satisfy the conditions 
    $$\frac{\log A}{\log x},\,\frac{\log B}{\log x} \to \infty.$$
    Then, as $E$ ranges over $\mathcal F_{A,B}$, the moments of  $\displaystyle\frac{1}{\sqrt{\pi(x)}}\sum_{p\leq x}a_E(p^j)$ converge to the Gaussian moments. 
\end{cor}

\begin{proof}
    Let $\alpha = 1$. Denote $C_{m,p,A,B}(1) = C_{m,p,A,B}$ and $D_m(1) = D_m$.  Choose $H(t) = U_j(\cos (\pi t))$ for $j \geq 1$ and apply Theorem \ref{2-periodic} with the conventions of Setup \ref{spaces-rvs-elliptic} with exceptional sets $\mathcal G_{p,A,B}$.  As in the proof of Theorem \ref{Elliptic-curves-smooth}, the moment form of the thinness condition \eqref{thin-sets} is satisfied by Proposition \ref{bad-red-points}\ref{omega_x-higher moments-ell}.  Since $H$ is a trigonometric polynomial of degree $j$ in $\cos(\pi t)$, \ref{Cond GC 1} holds trivially with $M=j$ and $Z_M = Z$.  
    For $m \neq 0$, $p \geq 5$ and $p \nmid ab\Delta(a,b)$, 
    $$\cos (\pi m \theta_E(p)) = \frac{1}{2}\left(a_E(p^{|m|}) - a_E(p^{|m| - 2})\right) \text{ with }a_E(p^{-1}) := 0.$$
    Applying Proposition \ref{Condition 2 elliptic}\ref{traces-elliptic} with $u=1$ and $m \neq 0$, and setting $S_p(p^{-1}) := 0$, we get
    $$C_{m,p,A,B}(1) = \frac{1}{2}\left(S_p(p^{|m|}) - S_p(p^{|m|-2})\right) + \O\left((|m| + 1)p\left(\frac{1}{A} + \frac{1}{B}\right)\right). $$
      We also have
    $$C_{0,p,A,B} = S_p(1) + \O\left(p\left(\frac{1}{A} + \frac{1}{B}\right)\right).$$
    If $m = 2l$, this is the estimate $\widetilde{c}_{l,p,A,B}(1)$ in Proposition \ref{Weyl limits elliptic}.  Thus, $D_m = d_l(1)$.  Since $S_0(p^k) = 0$ for odd $k \geq 1$ (Proposition \ref{S(n)-elliptic}\ref{S0pm}), we get, by Proposition \ref{S(n)-elliptic}\ref{Spm}, $S_p(p^k) = \O\left(k/p\right)$.  Hence, for odd $m$, 
    $$C_{m,p,A,B} \ll \frac{|m|}{p} + |m|p\left(\frac{1}{A} + \frac{1}{B}\right),$$
    and $D_m = 0$.  For all $m$, therefore,
    $$\abs{C_{m,p,A,B} - D_m} \ll  \frac{|m|+1}{\sqrt p} + (|m|+1)\,p\Big(\frac1A+\frac1B\Big),$$
    which is the one-prime bound also used in the proof of Theorem \ref{Elliptic-curves-smooth}.  This leads to \ref{Cond GC 3}.  Moreover, the proof of Proposition \ref{Condition 2 elliptic}\ref{independence-rv-ell} can be applied verbatim to derive \ref{Cond GC 2}.
    
      We have  
    $$\nu(t) = \sum_{l \in \Z} d_l \cos (2 \pi lt) = 2\sin^2( \pi t) = \mu_1(t).$$
By equation \eqref{eq:mean zero}, $\displaystyle\int_0^1 H(t) \nu(t) dt= 0$ and by equations \eqref{eq:mean zero} and \eqref{formula of U_n U_m}, $\displaystyle\int_0^1 H(t)^2 \nu(t) dt = 1.$  Thus, $V_H(\nu) = 1$.    By Theorem \ref{2-periodic} and the mechanism of passing from 
$\displaystyle\frac{1}{\sqrt{\pi(x)-2}}\sum_{5 \leq p\leq x}a_E(p^j)$ to $\displaystyle\frac{1}{\sqrt{\pi(x)}}\sum_{p\leq x}a_E(p^j)$ described in the proof of Theorem \ref{Elliptic-curves-smooth}, we obtain Corollary \ref{Nagoshi-revisited-elliptic}.
\end{proof}

\section{Sequences arising from Kloosterman sums}\label{sec: Kloosterman}

In this section, we prove Theorem \ref{Kloosterman-multiples-smooth}. The argument follows that of Section \ref{sec: Elliptic curves}, with Katz’s bounds (Proposition \ref{Katz}) for symmetric powers of the Kloosterman sheaf in place of Birch’s formula and the Ramanujan--Deligne bound.

For a prime $p$, let $\F_p$ denote the finite field of order $p$, and let $\F_p^{\times}$ denote the set of nonzero elements of $\F_p$. For an integer $a$, the Kloosterman sum $K(p,a)$ is defined as follows:
\begin{equation*}
    K(p,a) = \sum_{x \in \F_p^{\times}} e\left(\frac{x + ax^{-1}}{p}\right), \, e(t) = e^{2\pi i t}.
\end{equation*} 

Substituting $x$ with $-x$, we observe that $K(p,a)$ is a real number. Thus, using the Weil bound (see \cite{A_Weil}) for Kloosterman sums, we have that for each prime $p \nmid a$,
$$K(p,a) \in [-2\sqrt{p},2\sqrt{p}].$$
We denote
$$w(p,a) = \frac{K(p,a)}{2\sqrt{p}}.$$
We now choose $\theta_{p,a} \in [0,1]$ so that
$\cos (\pi \theta_{p,a}) = w(p,a).$
The following result states that as $p \to \infty$, the angles $\left\{ \theta_{p,a}:\,1 \leq a \leq p-1\right\}$
are equidistributed in the interval $[0,1]$  with respect to the Sato--Tate measure $2\sin^2(\pi t)dt$. 
\begin{prop}[Katz]\label{Katz}
For every $m \geq 1$, we have 
\begin{align}\label{Katz inequality}
    \biggl|\sum_{a\in\F_p^{\times}} U_m(\cos \pi \theta_{p,a}) \biggr| = \biggl|\sum_{a\in\F_p^{\times}}\frac{\sin((m+1)\pi\theta_{p,a})}{\sin(\pi\theta_{p,a})}\biggr|\leq\frac{(m+1)\sqrt{p}}{2}.
\end{align}
Thus, for any continuous function $h$ on $[0,1]$,
    \begin{equation}\label{VST}
       \lim_{p\to\infty} \frac{1}{p}\sum_{a = 1}^{p} h(\theta_{p,a}) = \int_{0}^{1} h(t) \, \mu(t)dt 
    \end{equation}
    where $\mu(t)=2\sin^2(\pi t)$.
    \end{prop}

\begin{proof}
The sum $\displaystyle \sum_{a\in\F_p^{\times}} U_m(\cos \pi \theta_{p,a})$ is the complete sum over $\F_p^{\times}$ of the normalized trace function of the $m$-th symmetric power of the Kloosterman sheaf.  
The inequality in \eqref{Katz inequality} is a well-known inequality due to Katz for this sum. We refer the readers to \cite[Example 13.6]{Katz} for more details.
To prove \eqref{VST}, it is sufficient to prove this asymptotic for $h(t) := U_m(\cos (\pi t))$ for every integer $m \geq 0$, since the linear span of $\{U_m(\cos (\pi t)):\,m\geq 0\}$ is dense in $C([0,1], \R)$.  By \eqref{Katz inequality}, if $m \geq 1$, then
$$\left| \frac{1}{p}\sum_{a\in\F_p^{\times}}U_m(\cos \pi \theta_{p,a})\right| \ll \frac{m}{\sqrt{p}} \to 0\text{ as }p \to \infty.$$
Thus, for each $m \geq 0$,
$\displaystyle \lim_{p\to\infty} \frac{1}{p}\sum_{a\in\F_p^{\times}}U_m(\cos \pi\theta_{p,a}) = \delta_{m,0} = \int_0^1 U_m(\cos \pi t) \mu(t) dt.$
\end{proof}

We will refer to \eqref{Katz inequality} as the Katz inequality and \eqref{VST} as the ``vertical'' Sato--Tate law for Kloosterman angles modulo $p$.  Katz also conjectured a  ``horizontal'' analogue of the above equidistribution theorem, namely that for a fixed positive integer $a$, the sequence $(\theta_{p,a})_{p \to \infty}$ is $\mu$-equidistributed in the interval $[0,1]$, where $\mu(t) = 2\sin^2(\pi t)$.  This conjecture remains open, and implies that for any integer $\alpha \geq 2$, the sequence $(\alpha\theta_{p,a}\, \tmod 1)_{p \to \infty}$ is uniformly distributed in the interval $[0,1]$ (see Remark \ref{dilated-Weyl}).  This motivates us to explore analogues of central limit theorems for the sums $\sum_{p \leq x} h(\alpha\theta_{p,a})$ for suitable test functions.  We now adapt Setup \ref{spaces-rvs} in the context of Kloosterman angles.  As we will see, the setup and structural properties of Kloosterman angles are very similar to those of Frobenius angles for families of elliptic curves.

\begin{setup}\label{Spaces-rvs-Kloosterman}
For $N \geq 1$, let 
$\displaystyle \mathcal F_N = \{ a:\,1 \leq a \leq N\}$, be as in Theorem \ref{Kloosterman-multiples-smooth}, regarded as a finite probability space with the uniform measure.  
Let $\mathcal A$ denote the set of all primes, so that $A(x) = \pi(x)$.
For each $p \in \mathcal A$, 
define the sets
$$\mathcal G_{p,N} := \{1 \leq a \leq N:\,p \mid a\} \text{ and }\mathcal H_{p,N} = \mathcal F_N \setminus  \mathcal G_{p,N}.$$
To each positive integer $\alpha$ and $a \in \mathcal F_N$, we associate the sequence $\{\alpha\theta_{p,a}\,\tmod\,1\}_{p \to \infty}$.

For $m \in \Z$, and $\alpha \in \N$, define
\begin{align}\label{cmp-Kl}
\widetilde{c}_{m,p,N}(\alpha) := \frac{1}{|\mathcal F_N|} \sum_{a \in \mathcal H_{p,N}} \cos (2\pi m  \alpha\theta_{p,a}),
\end{align}
and
\begin{align*}
d_m(\alpha) := \lim_{x \to \infty}\frac{1}{\pi(x)}\sum_{p \leq x} \frac{1}{|\mathcal F_N|} \sum_{a \in \mathcal F_N} \cos (2\pi  m\alpha \theta_{p,a}).
\end{align*}
Let $h \in \mathcal C(\R)$ and $\mu_{\alpha}(t) := \sum_{m \in \Z} d_m(\alpha)e(-mt)$.  In analogy with Setup \ref{spaces-rvs-elliptic}, we have:
\begin{align*}
    h(\alpha \theta_{p,a}) &= \textstyle\sum_{m \in \Z} \widehat{h}(m) e(m\alpha \theta_{p,a}),\\
    h_M(\alpha\theta_{p,a}) &= \textstyle\sum_{|m| \leq M}\widehat{h}(m) e(m\alpha \theta_{p,a}),\\
    Z(\alpha \theta_{p,a}) &= h(\alpha \theta_{p,a}) - \textstyle\int_0^1 h(t) \mu_{\alpha}(t) dt  = h(\alpha \theta_{p,a}) - \sum\nolimits_{m \in \Z}\widehat{h}(m) d_m(\alpha),\\
    Z_M(\alpha \theta_{p,a}) &= h_M(\alpha \theta_{p,a}) - \textstyle\int_0^1 h_M(t) \mu_{\alpha}(t) dt  = h_M(\alpha \theta_{p,a}) - \sum\nolimits_{|m| \leq M}\widehat{h}(m) d_m(\alpha),\\
    W(\alpha \theta_{p,a}) &\coloneqq \chi_{\mathcal H_{p,N}}(a)Z(\alpha \theta_{p,a}),\\
    W_M(\alpha \theta_{p,a}) &\coloneqq \chi_{\mathcal H_{p,N}}(a)Z_M(\alpha \theta_{p,a}).
\end{align*}

\end{setup}

With the above setup, we prove Theorem \ref{Kloosterman-multiples-smooth} as another application of Theorem \ref{CLT general smooth variant}.  We first note an estimate towards verifying the hypothesis of \eqref{thin-sets}.

\begin{lemma}\label{thin-set-sum-Kl}
For each $a \in \mathcal F_N$, define $\omega_x(a) = \#\{p \leq x:\,p \mid a\}$.  Then, for every integer $u \geq 1$,
$$ \frac{1}{N} \sum_{a=1}^N \omega_x(a)^u \ll_u (\log \log x)^u + \frac{ \pi(x)^u}{N}.$$

\end{lemma}

\begin{proof}
    Expand $\omega_x(a)^u$ as a sum over $u$-tuples of primes, not necessarily distinct.  The condition $p_i \mid a$ for each $i$ means that lcm$(p_1,p_2,\dots,p_u) \mid a$, and the density of such $a$ in $\mathcal F_N$ is 
    $$\frac{1}{\text{lcm}(p_1,p_2,\dots,p_u)} + \O \left(\frac{1}{N}\right).$$
    Thus,
    \begin{align*}
        \frac{1}{N}\sum_{a=1}^N \omega_x(a)^u & =\sum_{p_1,p_2,\dots,p_u \leq x}\left( \frac{1}{\text{lcm}(p_1,p_2,\dots,p_u)} + \O \left(\frac{1}{N}\right)\right)\\
        &\ll_u\sum_{j=1}^u \left(\sum_{p_1 < p_2 < \dots <p_j \leq x} \frac{1}{p_1p_2\dots p_j} + \O \left(\frac{1}{N}\right)\right)\\
        &\ll_u(\log \log x)^u + \frac{ \pi(x)^u}{N}.
    \end{align*}
\end{proof}

The following proposition presents the values of $c_{m,p,N}(\alpha)$ and $d_m(\alpha)$. 

\begin{prop}\label{Weyl limits Kloosterman}
We follow the notation of Setup \ref{Spaces-rvs-Kloosterman}.  For a prime $p$, and an integer $m$,
\begin{align}\label{cmp Kloosterman}
 \widetilde{c}_{m,p,N}(\alpha) = \frac{1}{p}\sum_{a = 1}^{p-1} \cos (2 \pi m\alpha \theta_{p,a}) + \begin{cases}
 \displaystyle  \O\left(\frac{1}{N}\right) &\text{ if }m = 0,\\
\displaystyle  \O\left(\frac{p}{N}\right)& \text{ if }|m| \geq 1.
\end{cases} 
\end{align}
Let $N = N(x)$ such that 
$$ \frac{x}{N} \to 0 \text{ as }x \to \infty.$$
Then, for an integer $m$,
\begin{align}\label{dm Kloosterman}
d_m(\alpha) = \begin{cases}
\phantom{-}1 &\text{ if }m = 0\\
    -\frac{1}{2} &\text{ if }|m| = 1 = \alpha\\
    \phantom{-}0 &\text{ if }|m| \geq 2 \text{ or }\alpha \geq 2.
    \end{cases}
\end{align}
\end{prop}
\begin{proof}
If $m = 0$, $\displaystyle \widetilde{c}_{m,p,N}(\alpha) = 1 - \frac{1}{N}\left\lfloor\frac{N}{p}\right\rfloor = 1 - \frac{1}{p} + \O\left(\frac{1}{N}\right)$.  

Let $|m| \geq 1$. The term $\cos (2\pi  m \alpha\theta_{p,a})$ is periodic with respect to $a$ modulo $p$.

Arranging $\mathcal{H}_{p,N}$ into $\lfloor N/p\rfloor$ many complete set of coprime residue classes modulo $p$ and using the trivial bound $|\cos \theta| \leq 1 $ for the remaining residue classes, we get
 \begin{align*}
     \widetilde{c}_{m,p,N}(\alpha) =&\frac{1}{N}\sum_{a \in \mathcal H_{p,N}}\cos (2 \pi m\alpha \theta_{p,a})\\
     = &\frac{1}{N}\left\lfloor\frac{N}{p}\right\rfloor\sum_{a = 1}^{p-1} \cos (2 \pi m\alpha \theta_{p,a}) + \frac{1}{N}\sum_{\substack {\lfloor N/p \rfloor p < a \leq N \\p \nmid a}} \cos (2 \pi m\alpha \theta_{p,a})\\
     =&\frac{1}{p}\sum_{a = 1}^{p-1} \cos (2 \pi m\alpha \theta_{p,a}) + \O\left(\frac{p}{N}\right).
     \end{align*}
This proves \eqref{cmp Kloosterman}.  

For each prime $p$, writing $\mathcal F_N = \mathcal H_{p,N} \cup \mathcal G_{p,N}$ and applying Lemma \ref{thin-set-sum-Kl}, we get 
\begin{align}\label{dm-Kl}
d_m(\alpha) = \lim_{x \to \infty}\frac{1}{\pi(x)}\sum_{p \leq x}\widetilde{c}_{m,p,N}(\alpha) .
\end{align}
We now apply the identity
$$\frac{1}{p}\sum_{a = 1}^{p-1}\cos (2 \pi m \alpha \theta_{p,a}) = \frac{1}{2p}\sum_{a=1}^{p-1} \left(U_{2|m|\alpha} (\cos (\pi \theta_{p,a})) - U_{2|m|\alpha-2} (\cos (\pi \theta_{p,a}))\right) \text{ if }m \neq 0.$$
 By \eqref{cmp Kloosterman} and the Katz inequality applied to non-zero frequencies,
 we have
 \begin{align}
\label{cmp-Kl-bounds}\widetilde{c}_{m,p,N}(\alpha) &= \begin{cases}
    \displaystyle 1 - \frac{1}{p} + \O\left(\frac{1}{N}\right)&\text{ if }m = 0,\\
    \displaystyle -\frac{1}{2} + \O\left(\frac{1}{\sqrt{p}}\right) + \O\left(\frac{p}{N}\right)&\text{ if }|m| = \alpha = 1,\\
\displaystyle \O\left(\frac{|m|\alpha}{\sqrt{p}}\right) +\O\left(\frac{p}{N}\right)&\text{ if }|m| \geq 2 \text{ or }\alpha \geq 2.
    \end{cases}
\end{align}
Summing \eqref{cmp-Kl-bounds} over $p\leq x$, using Lemma \ref{thin-set-sum-Kl} (with $u=1$) for the primes dividing $a$ and the estimate $\displaystyle \sum_{p \leq x} \frac{1}{\sqrt{p}} \ll \sqrt{\frac{\pi(x)}{\log x}}$, gives \eqref{dm Kloosterman}, since $x/N \to 0$.

\end{proof}

The following proposition is the exact analogue of Proposition \ref{Condition 1 elliptic} in the Kloosterman setting, and leads to the verification of \ref{Cond GC 1}.

\begin{prop}\label{Condition 1 Kloosterman}
Let $h$ be as defined in Theorem \ref{Kloosterman-multiples-smooth}.  Let $\displaystyle M = M(x) \asymp x^{\frac{1}{2(\kappa - 1)}}$. Then
$$\sup_{1 \leq a \leq N}\abs{\sum_{p \leq x} (W - W_M)(\alpha \theta_{p,a})} = \o\left(\sqrt{\pi(x)}\right).$$

 \end{prop}
 \begin{proof}
    We follow verbatim the method in the proof of Proposition \ref{Condition 1 elliptic} (replacing $\alpha\theta_E(p)$ by $\alpha\theta_{p,a}$).  
    \end{proof}
   Condition \ref{Cond GC 2} for Setup \ref{Spaces-rvs-Kloosterman} will follow from the Chinese remainder theorem, and
\ref{Cond GC 3} from the Katz inequality.  We now state a proposition towards \ref{Cond GC 2}. 
\begin{prop}\label{appl of crt}
    Let $p_1,p_2,\dots,p_u$ be $u$ distinct primes and let $q=p_1p_2\cdots p_u$.  For integers $m_1,m_2,\dots,m_u$, 
    $$\frac{1}{N} \sum_{\substack{a \in \mathcal F_N \\ (a,q) = 1}}\prod_{i=1}^u \cos (2\pi m_i \alpha \theta_{p_i,a}) = \prod_{i=1}^u \widetilde{c}_{m_i,p_i,N}(\alpha) + \O\left(\frac{q}{N}\right).$$
\end{prop}
\begin{proof}
Arranging $\mathcal F_N$ into $\lfloor N/q \rfloor$ complete residue classes modulo $q$ and applying trivial bounds for the remaining residue classes $a$ lying in the interval $[q\lfloor N/q \rfloor +1,N]$, 
    \begin{align}\label{cmp-red}
\frac{1}{N} \sum_{\substack{a \in \mathcal F_N \\ (a,q) = 1}}\prod_{i=1}^u \cos (2\pi m_i \alpha \theta_{p_i,a})  =& \frac{1}{N}\left\lfloor\frac{N}{q} \right \rfloor\sum_{\substack{a =1 \\ (a,q) = 1}}^q\prod_{i=1}^u \cos (2\pi m_i \alpha \theta_{p_i,a}) + \O\left(\frac{q}{N}\right)\notag\\
=&\prod_{i=1}^u\frac{1}{p_i}\sum_{a_i =1}^{p_i-1}\cos (2\pi m_i \alpha \theta_{p_i,a_i}) + \O\left(\frac{q}{N}\right),
\end{align}
where the last equality follows by the Chinese remainder theorem. By \eqref{cmp Kloosterman}, we have, for each prime $p_i$,
\begin{align*}
&\frac{1}{p_i}\sum_{a_i =1}^{p_i-1}\cos (2\pi m_i \alpha \theta_{p_i,a_i}) - \widetilde{c}_{m_i,p_i,N}(\alpha)= \begin{cases}
    \displaystyle  \O\left(\frac{1}{N}\right)&\text{ if }m_i = 0\\
    \displaystyle  \O\left(\frac{p_i}{N}\right) &\text{ if }|m_i| \geq 1.
\end{cases}
\end{align*}
For each $i$,
$$\abs{\frac{1}{p_i}\sum_{a_i =1}^{p_i-1}\cos (2\pi m_i \alpha \theta_{p_i,a_i})} \leq 1 \text{ and }\abs{\widetilde{c}_{m_i,p_i,N}} \leq 1,$$
as both are averages of cosines over sets of density at most 1.  Thus, by Lemma \ref{product-estimates} \ref{product-sum},
\begin{align*}
\abs{\prod_{i=1}^u\frac{1}{p_i}\sum_{a_i =1}^{p_i-1}\cos (2\pi m_i \alpha \theta_{p_i,a_i}) - \prod_{i=1}^u \widetilde{c}_{m_i,p_i,N}(\alpha)}
&\leq \sum_{i=1}^u \abs{\frac{1}{p_i}\sum_{a_i =1}^{p_i-1}\cos (2\pi m_i \alpha \theta_{p_i,a_i}) - \widetilde{c}_{m_i,p_i,N}(\alpha)} \\
&\ll \frac{1}{N}\sum_{i=1}^u p_i \ll \frac{q}{N}.
\end{align*}
Combining this with \eqref{cmp-red} and applying the triangle inequality gives the proposition.

\end{proof}

\subsection{Proof of Theorem \ref{Kloosterman-multiples-smooth}}\label{sec: Kloosterman CLT proof}
We apply Theorem \ref{CLT general smooth variant} to Setup \ref{Spaces-rvs-Kloosterman}.  By Lemma \ref{h in setup5}, $h$ satisfies the hypotheses of Setup \ref{spaces-rvs}.

Since $N = N(x)$ satisfies $\displaystyle \frac{\log N}{\log x} \to \infty,$
 we get that
 \begin{align}\label{growth of N}
 \lim_{x \to \infty} \frac{x^{\beta}}{N} = 0  \text{ for any }\beta >0.
 \end{align}
 By Lemma \ref{thin-set-sum-Kl} and the above bound,
 \begin{equation}\label{thin set Kloosterman}
    \left \langle\left(\sum_{ p \leq x }\chi_{\mathcal G_{p,N}}\right)^u \right \rangle_k= \o\left(\pi(x)^{u/2}\right), 
 \end{equation}
 and therefore, the hypothesis of \eqref{thin-sets} is satisfied.
  Assume that $M(x) \asymp x^{1/(2(\kappa-1))}$. 

Condition \ref{Cond GC 1} is an immediate consequence of Proposition \ref{Condition 1 Kloosterman}.

To verify condition \ref{Cond GC 2}, we apply Proposition \ref{appl of crt} and get
\begin{align}\label{GC-2-estimate-Kl}
&\sum_{\underline{p}}\sum_{\underline{m}}\abs{\left(\frac{1}{N} \sum_{\substack{a \in \mathcal F_N \\ (a,q) = 1}}\prod_{i=1}^u \cos (2\pi m_i \alpha \theta_{p_i,a})\right)  - \prod_{i=1}^u  \widetilde{c}_{m_i,p_i,N}(\alpha)} \\
\ll& \sum_{\underline{p}}\sum_{\underline{m}} \frac{p_1p_2\dots p_u}{N} \ll \frac{M^ux^u\pi(x)^u}{N}\notag = \o(\pi(x)^{u/2}),\notag
\end{align}
 by the growth condition in \eqref{growth of N}.  This verifies \ref{Cond GC 2} as $u \leq r$.

To verify condition \ref{Cond GC 3}, apply equations \eqref{dm Kloosterman} and \eqref{cmp-Kl-bounds}.  Thus,
\begin{align}\label{GC3-inner-Kl}
    \sum_{|m_i| \leq Mr_i} \abs{ \widehat{Z_M^{r_i}}(m_i) (\widetilde{c}_{m_i,p_i,N}(\alpha)-d_{m_i}(\alpha))} \ll_{\alpha} \sum_{|m_i| \leq Mr_i} \abs{ \widehat{Z_M^{r_i}}(m_i)}\left(\frac{|m_i| + 1}{\sqrt{p_i}} + \frac{p_i}{N}\right).
\end{align}
As in the proof of Theorem \ref{Elliptic-curves-smooth}, applying \eqref{mZmpowers} and Lemma \ref{product-estimates}\ref{sum-product}, we have
\begin{align}\label{GC3-estimate-Kl}
&\sum_{\underline{p}}\prod_{i=1}^u \sum_{|m_i| \leq Mr_i} \abs{ \widehat{Z_M^{r_i}}(m_i) (\widetilde{c}_{m_i,p_i,N}(\alpha)-d_{m_i}(\alpha))}\\
\ll_u{}&\sum_{\underline{p}}\prod_{i=1}^u \left(\frac{ 1}{\sqrt{p_i}} + \frac{p_i}{N}\right)\notag\\
\ll_u{}& \left(\frac{\pi(x)}{\log x}\right)^{u/2} + \ \frac{x^u\pi(x)^u}{N}\notag\\ = \o{}&\left(\pi(x)^{u/2}\right)\notag,
\end{align}
which verifies \ref{Cond GC 3} as $u \leq r$.
This proves Theorem \ref{Kloosterman-multiples-smooth}.

\begin{cor}\label{Nagoshi-type-Kloosterman}
Fix an integer $j \geq 1$, and let $N=N(x)$ satisfy $\log N/\log x\to\infty$.  Then, as $a$ ranges over $1 \leq a \leq N$, the moments of
$$\frac{1}{\sqrt{\pi(x)}}\sum_{p\leq x}U_j\left(\cos(\pi\theta_{p,a})\right) = \frac{1}{\sqrt{\pi(x)}}\sum_{p\leq x}U_j\left(\frac{K(p,a)}{2\sqrt{p}}\right)$$
converge to the Gaussian moments.  For $j=1$, this is the statistic
$$ \frac{1}{\sqrt{\pi(x)}}\sum_{p\leq x}\frac{K(p,a)}{\sqrt{p}}.$$
\end{cor}

\begin{proof}
     We apply Theorem \ref{2-periodic} with $\alpha = 1$ and $H(t) = U_j(\cos (\pi t))$, $\mathcal A$ the set of all primes, and the exceptional sets $\mathcal G_{p,N}$ of Setup \ref{Spaces-rvs-Kloosterman}.   The thinness condition was verified in the proof of Theorem \ref{Kloosterman-multiples-smooth}. Since $H$ is a trigonometric polynomial of degree $j$ in $\pi t$, condition \ref{Cond GC 1} holds trivially with $M=j$.   

    For every $m\neq0$, $\cos(\pi m\theta)=\frac{1}{2}\left(U_{|m|}-U_{|m|-2}\right)(\cos\pi\theta)$, with $U_{-1}:=0$.
     
     Since the Katz inequality \eqref{Katz inequality} holds for all $m \geq 1$, the argument of  Proposition \ref{Weyl limits Kloosterman} gives, for both even and odd frequencies,  
    \begin{equation}\label{Cmp1-Kl}
C_{m,p,N}(1) = D_m(1)+\O\left(\frac{|m|+ 1}{\sqrt{p}}\right) + \O\left(\frac{p}{N}\right), \quad D_m(1) = \delta_{m,0} -\frac{1}{2}\delta_{|m|,2}.
\end{equation}

 Proposition \ref{appl of crt} only uses the periodicity of $\theta_{p,a}$ modulo $p$ and the Chinese Remainder Theorem.  Thus, we get, for $q = p_1p_2\dots p_u$,
\begin{equation}\label{approx-ind-odd-Kl}
    \frac{1}{N} \sum_{\substack{a \in \mathcal F_N \\ (a,q) = 1}}\prod_{i=1}^u \cos (\pi m_i  \theta_{p_i,a}) = \prod_{i=1}^u {C}_{m_i,p_i,N}(1) + \O\left(\frac{q}{N}\right).
    \end{equation}
   Conditions \ref{Cond GC 2} and \ref{Cond GC 3} now follow exactly as in the proof of Theorem \ref{Kloosterman-multiples-smooth}. 
   
   Finally, $\nu(t) = \sum_{m \in \Z}D_m(1)\cos \pi mt = 2\sin^2\pi t$ , and by the orthogonality relation \eqref{eq:mean zero}, $\int_0^1 H(t) \nu(t) dt= 0$, and 
   $V_H(\nu) = \int_0^1 H(t)^2 \nu(t) dt = 1.$  Since $A(x) = \pi(x)$, Theorem \ref{2-periodic} gives the corollary.
\end{proof}

\section{Sequences arising from Modular forms}\label{sec: modular forms}

In this section, we prove Theorem \ref{Hecke-multiples-smooth}. The structure is similar to that of Sections \ref{sec: Elliptic curves} and \ref{sec: Kloosterman}.  The main tool in this section is the Eichler–Selberg trace formula (Proposition \ref{trace formula}), whose uniformity in the weight and the level allows both to vary. The proof of Theorem \ref{Hecke-multiples-smooth} is completed in Section \ref{sec: modular CLT proof}, where the cases \ref{polynomial decay}–\ref{compact-support}  are treated, and the section ends with Nagoshi-type statistics (Corollary \ref{Nagoshi-revisited}) and a comparison with earlier work.

Let $k$ be an even, positive integer, and let $S_k(\Gamma_0(N))$ denote the space of modular cusp forms of weight $k$ with respect to $\Gamma_0(N)$.  Let $S^{\text{new}}(N,k)$ denote the space of primitive cusp forms of weight $k$ with respect to $\Gamma_0(N)$, that is, the orthogonal complement, with respect to the Petersson inner product, of the subspace of oldforms in the space $S_k(\Gamma_0(N))$.  Let $\f_{N,k}$ be the set of normalized Hecke newforms in $S^{\text{new}}(N,k)$; by Atkin--Lehner theory, it forms an orthogonal basis of $S^{\text{new}}(N,k)$.  Each $f \in \f_{N,k}$ has a   Fourier expansion 
\begin{align*}
    f(z) = \sum_{n=1}^{\infty} n^{(k-1)/2}a_f(n)e(nz),\,z \in \H,
\end{align*}
where $a_f(1) = 1$, and the normalized Fourier coefficients $a_f(n)$ satisfy the property 
\begin{equation}\label{prod formula}
           a_f(n)a_f(m)=\sum_{\genfrac{}{}{0pt}{}{d|(m,n)}{(d,N)=1}}a_f\left(\frac{mn}{d^2}\right).
\end{equation}
In particular, 
\begin{align}\label{mult-FC-f}
a_f(mn)=a_f(m)a_f(n) \text{ whenever }(m,n) = 1.
\end{align}
For a prime $p \nmid N$ and for integers $m_1,\,m_2 \geq 0$, \eqref{prod formula} gives
\begin{equation}\label{prod of coeff prime power}
    a_f(p^{m_1})a_f(p^{m_2})=\sum_{l=0}^{\min(m_1,m_2)}a_f\left(p^{m_1+m_2-2l}\right),
\end{equation}
while for a prime $p \mid N$, only $d=1$ contributes to \eqref{prod formula}, so that $\displaystyle a_f(p^m)=a_f(p)^m$ for all $m \geq 0$.  Moreover, if $p \mid N$, then $a_f(p)^2=1/p$ when $p\,\|\,N$, and $a_f(p)=0$ when $p^2\mid N$ (see \cite{AL}, \cite{Winnie-Li}).  Thus, $\displaystyle |a_f(p)|\leq \frac{1}{\sqrt{p}}$.

By the Ramanujan--Petersson conjecture, proved by Deligne, we have 
\begin{equation}\label{Ramanujan--Deligne}
    |a_f(n)| \leq d(n) \text{ for all }(n,N) = 1.
\end{equation}
In fact, by the above discussion, this bound holds for all $n \geq 1$.  
For a prime $p \nmid N$,  $a_f(p)$ is a real number and $|a_f(p)| \leq 2$.  Thus, we write 
\begin{equation}\label{Hecke angle}
    a_f(p)=2\cos(\pi\thtp)\text{ for a unique }\thtp \in [0,1].
\end{equation}
 Taking $m_2=1$ in \eqref{prod of coeff prime power} gives $a_f(p^{m+1})=a_f(p)\,a_f(p^m)-a_f(p^{m-1})$ for $m\geq 1$. By \eqref{recurrence of Chebyshev poly}, this is the recurrence satisfied by $U_m(\cos (\pi \thtp))$, and the initial values agree: $a_f(1) = 1 = U_0(\cos (\pi \thtp))$, and $a_f(p) = 2\cos \pi \thtp = U_1(\cos (\pi \thtp))$.  Hence, for $p \nmid N$ and $m\geq 0$
\begin{align*}
a_f(p^m)=U_m(\cos(\pi\thtp)).
\end{align*}

We recall the following estimates of the Eichler--Selberg trace formula for the traces of Hecke operators acting on $S^{\text{new}}(N,k)$.  
\begin{prop}\label{trace formula}
Let $k$ be a positive, even integer and let $N$ be a positive integer such that $\mathcal F_{N,k}$ is non-empty. Let $n \geq 1$ be coprime to $N$. Then, for any real number $d' > 1$,
    \begin{equation}
        \frac{1}{|\f_{N,k}|}\sum_{f\in\f_{N,k}}a_f(n)=
    \begin{cases}
        \frac{1}{\sqrt{n}}+\O_{d'}\left(\frac{8^{\omega(N)}n^{d'}}{kN}\right)&\text{if $n$ is a square},\\
        \O_{d'}\left(\frac{8^{\omega(N)}n^{{d'}}}{kN}\right)&\text{otherwise},
    \end{cases}
    \end{equation}
    where the implied constant depends only on $d'$.  
    \end{prop}
\begin{proof}
For $n = 1$, the average on the left hand side is exactly 1.  For a proof of the case when $n> 1$, we refer to \cite[Proposition 15]{Mahajan-Sinha}.
\end{proof}

By the discussion in Section \ref{sec: ed-Weyl}, the Sato--Tate distribution law for Hecke angles \cite{BGHT} implies that  for a non-CM $f$, the sequence $(\alpha \thtp \,\tmod 1)_{p \to \infty}$ is $\mu_{\alpha}$-equidistributed in $[0,1]$, with $\mu_1$ the Sato--Tate density and $\mu_{\alpha} = 1$ for $\alpha \geq 2$.  We now explore the smooth statistics for these sequences for each $\alpha \geq 1$.

\begin{setup}\label{spaces-rvs-modular}
For positive integers $k$ and $N$ with $k$ even, let $\mathcal F_{N,k}$ be as in Theorem \ref{Hecke-multiples-smooth}, regarded as a finite probability space with the uniform measure; it is nonempty for $x$ sufficiently large, since $kN\to\infty$.  Let $\mathcal A$ denote the set of primes, so that $A(x) = \pi(x)$.  For each prime $p$, define the sets
$$\mathcal G_{p,N,k} = \begin{cases} 
\mathcal F_{N,k} &\text{ if }p \mid N\\
\emptyset &\text{ if }p \nmid N
\end{cases} \text{ and }\mathcal H_{p,N,k} = \mathcal F_{N,k} \setminus \mathcal G_{p,N,k}.$$
Thus,
$$\frac{1}{\abs{\mathcal F_{N,k}}}\sum_{f \in \mathcal F_{N,k}}\left(\sum_{p \leq x} \chi_{ G_{p,N,k}}(f)\right)^u =  \left(\sum_{\substack{ p\leq x \\ p\mid N}}1\right)^u \leq \omega(N)^u.$$
If we choose $N = N(x)$ such that $\omega(N) = \o\left(\sqrt{\pi(x)}\right)$, then
\begin{align}\label{thin-sets-modular}
\frac{1}{\abs{\mathcal F_{N,k}}}\sum_{f \in \mathcal F_{N,k}}\left(\sum_{p \leq x} \chi_{ \mathcal G_{p,N,k}}(f)\right)^u  = \o\left(\pi(x)^{u/2}\right).\end{align}

To each positive integer $\alpha$ and $f \in \mathcal F_{N,k}$, we associate the sequence $(\alpha \theta_f(p) \, \tmod 1)_{p \to \infty}$.

Adapting the notation of Theorems \ref{CLT general smooth} and  \ref{CLT general smooth variant}, we have, 
\begin{align}\label{WM-modular}
(W - W_M)(\alpha\thtp) = \begin{cases}
0 &\text{ if }p \mid N,\\
    \sum_{|m| > M} \widehat{Z}(m) \cos (2\pi m\alpha\thtp) &\text{ if }p \nmid N.
    \end{cases}    
\end{align}

For each $m \in \Z$ and $\alpha \in \N$,
define
\begin{align}
\widetilde{c}_{m,p,N,k}(\alpha) :={}& \frac{1}{|\mathcal F_{N,k}|}\sum_{f \in \mathcal H_{p,N,k}} \cos 2\pi m \alpha\thtp \notag\\
\label{cmp-modular}
={}& \begin{cases}
    0 &\text{ if }p \mid N\\
    \displaystyle \frac{1}{|\mathcal F_{N,k}|}\sum_{f \in \mathcal F_{N,k}}\cos (2\pi m \alpha \thtp) &\text{ if }p \nmid N,
\end{cases}\\
\intertext{and}
\label{dm-modular}
d_m(\alpha) :={}&  \lim_{x \to \infty} \frac{1}{\pi(x)} \sum_{p \leq x} \frac{1}{|\mathcal F_{N,k}|}\sum_{f \in \mathcal F_{N,k}}\cos (2\pi m \alpha \thtp) .
\end{align}

\end{setup}

We now state some propositions in order to verify \ref{Cond GC 1} - \ref{Cond GC 3} of Theorem \ref{CLT general smooth variant}.  We first evaluate the family averages $\tilde{c}_{m,p,N,k}$ and the Weyl limits $d_m$ for the setup of Hecke angles. 

\begin{prop}\label{Weyl-limits-modular}
Let $\widetilde{c}_{m,p,N,k}(\alpha)$ and $d_m(\alpha)$ be as defined in \eqref{cmp-modular} and \eqref{dm-modular} respectively.  
Then, for $p \nmid N$,
\begin{align}\label{cmp-modular-calculated}
\widetilde{c}_{m,p,N,k} (\alpha)  = 
\begin{cases}
    \displaystyle \frac{1}{2}\left(\frac{1}{p^{|m|\alpha}} - \frac{1}{p^{|m|\alpha - 1}}\right) + \O_{d'}\left(\frac{8^{\omega(N)}p^{2|m|\alpha d'}}{kN}\right) &\text{ if }|m| \geq 1\\
    1 &\text{ if }m = 0.
\end{cases}
\end{align}
Let $k = k(x)$ and $N = N(x)$ with $k$ even such that 
\begin{align}\label{growth-conditions-N-and-k}
\frac{\displaystyle\log (kN/8^{\omega(N)})}{\log x} \to \infty \text{ and } \omega(N) = \o\left(\pi(x)\right).
\end{align}
Then,
\begin{align}\label{dm-modular-calculated}
d_m(\alpha) = \begin{cases}
 \phantom{-}1 &\text{ if }m = 0\\
-\frac{1}{2}&\text{ if }|m| = 1 = \alpha\\
 \phantom{-}0 &\text{ if }|m|\geq 2 \text{ or }\alpha \geq 2.
\end{cases}
\end{align}
\end{prop}

\begin{proof}
Let $p$ be a prime such that $(p,N) = 1$.  Then, for $m \neq 0$,
$$\cos (2\pi m \alpha \thtp) = \frac{1}{2}\left(a_f(p^{2|m|\alpha}) - a_f(p^{2|m|\alpha - 2})\right).$$
By Proposition \ref{trace formula}, for any $d' >1$ and $m \neq 0$,
\begin{align*}
&\frac{1}{|\mathcal F_{N,k}|}\sum_{f \in \mathcal F_{N,k}} \Big(a_f(p^{2|m|\alpha}) - a_f(p^{2|m|\alpha - 2})\Big)
=\frac{1}{p^{|m|\alpha}} - \frac{1}{p^{|m|\alpha - 1}} + \O_{d'}\left(\frac{8^{\omega(N)}p^{2|m|\alpha d'}}{kN}\right).
\end{align*}
This proves \eqref{cmp-modular-calculated}.

Fix $d' > 1$.  By virtue of \eqref{thin-sets-modular} for $u=1$ and the trivial bound $\abs{\cos (2\pi m \alpha \thtp)} \leq 1$ for $p \mid N$, we have
\begin{align*} 
&\frac{1}{\pi(x)} \sum_{p \leq x} \frac{1}{|\mathcal F_{N,k}|}\sum_{f \in \mathcal F_{N,k}}\cos (2\pi m \alpha \thtp)\\
={}& 
\frac{1}{\pi(x)} \sum_{\substack {p \leq x \\(p,N) = 1}} \frac{1}{|\mathcal F_{N,k}|}\sum_{f \in \mathcal F_{N,k}}\cos (2\pi m \alpha \thtp) + \O\left(\frac{\omega( N)}{\pi(x)}\right)\\
={}&\frac{1}{\pi(x)} \sum_{\substack {p \leq x \\(p,N) = 1}}\frac{1}{2}\left(\frac{1}{p^{|m|\alpha}} - \frac{1}{p^{|m|\alpha - 1}} \right)+ \O\left(\frac{8^{\omega(N)}p^{2|m|\alpha d'}}{kN}\right) + \O\left(\frac{\omega( N)}{\pi(x)}\right)\\
={}& -\frac{1}{2}\delta_{|m|\alpha,1} + \O\left(\frac{\log \log x}{\pi(x)}\right) + \O\left(\frac{8^{\omega(N)}x^{2|m|\alpha d'}}{kN}\right) + \O\left(\frac{\omega( N)}{\pi(x)}\right).
  \end{align*}
 Applying the growth condition in \eqref{growth-conditions-N-and-k}, we have
 $$\lim_{x \to \infty} \frac{8^{\omega(N)}x^{2|m|\alpha d'}}{kN} = 0 \text{ for every }|m| \geq 1,$$
 and 
 $$\lim_{x \to \infty}\frac{\omega( N)}{\pi(x)} = 0.$$
 Thus, we derive \eqref{dm-modular-calculated}.
\end{proof}

Next, we state a proposition that leads to \ref{Cond GC 1}.
\begin{prop}\label{Condition 1 modular}
Let $h$ be as defined in Theorem \ref{Hecke-multiples-smooth}.  Suppose $M = M(x)$ satisfies 
\begin{align}\label{FC-bound-modular}
\sum_{m > M} \abs{\widehat{h}(m)} = \o\left(\frac{1}{\sqrt{\pi(x)}}\right).
\end{align}
 Then,
 $$\sup_{f \in \mathcal F_{N,k}}\abs{\sum_{p \leq x }(W - W_M)(\alpha\thtp)} = \o\left(\sqrt{\pi(x)}\right).$$

\end{prop}
\begin{proof}
We have, for any $f$,  $\displaystyle \biggl|\sum_{p \leq x }(W - W_M)(\alpha\thtp)\biggr| \leq \sum_{p \leq x} \sum_{|m| > M} \abs{\widehat{h}(m)} = \o\left(\sqrt{\pi(x)}\right).$

\end{proof}

The following proposition will help us to verify condition \ref{Cond GC 2}.   
\begin{prop}\label{Towards GC 2}
 Let $p_1,p_2,\dots, p_u$ be distinct primes, and let $m_1,m_2,\dots, m_u$ be integers.  For any fixed $d'>1$, we have 
           \begin{align*} 
\abs{\prod_{i=1}^u \widetilde{c}_{m_i,p_i,N,k}(\alpha)-\left \langle \prod_{i=1}^u\begin{cases}\cos (2\pi m_i\alpha\thtpi) &\text{ if }f \in \mathcal H_{p_i,N,k}\\
0 &\text{ if }f \in \mathcal G_{p_i,N,k}
\end{cases} \right \rangle_k } =  \O_u\left(\frac{8^{\omega(N)}}{kN}\prod_{i=1}^u p_i^{2|m_i|\alpha d'}\right).
            \end{align*}

\end{prop}
\begin{proof}
If $p_i \mid N$ for some $i$, then there is nothing to prove as the left hand side is zero.  Let $( p_1p_2\dots p_u,N) = 1$.  Then,
 \begin{align}\label{avg-prod-mod}
  &\frac{1}{|\mathcal F_{N,k}|} \sum_{f \in \mathcal F_{N,k}} \prod_{i=1}^u \cos (2\pi m_i \alpha \thtpi)\\
  ={}&  \frac{1}{|\mathcal F_{N,k}|} \sum_{f \in \mathcal F_{N,k}} \left(\prod_{i=1}^u\begin{cases}
      1 &\text{ if }m_i = 0\\
      \frac{1}{2}\left(a_f\left( p_i^{2|m_i|\alpha}\right) - a_f\left( p_i^{2|m_i|\alpha-2}\right)\right) &\text{ if }|m_i| \geq 1\\
      \end{cases}\right)\notag\\
     ={}& \left(\prod_{i=1}^u \begin{cases}
      1 &\text{ if }m_i = 0\\
   \displaystyle\frac{1}{2}\left(\frac{1}{p_i^{|m_i|\alpha}} - \frac{1}{p_i^{|m_i|\alpha - 1} } \right) &\text{ if }|m_i| \geq 1   
 \end{cases}  \right) + \O_{d',u}\left(\frac{8^{\omega(N)}}{kN}\prod_{i=1}^u p_i^{2|m_i|\alpha d'}\right)\notag,
 \end{align}
where the last equality follows by the multiplicativity of $a_f(n)$ as well as the main term of Proposition \ref{trace formula} as functions of $n$. The error term also follows from Proposition \ref{trace formula}.  Let us denote
$$T_i = \begin{cases}
      1 &\text{ if }m_i = 0\\
   \displaystyle\frac{1}{2}\left(\frac{1}{p_i^{|m_i|\alpha}} - \frac{1}{p_i^{|m_i|\alpha - 1} } \right) &\text{ if }|m_i| \geq 1.  
 \end{cases}$$
By \eqref{cmp-modular-calculated},
$$\widetilde{c}_{m_i,p_i,N,k}(\alpha) = T_i + \mathcal E_i,$$
where $\displaystyle \abs{\mathcal E_i} \ll \frac{8^{\omega(N)}}{kN} p_i^{2|m_i|\alpha d'}$ (for $m_i = 0, \text{define } \mathcal E_i = 0$).  For each $i$,  $\widetilde{c}_{m_i,p_i,N,k}(\alpha)$ and $T_i$ have absolute value at most 1.  Thus, by Lemma \ref{product-estimates}\ref{product-sum},
\begin{align*}
\abs{\prod_{i=1}^u \widetilde{c}_{m_i,p_i,N,k}(\alpha) - \prod_{i=1}^u T_i}  &\leq \sum_{i=1}^u  \abs{\widetilde{c}_{m_i,p_i,N,k}(\alpha) - T_i}\\
&\leq\sum_{i=1}^u \abs{\mathcal E_i} \ll_{d',u}\sum_{i=1}^u \frac{8^{\omega(N)}}{kN} p_i^{2|m_i|\alpha d'}.
\end{align*}
Henceforth, let us fix $d'>1$.  Combining the above inequality with \eqref{avg-prod-mod} and applying the triangle inequality, we derive the proposition.
\end{proof}

\begin{remark}\label{Lau-Wang-Tang-Wang}
    Proposition \ref{Towards GC 2} is a quantitative form of the approximate independence of Hecke angles at distinct primes, in the vertical aspect. For Hecke eigenforms of level $1$, Lau and Wang \cite[Theorem 2]{LW} proved a quantitative joint equidistribution theorem for $(a_f(p_1),\dots,a_f(p_u))$ as $k\to\infty$, with respect to the product of the $p_n$-adic Plancherel measures (see \eqref{p-adic plancherel measure}), with an error term of order $u\log(p_1\cdots p_u)/\log k$ for boxes. Tang and Wang \cite[Theorems 1 and 2]{TANG2016295} studied the distributions of $a_f(p^{r})$ and of $a_f(p^{r})-a_f(p^{r-2})$ as $k\to\infty$, with quantitative versions. Since $a_f(p^{r})-a_f(p^{r-2})=2\cos(\pi r\theta_f(p))$, the latter is the vertical distribution of the dilated angles $r\theta_f(p)$, which appear in Theorem \ref{Hecke-multiples-smooth} (for $r = 2\alpha$ even) and in the 2-periodic extension (for odd $r$). In contrast with these results, which concern indicator functions of boxes and have logarithmic error terms, Proposition \ref{Towards GC 2} concerns the product $\prod_i\cos(2\pi m_i\alpha\theta_f(p_i))$, which are what the moment computations require; for these, the trace formula gives errors with a power saving in $kN$, uniformly in both the weight and the level.
\end{remark}

We now state a proposition towards verifying condition \ref{Cond GC 3}.
\begin{prop}
 \label{Towards GC 3} Let $r_1,r_2,\dots, r_u$ be positive integers, and denote $r := r_1 + r_2 + \dots +r_u$.  Suppose $\underline p$ runs over $u$-tuples of distinct primes $p_i\leq x$.  Then,
\begin{align*}
&\sum_{\underline{p}}\prod_{i=1}^u \sum_{|m_i| \leq Mr_i} \abs{\widetilde{c}_{m_i,p_i,N,k}(\alpha) - d_{m_i}(\alpha)} 
\ll_u(\omega(N) + \log\log x)^u+ \left(x^{2Mr\alpha d'}\pi(x)\frac{8^{\omega(N)}}{kN}\right)^u.
\end{align*}   
\end{prop}

\begin{proof}
 If $p_i \mid N$, then $\widetilde{c}_{m_i,p_i,N,k}(\alpha) = 0$ for all $m_i$ (since $\mathcal H_{p_i,N,k} = \emptyset$).  Thus, if $p_i \mid N$, then
 \begin{align}\label{p dividing N}
 \sum_{|m_i| \leq Mr_i} \abs{\widetilde{c}_{m_i,p_i,N,k}(\alpha) - d_{m_i}(\alpha)} = \begin{cases} 2 &\text{ if }\alpha = 1\\
 1 &\text{ if }\alpha \geq 2.
 \end{cases}
 \end{align}
 If $p_i \nmid N$, then
\begin{align*}
&\abs{\widetilde{c}_{m_i,p_i,N,k}(\alpha) - d_{m_i}(\alpha)} =\begin{cases}
     \displaystyle\frac{1}{2p_i} + \O\left(\frac{8^{\omega(N)}}{kN} p_i^{2 d'}\right)&\text{ if }|m_i|\alpha = 1\\
     \displaystyle \frac{p_i-1}{2p^{|m_i|\alpha}} + \O\left(\frac{8^{\omega(N)}}{kN} p_i^{2|m_i|\alpha d'}\right)&\text{ if }|m_i|\alpha > 1.
 \end{cases}
\end{align*}
Thus, for $p_i \nmid N$,
\begin{align}\label{p not dividing N}
\sum_{|m_i| \leq Mr_i} \abs{\widetilde{c}_{m_i,p_i,N,k}(\alpha) - d_{m_i}(\alpha)} 
    \ll \frac{1}{p_i} + \frac{8^{\omega(N)}}{kN}p_i^{2Mr_i\alpha d'}.
\end{align}
Summing over the primes, combining \eqref{p dividing N} and \eqref{p not dividing N},
\begin{align}\label{sum-product-cmp-mod}
\sum_{p \leq x} \left(\sum_{|m_i| \leq Mr_i} \abs{\widetilde{c}_{m_i,p,N,k}(\alpha) - d_{m_i}(\alpha)} \right)
\ll \omega(N) + \log \log x + \frac{8^{\omega(N)}}{kN}x^{2Mr\alpha d'}\pi(x).
\end{align}
Applying Lemma \ref{product-estimates}\ref{sum-product} to 
$\displaystyle F_i(p) = \sum_{|m_i| \leq Mr_i} \abs{\widetilde{c}_{m_i,p,N,k}(\alpha) - d_{m_i}(\alpha)}$ and using \eqref{sum-product-cmp-mod}, we have
\begin{align*}
\sum_{\underline{p}}\prod_{i=1}^u \sum_{|m_i| \leq Mr_i} \abs{\widetilde{c}_{m_i,p_i,N,k}(\alpha) - d_{m_i}(\alpha)} 
& \leq \prod_{i=1}^u\sum_{p\leq x} \left(\sum_{|m_i| \leq Mr_i} \abs{\widetilde{c}_{m_i,p,N,k}(\alpha) - d_{m_i}(\alpha)}\right)\\
& \ll \left(\omega(N) + \log \log x + \frac{8^{\omega(N)}}{kN}x^{2Mr\alpha d'}\pi(x)\right)^u.
\end{align*}   
By the elementary inequality
$(a+b)^u \leq 2^{u-1}(a^u + b^u)$ for $a,b \geq 0$, we have
\begin{align*}
\left(\omega(N) + \log \log x + \frac{8^{\omega(N)}}{kN}x^{2Mr\alpha d'}\pi(x)\right)^u \ll_u (\omega(N) + \log\log x)^u+ \left(x^{2Mr\alpha d'}\pi(x)\frac{8^{\omega(N)}}{kN}\right)^u.    
\end{align*}
\end{proof}

\subsection{Proof of Theorem \ref{Hecke-multiples-smooth}}\label{sec: modular CLT proof}

We apply Theorem \ref{CLT general smooth variant} to Setup \ref{spaces-rvs-modular}. By Lemma \ref{h in setup5}\ref{hm}, $h$ satisfies the hypotheses of Setup \ref{spaces-rvs}. The hypothesis of equation \eqref{thin-sets} has already been verified (see \eqref{thin-sets-modular}).  Choose $M = M(x)$ such that
\begin{equation}\label{growth conditions FC and families}
\sum_{m > M} \abs{\widehat{h}(m)} = \o\left(\frac{1}{\sqrt{\pi(x)}}\right) \text{ and }\frac{\log (kN/8^{\omega(N)})}{M\log x} \to \infty.
\end{equation}
By Proposition \ref{Condition 1 modular}, condition \ref{Cond GC 1} holds for every $l \geq 1$. 

We observe the following bound:
\begin{align}\label{xMC}
\text{ If }\frac{\log (kN/8^{\omega(N)})}{M\log x} \to \infty, \text{ then }
x^{MC} \ll kN/8^{\omega(N)} \text{ for every fixed }C >0.
\end{align}
To verify condition \ref{Cond GC 2}, we use Proposition \ref{Towards GC 2}.  Since $m_i \leq Mr_i$ for each $i$, and $r = r_1 + r_2 + \dots + r_u$, we deduce that
\begin{align}\label{GC-2-mod}
&\sum_{\underline{p}}\sum_{\underline{m}}\abs{\prod_{i=1}^u \widetilde{c}_{m_i,p_i,N,k}(\alpha)-\left \langle \prod_{i=1}^u\begin{cases}\cos (2\pi m_i\alpha\thtpi) &\text{ if }f \in \mathcal H_{p_i,N,k}\\
0 &\text{ if }f \in \mathcal G_{p_i,N,k}
\end{cases} \right \rangle_k }\\
\ll{}&\sum_{\underline{p}}\sum_{\underline{m}}\frac{8^{\omega(N)}}{kN}\prod_{i=1}^u p_i^{2|m_i|\alpha d'}
\ll\frac{8^{\omega(N)}}{kN}x^{2Mr\alpha d'}\pi(x)^u\notag\\
={}&\o\left(\pi(x)^{r/2}\right),\notag
\end{align}
where the last asymptotic follows by \eqref{xMC}.  This verifies \ref{Cond GC 2}. To verify condition \ref{Cond GC 3}, we invoke  Proposition \ref{Towards GC 3} and equation \eqref{xMC}.  Since
$\omega(N) = \o\left(\sqrt{\pi(x)}\right)$,  
we deduce that 
$$(\omega(N) + \log\log x)^u+ \left(x^{2Mr\alpha d'}\pi(x)\frac{8^{\omega(N)}}{kN}\right)^u = \o\left(\pi(x)^{r/2}\right).$$
Combining the inequality $\abs{\widehat{Z_M^{r_i}}(m)} \ll 1$ (Lemma \ref{properties of powers of Z_M}\ref{ZM-power-F}) with \eqref{xMC} and Proposition \ref{Towards GC 3}, 
\begin{align}\label{GC-3-mod}
&\sum_{\underline{p}}\prod_{i=1}^u \sum_{|m_i| \leq Mr_i} \abs{\widehat{Z_M^{r_i}}(m_i)(\widetilde{c}_{m_i,p_i,N,k}(\alpha) - d_{m_i}(\alpha))} \\
&\ll_u(\omega(N) + \log\log x)^u+ \left(x^{2Mr\alpha d'}\pi(x)\frac{8^{\omega(N)}}{kN}\right)^u\notag\\
&=\o\left((\pi(x))^{r/2}\right). \notag
\end{align} 
By Theorem \ref{CLT general smooth variant}, the moments of $\mathcal S_x$ converge to the Gaussian moments. 

We now discuss the cases \ref{polynomial decay}-\ref{compact-support} of Theorem \ref{Hecke-multiples-smooth}.  

\ref{polynomial decay} Suppose $\displaystyle \widehat{g}(t) \ll \frac{1}{(1 +|t|)^{\kappa}}$ for some $\kappa > 1$.  Choose $$M(x) = x^{\frac{1}{2\kappa - 2} + \epsilon}\text{ and  }\displaystyle \log (kN/8^{\omega(N)}) \geq x^{1/(\kappa - 1)}.$$
The conditions  of \eqref{growth conditions FC and families} hold since 
$$ \sum_{|m| > M} \abs{\widehat{h}(m)} \ll \frac{1}{M^{\kappa - 1}} \text{ and }\frac{\log (kN/8^{\omega(N)})}{M\log x} \geq  \frac{x^{1/(\kappa - 1)}}{x^{\frac{1}{2\kappa - 2} + \epsilon}\log x} \to \infty. $$

 \ref{Schwartz} Suppose $\displaystyle \widehat{g}(t) \ll \frac{1}{(1 +|t|)^B}$ for every $B > 0$.  Choose 
 $$M(x) = x^{\eta/2}\text{ for some } \eta >0\text{ and }\log (kN/8^{\omega(N)}) \geq x^{\eta}.$$
 Taking $B > 1 + \frac{1}{\eta}$, we get
 $$\sum_{|m| > M} \abs{\widehat{h}(m)} \ll \sum_{|m| > x^{\eta/2}}\frac{1}{m^B} \ll \frac{1}{x^{\frac{\eta}{2}(B-1)}} \ll \frac{1}{x^{1/2}}$$
and  
$$\frac{\log (kN/8^{\omega(N)})}{M\log x} = \frac{\log (kN/8^{\omega(N)})}{x^{\eta /2}\log x} \to \infty,$$
thus deriving the conditions in \eqref{growth conditions FC and families}.

 \ref{exp-decay} Suppose $\displaystyle \widehat{g}(t) \ll e^{-\lambda|t|^{\gamma}}$ for some $\lambda,\gamma >0.$
 Choose $$M = L\left(\frac{1}{\lambda}\log x\right)^{1/\gamma}  \text{ and suppose }\frac{\log (kN/8^{\omega(N)})}{(\log x)^{1 + 1/\gamma}} \to \infty.$$
Then, for $m > M$,  we have 
$$e^{-\frac{\lambda}{2}\left|\frac{m}{L}\right|^{\gamma}} < e^{-\frac{1}{2}\log x} = \o\left(\frac{1}{\sqrt{\pi(x)}}\right).$$
Since $\displaystyle \sum_{|m| > M}e^{-\frac{\lambda}{2}\left|\frac{m}{L}\right|^{\gamma}} < \infty$, we get
$$\sum_{|m| > M}e^{-\lambda\left|\frac{m}{L}\right|^{\gamma}} = \sum_{|m| > M}e^{-\frac{\lambda}{2}\left|\frac{m}{L}\right|^{\gamma}}e^{-\frac{\lambda}{2}\left|\frac{m}{L}\right|^{\gamma}} = \o\left(\frac{1}{\sqrt{\pi(x)}}\right).$$
We also have
$$\frac{\log (kN/8^{\omega(N)})}{L(1/\lambda)^{1/\gamma}(\log x)^{1 + 1/\gamma}} \to \infty,$$
by the family growth hypothesis in \ref{exp-decay}.  Therefore, the conditions of \eqref{growth conditions FC and families} hold.

 \ref{compact-support} Suppose the support of $\widehat{g}$ is contained in $[-K,K]$ for a fixed constant $K \geq 1$.  Then $\widehat{h}(m) = 0$ for $|m| > KL$.  Choose $M =  KL $.  Thus, the first condition of \eqref{growth conditions FC and families} holds.    By the family growth hypothesis of \ref{compact-support}, the second condition of \eqref{growth conditions FC and families} holds since $M$ is now a fixed constant.
 
The conclusion of Theorem \ref{Hecke-multiples-smooth}, therefore, holds for \ref{polynomial decay}-\ref{compact-support}.

\begin{remark}\label{natural families S(N,k)}
    The framework of our results applies equally to other natural families.  For instance, consider the Hecke eigenbasis of $S_k(\Gamma_0(N))$, that is, a basis of $S_k(\Gamma_0(N))$ consisting of simultaneous eigenforms of the Hecke operators $T_n$, $(n,N)=1$, with $\theta_f(p)$ defined through the eigenvalues of $T_p$, $p\nmid N$.  Averages of the eigenvalues over such a basis are normalized traces of Hecke operators on $S_k(\Gamma_0(N))$, and are independent of the choice of basis.  The eigenvalues are multiplicative at arguments coprime to $N$, and the primes dividing $N$ can be discarded as in Setup \ref{spaces-rvs-modular}.  Equidistribution laws and discrepancies with respect to this family have been explored in \cite{Serre} and \cite{MS}.  As in  Proposition \ref{trace formula}, the Eichler--Selberg trace formula estimate for $T_n$ acting on  $S_k(\Gamma_0(N))$ is uniform in $k$ and $N$ (see \cite[Section 4]{Serre}, \cite[Theorem 29.5]{Knightly-Li} and \cite[Section 7-8]{MS}).  The proof of Theorem \ref{Hecke-multiples-smooth} carries over verbatim, the growth conditions determined by the error term of that estimate.
\end{remark}

\begin{cor}\label{Nagoshi-revisited}
    Fix an integer $j\geq 1$, and let $N=N(x)$ and $k=k(x)$ satisfy the conditions 
    $$\omega(N) = \o\left(\sqrt{\pi(x)}\right)\text{ and } \frac{\log (kN/8^{\omega(N)})}{\log x} \to \infty.$$
    Then the moments of $\displaystyle\frac{1}{\sqrt{\pi(x)}}\sum_{p\leq x}a_f(p^j)$ over $f\in\mathcal F_{N,k}$ converge to the Gaussian moments. 
\end{cor}
\begin{proof}
     We apply Theorem \ref{2-periodic} to $H(t) = U_j(\cos (\pi t))$, $\alpha = 1$, and the setup and exceptional sets of Setup \ref{spaces-rvs-modular}.  The thinness condition holds since $\omega(N)=\o(\sqrt{\pi(x)})$ (see \eqref{thin-sets-modular}).  For $p \nmid N$, $H(\thtp) = a_f(p^j)$.  For $p \mid N$, we have $|a_f(p^j)|\le1$, so these primes contribute at most $\omega(N)=\o(\sqrt{\pi(x)})$ to the sum, uniformly in $f$.   Since $H$ is a trigonometric polynomial of degree $j$ in $\pi t$, condition \ref{Cond GC 1} holds trivially with $M=j$.  
     
     For $p\nmid N$ and every $m\neq 0$, $\displaystyle \cos(\pi m\theta_f(p))=\frac{1}{2}\left(a_f(p^{|m|})-a_f(p^{|m|-2})\right)$, with $a_f(p^{-1}):=0$.
     For even $m$, the averages $C_{m,p,N,k}(1)$ are the averages $\tilde c_{m/2,p,N,k}$ of Proposition \ref{Weyl-limits-modular}.  For odd $m$, the arguments $p^{|m|}$ and $p^{|m|-2}$ are not squares, so Proposition \ref{trace formula} gives
     $$C_{m,p,N,k}(1)=\O\left(\frac{8^{\omega(N)}}{kN}\,p^{|m|d'}\right),$$ and hence $D_m(1)=0$.  In particular, 
     $$D_m(1) = \delta_{m,0} -\frac{1}{2}\delta_{|m|,2}.$$
     The proof of Proposition \ref{Towards GC 2} uses only the multiplicativity of $a_f(n)$ and Proposition \ref{trace formula}.  So, it applies to $m_i$'s of any parity, and the proof of Proposition \ref{Towards GC 3} applies with $M = j$ fixed.  Conditions \ref{Cond GC 2} and \ref{Cond GC 3} then follow as in Section \ref{sec: modular CLT proof}, under the growth condition \eqref{xMC}.

     Finally, $\nu(t) = \sum_{m \in \Z}D_m(1)\cos \pi mt = 2\sin^2\pi t$ , and by the orthogonality relation \eqref{eq:mean zero}, $\int_0^1 H(t) \nu(t) dt= 0$, and 
   $V_H(\nu) = \int_0^1 H(t)^2 \nu(t) dt = 1.$  Since $A(x) = \pi(x)$, Theorem \ref{2-periodic} gives the corollary.
     
\end{proof}

\begin{remark}

    For $N=1$, $j = 1$ and $\log k/\log x \to \infty$, Corollary \ref{Nagoshi-revisited} is a theorem of Nagoshi \cite{Nagoshi}, and for fixed $N$ and any $j \geq 1$, it recovers \cite[Theorem 6.1]{murty-prabhu}.  In the level aspect,  \cite[Theorem 4.1]{Cho-Kim} treats Nagoshi's statistic ($j=1$) for fixed $k$ as $\log N/\log x \to \infty$, averaging over a basis of $S_k(\Gamma_0(N))$ consisting of eigenforms of the Hecke operators $T_n, (n,N) = 1$; equivalently, over newforms of levels $M \mid N$, each counted with multiplicity $d(N/M)$.  Corollary \ref{Nagoshi-revisited} concerns newforms of level exactly $N$, and treats every $j\geq 1$, with $k$ and $N$ varying simultaneously.  We have restricted our treatment to these families for simplicity of exposition, but the framework applies equally to the  statistic considered in \cite{Cho-Kim} (see Remark \ref{natural families S(N,k)}).  Since the Hecke relations at a prime $p$ hold only for $p \nmid N$, the primes dividing $N$ are treated separately; in Corollary \ref{Nagoshi-revisited}, as in Theorem \ref{Hecke-multiples-smooth}, their contribution is negligible under the condition $\omega(N) = \o\left(\sqrt{\pi(x)}\right)$.

\end{remark}

\section{Smooth statistics: comparison of the three families}\label{sec: smooth comparison}

We conclude our treatment of smooth statistics by comparing how the hypotheses of Theorem \ref{CLT general smooth variant} are verified in the three families. This explains, in particular, why the growth conditions and the admissible test functions in Theorems \ref{Elliptic-curves-smooth}-\ref{Hecke-multiples-smooth} take different forms.

\begin{enumerate}
\item Theorems \ref{Elliptic-curves-smooth}, \ref{Kloosterman-multiples-smooth} and \ref{Hecke-multiples-smooth} follow from Theorem \ref{CLT general smooth variant}, and in each case, all moments of $\mathcal S_x$ converge to the Gaussian moments.  In each family, the tail $\sum_p(Z-Z_M)$ is $\o(\sqrt{\pi(x)})$.  This follows from the decay of $\widehat{h}$, since $|\cos (\pi m \theta)| \leq 1$.  Thus, \ref{Cond GC 1} holds for every $l \geq 1$.  The exceptional sets are handled by estimating the moments of $\#\{p \leq x:\,p \mid ab\Delta(a,b)\}$ for elliptic curves and $\#\{p \leq x:\,p \mid a\}$ for Kloosterman angles.  They are handled trivially for Hecke angles because the exceptional primes are at most $\omega(N)$ primes dividing $N$.

\item The growth conditions on each family are what make the error term in \ref{Cond GC 2} and \ref{Cond GC 3} negligible for the chosen truncation level $M$.  The decay of $\widehat{g}$ determines how large $M$ must be.  For Hecke angles, the trace formula (Proposition \ref{trace formula}) gives an error of size $8^{\omega(N)}n^{d'}/kN$ at $n = p^{2m\alpha}$, which grows exponentially with respect to $m$.  Summed for $m$ up to $Mr$, this forces $\log (kN) \gg M\log x$.  Slower decay for $\widehat{g}$ requires a larger $M$, and hence a faster growth condition for $kN$.  This is demonstrated in the gradation \ref{polynomial decay}-\ref{compact-support} of Theorem \ref{Hecke-multiples-smooth}.  For Frobenius and Kloosterman angles, the averages are periodic modulo $p$ irrespective of the exponents of $p$, and so, $m$ enters the error terms only polynomially.  The constraint here is the modulus $q = p_1p_2\dots p_u$ in the Chinese Remainder Theorem argument for \ref{Cond GC 2}, which requires $A$ and $B$ (for Theorem \ref{Elliptic-curves-smooth}) and $N$ (for Theorem \ref{Kloosterman-multiples-smooth}) to grow faster than any power of $x$, independently of the test function.  Under this condition, any polynomial truncation level of the form $M = x^{\eta}$ works, and a single growth condition covers all test functions with $\kappa > 2$.

\item In the proof of Theorem \ref{Hecke-multiples-smooth}, the trace formula estimates in Proposition \ref{trace formula} are uniform with respect to both $N$ and $k$.  This allows us to vary $N(x)$ and $k(x)$ simultaneously under the restriction $\omega(N) = \o\left(\sqrt{\pi(x)}\right)$.  The primes dividing $N$ at which the Hecke relations at prime powers do not hold are removed by Theorem \ref{CLT general smooth variant} with $\mathcal G_{p,N,k} = \mathcal F_{N,k}$ for $p \mid N$.

\item For $p\nmid N$, the vertical limiting distribution of $\theta_f(p)$, as $f$ varies over $\mathcal F_{N,k}$ with $Nk \to \infty$, is the $p$-adic Plancherel measure \begin{equation}\label{p-adic plancherel measure}
    \mu(t) = 2\frac{(p+1)\sin^2(\pi t)}{\left(\sqrt{p} + \frac{1}{\sqrt{p}}\right)^2 - 4\cos^2(\pi t)}.
\end{equation} 
Its Fourier coefficients tend to the Sato–Tate coefficients as $p \to \infty$. Hence, the Weyl limits along the primes are the same as in the fixed-level case, and varying $N$ does not change $\mu_\alpha$.

\end{enumerate}

\section{Sharp statistics: scope and limits}\label{sec:linear-statistics}

We conclude this article by discussing how far the framework of Section \ref{sec: Smooth statistics} extends to sharp statistics, that is, to the number of angles lying in a fixed interval.  We prove a general central limit theorem for such statistics (Theorem \ref{CLT gen linear}), apply it to Hecke angles (Theorem \ref{Prabhu-Sinha type CLT, diff}) and explain the limitations of the method when applied to Frobenius and Kloosterman angles.

For a positive real number $x$ and a fixed interval $I=[a_1,b_1]\subseteq (0,1)$, consider the random variable $\mathcal N_{I,x}$ defined on $\f_k$ as follows:
\begin{equation}\label{N_(I,x)}
   \mathcal N_{I,x}(X)=\sum_{n \in \mathcal{A} \cap [1,x]}\chi_I(X_n).
\end{equation}
In this section, we address the distribution of
\begin{equation}\label{S_(I,x)}
  \mathcal{S}_{I,x}(X) :=\frac{\mathcal N_{I,x}(X) - A(x) \nu(I)}{\sqrt{A(x)(\nu(I) - \nu(I)^2)}},  
\end{equation}

where $\nu(t) = \sum_{m \in \Z} D_m \cos(\pi mt)$ (see Theorem \ref{2-periodic}) and $\nu(I) := \int_I \nu(t)\, dt$.

\subsection{Approximation by trigonometric polynomials}
We address the sharp statistics with the help of trigonometric polynomials which approximate the indicator function of an interval. We recall this approximation in the following lemma.
\begin{lemma}\label{Selberg Beurling}
Let $I = [a_1,b_1] \subset (0,1)$, and let $M$ be a positive integer.  There exist real and even, 2-periodic trigonometric polynomials 
$$S_{M}^{\pm}(t) = \sum_{|n| \leq M}\widehat{S_M^{\pm}}(n) e^{i\pi nt}$$
with Fourier coefficients $$\widehat{S_M^{\pm}}(n) := \frac{1}{2}\int_{-1}^ 1 S_M^{\pm}(t) e^{-i\pi nt}dt =  \int_0^1 S_M^{\pm}(t)\cos(\pi nt) dt$$
such that
\begin{enumerate}[label=\normalfont(\alph*)] 
    \item \label{upper-lower-bounds} For all $t \in [0,1]$,
    $\displaystyle S_{M}^{-}(t) \leq \chi_I(t) \leq S_{M}^{+}(t).$
    \item \label{dist-SB-chi} $\displaystyle\|S_{M}^{\pm}(t) - \chi_I(t)\|_{L^1} \leq \frac{2}{M+1}.$
    \item \label{FC-SB-even}
    The Fourier coefficients $\widehat{S_{M}^{\pm}}(n)$ are even.
    \item \label{FC-0-SB}
    $\displaystyle \abs{\widehat{S_{M}^{\pm}}(0) - (b_1-a_1)} \leq \frac{2}{M+1}.$
    \item \label{FC-SB-bounds}For all integers $n$ such that $1 \leq |n| \leq M$, we have
    $$\widehat{S_{M}^{\pm}}(n) \ll \frac{1}{M+1} + \min\left\{b_1-a_1,\frac{1}{|n|}\right\}.$$

    \item \label{SM uniformly bounded}
    Denote $F_M(t) := S_M^+(t) - S_M^-(t)$.  Then, $0 \leq \displaystyle F_M(t) \ll 1$ and $\displaystyle \abs{S_M^{\pm}(t)} \ll 1$ uniformly in $M$ and $t$.
    
\end{enumerate}
\end{lemma}

\begin{proof}
    We refer the reader to \cite[Chapter 1]{Mont} for a  construction of the Beurling--Selberg trigonometric polynomials, which approximate the indicator function of $I$.  However, the Fourier coefficients of these polynomials are not necessarily even.  These polynomials can be suitably adapted into trigonometric polynomials with even Fourier coefficients in order to prove \ref{upper-lower-bounds}, \ref{dist-SB-chi} and \ref{FC-SB-even} in the above lemma (see, for example, \cite[Section 3]{Prabhu-Sinha} and \cite[Lemma 3.1]{RT}).  

    Part \ref{FC-0-SB} follows by \ref{upper-lower-bounds} and \ref{dist-SB-chi}, as
    $$\widehat{S_M^\pm}(0)=\int_0^1 S_M^\pm(t)\, dt=|I|+\int_0^1(S_M^\pm-\mathbf 1_I)(t)\, dt.$$
    Similarly,
    $$\biggl|\int_0^1 S_M^{\pm}(t) \cos (\pi nt)\,dt\biggr| \leq \|S_M^{\pm} - \chi_I\|_{L^1} + \biggl|\int_I \cos (\pi nt )dt\biggr|,$$
    which leads to \ref{FC-SB-bounds}. 
    
    Finally, $F_M(t) = S_M^+(t) - S_M^-(t) \geq 0$ and hence by \ref{dist-SB-chi}, we have $$\int_0^1 F_M(t) dt \leq \frac{4}{M+1}.$$  Thus, for each $n$, $$\abs{\widehat{F_M}(n)} \leq 4/(M+1),$$ and hence, $$0 \leq F_M(t) \leq \sum_{|n| \leq M}\widehat{F_M}(n) \leq 8.$$ By \ref{upper-lower-bounds}, we have
    $$0 \leq \chi_I \leq S_M^+ = S_M^- + F_M \leq \chi_I + F_M \leq 9$$
    and 
    $$-8 \leq -F_M \leq S_M^+ - F_M = S_M^- \leq \chi_I \leq 1,$$
    which proves the uniform boundedness of $S_M^\pm(t)$ in $M$ and $t$.
\end{proof}
\begin{remark} The function $S_{M}^{\pm}(t)$ and its Fourier coefficients depend on $M$ as well as the interval $I$.  As the interval $I$ remains fixed, we have chosen to not indicate it in the notation.
\end{remark}

\begin{lemma}\label{discrepancy-SB-polynomial}
  
    Let $\nu$ be a non-negative, bounded, measurable function on $[0,1]$. For an interval $I = [a_1,b_1] \subset (0,1)$,  define
   $$\nu(I) := \int_0^1 \chi_I(t) \nu(t) dt \text{ and }V_I := \nu(I) - \nu(I)^2,$$ 
   and for a bounded function $G$ on $[0,1]$, write  
$$V_{G}(\nu) = \int_0^1 G(t)^2 \nu(t) dt - \left(\int_0^1 G(t)\nu(t) dt\right)^2.$$
Let $S_M^{\pm}$ be as in Lemma \ref{Selberg Beurling}, and put $F_M := S_M^+ - S_M^-$.  
Then, for any positive integer $M$,
\begin{enumerate}[label=\normalfont(\alph*)]
\item \label{mu(S_M)}
$$\int_0^1 S_M^{\pm}(t)\nu(t) dt = \nu(I) + \O\left(\frac{1}{M}\right),$$
\item \label{Variance: Z_M gen lin} $\displaystyle V_{S_M^{\pm}}(\nu) = V_I + \O\left(\frac{1}{M}\right).$
\item \label{FM} $0 \leq F_M(t) \ll 1$ for all $t \in [0,1]$ and 
$$\int_0^1 F_M(t)\nu(t) dt \ll \frac{1}{M} \text{ and }V_{F_M}(\nu) \ll \frac{1}{M}.$$
\end{enumerate}
The implied constants in all the above estimates depend only on $\|\nu\|_{\infty}$.
\end{lemma}

\begin{proof}
To prove \ref{mu(S_M)}, we observe that
$$\int_0^1 S_M^{\pm}(t)\nu(t) dt = \int_0^1 (S_M^{\pm}(t) - \chi_I(t))\nu(t) dt + \int_0^1 \chi_I(t)\nu(t) dt.  $$
\ref{mu(S_M)} follows by Lemma \ref{Selberg Beurling}\ref{dist-SB-chi} as
    $$\abs{\int_0^1 (S_M^{\pm}(t) - \chi_I(t))\nu(t) dt } \leq \frac{2 \|\nu\|_{\infty}}{M+1}.$$

 To prove \ref{Variance: Z_M gen lin}, note that
$$\int_0^1 S_M^{\pm}(t)^2\nu(t) dt = \int_0^1 (S_M^{\pm}(t)^2 - \chi_I(t))\nu(t) dt + \int_0^1 \chi_I(t)\nu(t) dt,  $$
 and
 \begin{align*}
&\left|\int_0^1 (S_M^{\pm}(t)^2 - \chi_I(t))\nu(t) dt\right| = \left|\int_0^1 (S_M^{\pm}(t)^2 - \chi_I(t)^2)\nu(t) dt\right|\\ &\leq \|\nu\|_{\infty}\int_0^1 \left|\Bpm(t)-\chi_I(t)\right|\left|\Bpm(t)+\chi_I(t)\right|\,dt =\O\left(\frac{1}{M+1}\right), 
 \end{align*}
where the last inequality follows by Lemma \ref{Selberg Beurling}\ref{dist-SB-chi} and the fact that $S_{M}^{\pm}(t)$ is absolutely bounded (see Lemma \ref{Selberg Beurling}\ref{SM uniformly bounded}). Combining the above inequality with \ref{mu(S_M)}, we get \ref{Variance: Z_M gen lin}.  

 For \ref{FM}, the nonnegativity follows by Lemma \ref{Selberg Beurling}\ref{upper-lower-bounds}.  The boundedness of $F_M$ follows from Lemma \ref{Selberg Beurling}\ref{SM uniformly bounded}.  We also have 
 \begin{align*}
\int_0^1 F_M(t) \nu(t) dt \leq \|\nu\|_{\infty}\left(\|S_M^+ - \chi_I\|_{L^1}  + \|\chi_I - S_M^-\|_{L^1}\right)   \ll \frac{1}{M},
 \end{align*}
 and
 \begin{align*}
     V_{F_M}(\nu) \leq \int_0^1 F_M^2(t)\nu(t) dt \leq \|F_M\|_{\infty}\int_0^1 F_M(t) \nu(t) dt \ll \frac{1}{M}.
 \end{align*}
 
\end{proof}

\subsection{A general CLT for sharp statistics}
In this subsection, we state a general CLT  for the sharp statistics $\mathcal N_{I,x}$.

\begin{prop}\label{rem:uniform}
Assume the setup of Theorem \ref{2-periodic} with $\nu(t) = \sum_{m \in \Z} D_m \cos (\pi mt)$.  Let $G=G_x$ be real, even, $2$-periodic trigonometric polynomials of degree at
most $M(x)$ with $\sup_x\|G_x\|_\infty<\infty$, and put
$$\mathfrak C_G(X):=\sum_{n\in\mathcal A\cap[1,x]}G(X_n)-A(x)\int_0^1 G(t)\nu(t)\,dt.$$
Suppose \eqref{thin-sets}, \ref{Cond GC 2} and \ref{Cond GC 3} hold for $G_x$ with this $M$.  Then for every $r\geq 1$,
\[
\langle \mathfrak C_G^r\rangle_k=\gamma_r\big(A(x)V_G(\nu)\big)^{r/2}+\o\big(A(x)^{r/2}\big).
\]
\end{prop}

\begin{proof}
Since the degree of $G_x$ is at most $M(x)$, we take $Z_M = Z$.  Thus, \ref{Cond GC 1} holds trivially.  In the proofs of Theorems \ref{CLT general smooth}, \ref{CLT general smooth variant} and \ref{2-periodic}, the test function enters only through the bounds of Lemma \ref{properties of powers of Z_M}\ref{ZM-power-F}, \ref{DMr-all} and $D(M,2)$. For $h = G_x$, we have
\begin{align*}
&\left|\widehat{Z_M^r}(m)\right|\leq (2\sup_x\|G_x\|_\infty)^r, \text{ for each }m,\\
&D(M,r) \leq (2\sup_x\|G_x\|_\infty)^r, \text{ and }\\
& D(M,2) = V_G(\nu).
\end{align*}
Applying Steps \hyperlink{step1}{1} - \hyperlink{step4}{4} in the proof of Theorem \ref{CLT general smooth} (in the setting of Theorem \ref{CLT general smooth variant}), we get the proposition.

\end{proof}
\begin{unmremark}

    $G_x$ supplies every property that the proof of Theorem \ref{CLT general smooth} (and of Theorems \ref{CLT general smooth variant} and \ref{2-periodic}) actually uses, even though it does not satisfy the literal hypothesis of these theorems of being a fixed function.  However, in contrast to Theorems \ref{CLT general smooth} etc, $V_{G_x}(\nu)$ varies with $x$ and may even go to zero (for example, $V_{F_M}(\nu)$ from the previous lemma).  That is why, in the above proposition, we take the unnormalized $r$-th moments of $\mathfrak C_G$.

\end{unmremark}

\begin{thm}\label{CLT gen linear}
    Assume the setup of Theorem \ref{2-periodic} with $\nu(t) = \sum_{m \in \Z} D_m \cos (\pi mt)$.  Let $I \subset (0,1)$ be a fixed interval, let $S_M^{\pm}$ be as in Lemma \ref{Selberg Beurling}, and let $M = M(x)$ satisfy $\displaystyle {M(x)}/{\sqrt{A(x)}} \to \infty.$  Suppose that the hypotheses of Theorem \ref{2-periodic} hold for the test functions $S_M^{-}$ and $F_M = S_M^+ - S_M^-$ with this $M$.  Let $\nu(I)$ and $V_I$ be as defined in Lemma \ref{discrepancy-SB-polynomial}, and assume that $0 < \nu(I) < 1$, so that $V_I >0$. Then, all moments of 
    \begin{equation*}
  \mathcal{S}_{I,x}(X) :=\frac{\mathcal N_{I,x}(X) - A(x) \nu(I)}{\sqrt{A(x)V_I}},  
    \end{equation*}    
converge to the moments of the Gaussian distribution.  That is, for any integer $r\geq 1$, we have
      \[
  \lim_{x \to \infty} \bigl\langle \mathcal S_{I,x}(X)^r \bigr\rangle_{k(x)} = \gamma_r.
      \]

\end{thm}

\begin{proof}
  Denote $\sigma_x = \sqrt{A(x)V_I}$, where $V_I = \nu(I) - \nu(I)^2$.  For any test function $G$, denote
$$\mathfrak C_G(X) = \sum_{n \in \mathcal A \cap [1,x]}G(X_n) - A(x)\int_0^1 G(t)\nu(t)dt.$$
We write
\begin{align}\label{decompose-N}
\frac{\mathcal N_{I,x}(X) - A(x)\nu(I)}{\sigma_x} = \frac{1}{\sigma_x}\mathfrak C_{S_M^-}(X) + E_x(X),
\end{align}
where
\begin{align}\label{ExX}
E_x(X) = \frac{1}{\sigma_x}\sum_{n \in \mathcal A \cap [1,x]}(\chi_I - S_M^-)(X_n) + \frac{A(x)}{\sigma_x}\left(\int_0^1 S_M^-(t) \nu(t) dt - \nu(I)\right).
\end{align}
By Lemma \ref{Selberg Beurling}\ref{SM uniformly bounded}, the function $S_M^-$ satisfies the hypotheses of Proposition \ref{rem:uniform}.  Thus, 
\begin{align*}
\left \langle \left(\mathfrak C_{S_M^-}(X)\right)^r \right \rangle = \gamma_r (A(x) V_{S_M^-}(\nu))^{r/2} + \o\left( A(x)^{r/2}\right).
\end{align*}
By Lemma \ref{discrepancy-SB-polynomial}\ref{Variance: Z_M gen lin}, $V_{S_M^-}(\nu) = V_I + \O(1/M)$.
Thus, for every $r \geq 1$,
\begin{align}\label{lower-Selberg-moments}
\left \langle \left(\frac{1}{\sigma_x}\mathfrak C_{S_M^-}(X)\right)^r \right \rangle \to  \gamma_r \text{ as }x \to \infty.
\end{align}
Now, by Lemma \ref{Selberg Beurling}\ref{upper-lower-bounds}, $0 \leq \chi_I - S_M^- \leq F_M$.  So, by Minkowski's inequality, for every $j \geq 1$,
\begin{align}\label{Exq}
\|E_x\|_{2j} \leq \frac{1}{\sigma_x}\left\|\mathfrak C_{F_M}\right\|_{2j} + \frac{A(x)}{\sigma_x} \left(\int_0^1 F_M(t) \nu(t) dt + \abs{\int_0^1 S_M^-(t)\nu(t)dt -\nu(I)}\right).   \end{align}
By Lemma \ref{discrepancy-SB-polynomial}\ref{mu(S_M)} and \ref{FM}, the second term is $\ll \sqrt{A(x)}/M = \o(1)$.  By Proposition \ref{rem:uniform} applied to $G = F_M$, we have
\begin{align*}
   \left\|\mathfrak C_{F_M}\right\|_{2j}^{2j} =  \left \langle \left(\mathfrak C_{F_M}\right)^{2j}\right \rangle = \gamma_{2j}(A(x)V_{F_M}(\nu))^{j} + \o\left(A(x)^{j}\right)
    \end{align*}
and further, since $V_{F_M}(\nu) \ll 1/M$ by Lemma \ref{discrepancy-SB-polynomial}\ref{FM},
\begin{align*}
    \left \langle \left(\mathfrak C_{F_M}\right)^{2j}\right \rangle = \o\left(A(x)^{j}\right).\end{align*}
Combining the above inequality and \eqref{Exq}, we get $\|E_x\|_{2j} \to 0$ for each $j \geq 1$.  We repeat Step \hyperlink{step5}{5} of the proof of Theorem \ref{CLT general smooth} with
$$\frac{1}{\sigma_x}\mathfrak C_{S_M^-} \text{ and } E_x.$$
By \eqref{decompose-N} and \eqref{lower-Selberg-moments}, we get that
$$\left \langle\left(\frac{\mathcal N_{I,x}(X) - A(x)\nu(I)}{\sigma_x}\right)^r\right \rangle \to \gamma_r$$
as $x \to \infty$, proving the theorem.

\end{proof}

\begin{cor}\label{CLT gen linear union}
Assume the setup of Theorem \ref{2-periodic} with $\nu(t) = \sum_{m \in \Z} D_m \cos (\pi mt)$.  Let $J \geq 1$ and let $\mathcal I=\bigcup_{j=1}^JI_j$, where $I_1,\dots,I_J\subset(0,1)$ are fixed, pairwise disjoint closed intervals. Let $S_{M,j}^\pm$ be the polynomials of Lemma \ref{Selberg Beurling} for the interval $I_j$, and put
$$\mathbb S_M^\pm:=\sum_{j=1}^JS_{M,j}^\pm,\qquad\mathbb F_M:=\mathbb S_M^+-\mathbb S_M^-.$$
Let $M=M(x)$ satisfy $M(x)/\sqrt{A(x)}\to\infty$, and suppose that the hypotheses of Theorem \ref{2-periodic} hold for the test functions $\mathbb S_M^-$ and $\mathbb F_M$ with this $M$.  Put $\nu(\mathcal I):=\int_{\mathcal I}\nu$ and $V_{\mathcal I}:=\nu(\mathcal I)-\nu(\mathcal I)^2$, and assume that $0<\nu(\mathcal I)<1$. Then, for every integer $r\geq 1$,
$$\lim_{x\to\infty}\Bigg\langle\Bigg(\frac{\sum_{n\in\mathcal A\cap[1,x]}\mathbf 1_{\mathcal I}(X_n)-A(x)\nu(\mathcal I)}{\sqrt{A(x)V_{\mathcal I}}}\Bigg)^r\Bigg\rangle_{k(x)}=\gamma_r .$$

\end{cor}

\begin{proof}
    Since the intervals $I_j$ are disjoint, $\chi_{\mathcal I}=\sum_j\chi_{I_j}$. Summing the properties in Lemma \ref{Selberg Beurling} over $j$ shows that $\mathbb S_M^\pm$ are real, even, 2-periodic trigonometric polynomials of degree at most $M$ such that

    \begin{itemize}

\item $\displaystyle \mathbb S_M^-\le\chi_{\mathcal I}\le\mathbb S_M^+$ on $[0,1]$,
\item $\displaystyle\|\mathbb S_M^\pm-\chi_{\mathcal I}\|_{L^1}\leq\frac{2J}{M+1}$, and
\item $\displaystyle|\mathbb S_M^\pm(t)|\leq 9J$ and $0\leq\mathbb F_M(t)\leq 8J$ for all $t$.

\end{itemize}
The proof of Lemma \ref{discrepancy-SB-polynomial} uses only these three properties, together with $\chi_{\mathcal I}^2=\chi_{\mathcal I}$. Hence Lemma \ref{discrepancy-SB-polynomial} holds with $I$, $S_M^\pm$, $F_M$ replaced by $\mathcal I$, $\mathbb S_M^\pm$, $\mathbb F_M$, with implied constants depending only on $J$ and $\|\nu\|_\infty$. Since $\sup_x\|\mathbb S_M^-\|_\infty$ and $\sup_x\|\mathbb F_M\|_\infty$ are finite, Proposition \ref{rem:uniform} applies to $G=\mathbb S_M^-$ and $G=\mathbb F_M$. The proof of Theorem \ref{CLT gen linear} now applies mutatis mutandis.
\end{proof}

\subsection{Application to Hecke angles}
As an immediate consequence of Corollary \ref{CLT gen linear union}, we derive a CLT for the sharp statistics for Hecke angles.  This generalizes \cite[Theorem 1.1]{Prabhu-Sinha} to the reduced dilated angles $\{\alpha \thtp\}$ for each fixed integer $\alpha \geq 1$, and to varying level $N$.

\begin{thm}\label{Prabhu-Sinha type CLT, diff}
    Let $I=[a_1,b_1]$ be a fixed interval in $(0,1)$. Fix a positive integer $\alpha$.  For  a positive real number $x$ and $f\in \f_{N,k}$, let
    \begin{equation*}\label{N_(I,x)(f,a)}
        \mathcal{N}_{I, x}(f, \alpha)=\#\{p\leq x\;:\; \{\alpha\thtp\}\in I\},
    \end{equation*}  
    where $\{t\}$ denotes the fractional part of a real number $t$.  Denote $$\mu_{\alpha}(I) :=  \int_I \mu_{\alpha}(t) dt  = \begin{cases} \displaystyle\int_I 2\sin^2 (\pi t) \,dt &\text{ if }\alpha = 1\\
    \\
    \displaystyle \int_I dt =b_1 - a_1 &\text{ if } \alpha \geq 2.
    \end{cases}$$
   
    Suppose that $N = N(x)$ and $k=k(x)$ with $k$ even satisfy the growth conditions 
    \begin{equation}\label{gc-sharp}
    \frac{\log (kN/8^{\omega(N)})}{\sqrt{x}\log x}\to\infty\text{ and }\omega(N) = \o\left(\sqrt{\pi(x)}\right)\text{ as }x\to\infty.
    \end{equation}
    Then for any integer $r\geq 1$, we have 
    \begin{equation*}
        \lim_{x\to\infty}\frac{1}{\abs{\mathcal F_{N,k}}}\sum_{f\in\mathcal F_{N,k}}\left(\frac{\mathcal{N}_{I, x}(f, \alpha)-\pi(x){\mu}_{\alpha}(I)}{\sqrt{\pi(x)\left({\mu}_{\alpha}(I)-{\mu}_{\alpha}(I)^2\right)}}\right)^r=\frac{1}{\sqrt{2\pi}}\int_{-\infty}^\infty t^r e^{-t^2/2}\,dt.
    \end{equation*}
\end{thm}
\begin{proof}

    For a positive integer $\alpha$, denote
$$\mathcal I := \bigcup_{j=0}^{\alpha -1}\frac{I + j}{\alpha} \text{ and }\tilde{\mu}_{\alpha}(\mathcal I) :=  \sum_{j=0}^{\alpha - 1} \mu_{1}\left(\frac{I + j}{\alpha}\right).$$

We observe that 
$$\{\alpha \theta_f(p)\} \in I \iff \theta_f(p) \in \mathcal I.$$  
Thus,
 \begin{align}\label{Step-1-sharp}
 \mathcal{N}_{I, x}(f, \alpha) = \sum_{p \leq x} \chi_{\mathcal I}(\theta_f(p)).
 \end{align}
 We apply Corollary \ref{CLT gen linear union} to the Hecke angles $(\thtp)_p$ with $\mathcal A$ and $\mathcal G_{p,N,k}$ as defined in Setup \ref{spaces-rvs-modular}, and to $\mathcal I$, a disjoint union of intervals.   
Denote $S_{M,j}^{\pm}$ to be the Beurling--Selberg polynomial of degree $M$ corresponding to the interval $(I + j)/\alpha$, and let
$\mathbb S_M^{\pm}(t) := \sum_{j=0}^{\alpha - 1}S_{M,j}^{\pm}(t)$.  

We take $M(x) = \lceil \sqrt{\pi(x) \log \log x} \rceil$, and verify the hypotheses of Theorem \ref{2-periodic}.  The verification of \eqref{thin-sets} follows by \eqref{thin-sets-modular} since $\omega(N) = \o\left(\sqrt{\pi(x)}\right)$.  \ref{Cond GC 1} holds trivially since $\mathbb S_M^{\pm} $ is a trigonometric polynomial.  Similar to the proofs of Theorem \ref{Hecke-multiples-smooth} and Corollary \ref{Nagoshi-revisited}, in the verification of \ref{Cond GC 2} and \ref{Cond GC 3}, the trace formula estimates give us error terms of the form 
$\displaystyle \frac{8^{\omega(N)}}{kN}\,(Mr)^u\,\pi(x)^u\,x^{Mrd'}.$
Since $Mr\log x\ll_r\sqrt{x\log x\,\log\log x}=\o(\sqrt x\log x)$, the first growth condition in \eqref{gc-sharp} gives $x^{CMr}=\o\left(kN/8^{\omega(N)}\right)$ for every fixed $C>0$. Hence the right hand side of the expressions in both \ref{Cond GC 2} and \ref{Cond GC 3} can be shown to be $\o(\pi(x)^{r/2})$.  In \ref{Cond GC 3}, the main terms contribute $\ll(\omega(N)+\log\log x)^u=\o(\pi(x)^{u/2})$, since for $p\nmid N$ the differences $|C_{m,p}-D_m|$ decay geometrically in $m$, as in the proof of Proposition \ref{Towards GC 3}.

The verification of \ref{Cond GC 1}-\ref{Cond GC 3} for $\mathbb F_M = \mathbb S_M^+ - \mathbb S_M^-$ is identical. 

By Corollary~\ref{CLT gen linear union}, we get that
\begin{equation}\label{linear-stats-mod}
\lim_{x \to \infty}\left \langle\left(\frac{\sum_{p \leq x}\chi_{\mathcal I}(\theta_f(p)) -  \pi(x)\tilde{\mu}_{\alpha}(\mathcal I)}{\sqrt{\pi(x) (\tilde{\mu}_{\alpha}(\mathcal I) - \tilde{\mu}_{\alpha}(\mathcal I)^2)}}\right)^r \right \rangle =\begin{cases}
                        \frac{r!}{2^{r/2} (r/2)!} &\text{ if }r \text{ is even},\\
                        0 &\text{ if }r \text{ is odd.}
                \end{cases}
\end{equation}

Finally, we show that $\tilde{\mu}_{\alpha}(\mathcal I) = \mu_{\alpha}(I)$ for each $\alpha \geq 1$.  This assertion is immediate for $\alpha = 1$. For $\alpha \geq 2$,
 \begin{align*}
 \tilde{\mu}_{\alpha}(\mathcal I) &= \sum_{j=0}^{\alpha - 1} \int_{(I+j)/\alpha}2 \sin^2(\pi t )\,dt  =  \sum_{j=0}^{\alpha - 1} \frac{1}{\alpha}\int_I 2 \sin^2\left(\frac{\pi (s+j)}{\alpha}\right) \,ds\\
&=\sum_{j=0}^{\alpha-1}\frac{1}{\alpha} \int_I \left(1 - \cos \left(\frac{2\pi(s+j)}{\alpha}\right)\right)ds = \int_I ds = \mu_{\alpha}(I),
 \end{align*}
 since
 \begin{align*}
 \sum_{j=0}^{\alpha - 1} \cos\left(\frac{2\pi(s+j)}{\alpha}\right) = 0.
 \end{align*}
 Combining \eqref{Step-1-sharp} and \eqref{linear-stats-mod} proves the theorem.
\end{proof}

\begin{remark}\label{limitations of Theorem 52}\leavevmode
  \begin{enumerate}[label=(\alph*)]
    \item Theorem \ref{Prabhu-Sinha type CLT, diff} is the only case in this article in which our method yields a central limit theorem for sharp statistics; we compare it with earlier work and explain the growth condition of this theorem.
            \begin{enumerate}[label=(\arabic*)]
            \item For $\alpha = 1$ and fixed $N$, Theorem \ref{Prabhu-Sinha type CLT, diff} recovers \cite[Theorem 1.1]{Prabhu-Sinha}, under the same growth condition $\log k/(\sqrt{x}\log x) \to \infty$.  Theorem \ref{Prabhu-Sinha type CLT, diff} extends it in two directions: to varying level provided $\omega(N) = \o\left(\sqrt{\pi(x)}\right)$ and to the reduced dilated angles $\{\alpha \thtp\}$ for all $\alpha \geq 2$.  These are uniformly distributed as the Weyl limits  $d_m(\alpha) = d_{m\alpha}(1) = 0$ for all $|m| \geq 1$.  Thus, the mean and variance are asymptotic to $\pi(x)|I|$ and $\pi(x)(|I| - |I|^2)$ respectively.
            \item The growth condition on $kN$ needed for Theorem \ref{Prabhu-Sinha type CLT, diff}, $\log(kN/8^{\omega(N)})/(\sqrt{x}\log x) \to \infty$, is much stronger  than the one needed in  Theorem \ref{Hecke-multiples-smooth}.  This is because the trace formula errors grow like $x^{Md'}$, and in order to approximate $\chi_I$ to within an error bound of order $\o\left(\sqrt{\pi(x)}\right)$, we have to choose a truncation degree $M \gg \sqrt{\pi(x)}$.
            \end{enumerate}
    \item \label{obstruction-linear} For Frobenius angles of elliptic curves and for Kloosterman angles, the method of this section does not lead to a central limit theorem for sharp statistics. 
    
            \begin{enumerate}[label=(\arabic*)]
            
                \item   The obstruction is in Condition \ref{Cond GC 3}.   The Fourier coefficients $\widehat{S_M^{\pm}}(m)$ are of size $\asymp 1/|m|$ for $1 \leq |m| \leq M$, and the available ``one-prime'' bounds (Proposition \ref{S(n)-elliptic}\ref{Spm-bounds} for elliptic curves and the Katz bound \eqref{Katz inequality} for Kloosterman sums) give us  $\abs{\tilde{c}_{m,p} -d_m} \ll |m|/\sqrt{p}$.  For $r=1$, the weighted sum in \ref{Cond GC 3} is bounded by $\sum_{|m| \leq M}|m|^{-1} |m|/\sqrt{p} \ll M/\sqrt{p}$.  Summed over the primes up to $x$, this is of order $M\sqrt{\pi(x)/\log x}$, which suffices   only for $M = \o(\sqrt{\log x})$, but the approximation for sharp statistics requires $M \gg \sqrt{\pi(x)}$.  On the other hand, for the Hecke angles, the differences $\abs{\tilde{c}_{m,p} - d_m}$ decay geometrically in $m$, which makes Theorem \ref{Prabhu-Sinha type CLT, diff} possible.
                
                \item    For a fixed prime $p$, the angle $\theta_E(p)$ can take at most $4 \sqrt{p} + 1$ distinct values as $(a,b)$ varies, since $\tilde{a}_E(p)$ is an integer with $\abs{\tilde{a}_E(p)} \leq 2\sqrt{p}$.  Similarly, the  Kloosterman angle $\theta_{p,a}$ depends only on $a \,  \tmod p$.  Let $\nu_p$ denote the distribution of these angles over the family at the prime $p$.  For the test functions $S_M^{\pm}$ of this section, the quantity controlled by \ref{Cond GC 3} at a single prime $p$ (when $u = r = 1$) is
            \begin{equation*}
                   \sum_{|m| \leq M} \widehat{S_M^{\pm}}(m)(\tilde{c}_{m,p} - d_m) = \int_0^1 S_M^{\pm} d\nu_p - \int_0^ 1 S_M^{\pm} d\nu, 
                \end{equation*}
                which is a ``smooth'' analogue of the deviation $\nu_p(I) - \nu(I)$ of the local distribution at $p$ from the limiting measure.  Condition \ref{Cond GC 3} therefore requires 
                \[\sum_{p \leq x}|\nu_p(I) - \nu(I)| = \o\left(\sqrt{\pi(x)}\right).\]
                Since $\sum_{p\le x} p^{-1/2}\asymp \sqrt{\pi(x)/\log x}$, the above requirement amounts to local discrepancy bounds of square-root strength, $|\nu_p(I) - \nu(I)| \ll p^{-1/2 }(\log p)^{1/2- \delta}$ for some $\delta >0$.  The best known bounds are $\O(p^{-1/4})$ in both families: for Kloosterman angles this follows from Katz’s estimates \cite{Nied}, and for the family of all elliptic curves modulo $p$, it was obtained by \cite{Banks-Shpar} and later, also by \cite{Murty-Prabhu-BAuMS} from the same inputs as in Section \ref{sec: Elliptic curves} (Deuring’s theorem, the Eichler–Selberg trace formula and Deligne’s bound). See also \cite{Miller-Murty}, where a logarithmic saving is obtained only using Birch’s moments.

            \end{enumerate}
    
\end{enumerate}

\end{remark}

\section{Notation}\label{sec: Notation}
We collect the general notation used throughout the paper.

\begin{xltabular}{\textwidth}{@{}p{2.5cm}X@{}}
$\mathcal{A} \subseteq \N$ & An infinite subset of $\N$\\
$A(x)$ & $\left|\mathcal{A} \cap [1,x]   \right|
$\\
$\mathcal{F}_k$, $k \in \N$ & A finite set equipped with the uniform probability measure. $| \mathcal{F}_k| \to \infty$ as $k \to \infty$. \\
$X =(X_n)_{n \in \mathcal A}$ & A sequence with $X_n \in [0,1]$ for all $n \in \mathcal A$. \\
$\mathcal{P}$ & Set of all prime numbers. \\
$\mathcal{P}_{\geq n}$ & Set of all prime numbers greater than or equal to $n$.\\
$\langle \Phi \rangle_k$ & Expectation/average value of a random variable $\Phi \colon \f_k \to \C$ under the uniform measure; see Setup \ref{spaces-rvs}. \\ 
$h \colon \R \to \R$ & Continuous, even, one-periodic function with $\sum_{m \in \Z} |\widehat{h}(m)| < \infty$; see Setup $5$, equation \eqref{h-FE}. \\
$h_n \colon\mathcal F_k \to \R$ & $h_n(X) := h(X_n)$. \\
$c_{m,n,k}$ & $m$-th family average of $\f_k$ at index $n$; see \eqref{Weyl sum of X}. \\
$d_m$ & $m$-th Weyl limit of the collection $(\mathcal F_k)_{k \geq 1}$ along $\mathcal A$; see \eqref{Weyl limit d_l}. \\
$\mu(t)\,dt$ & Symmetric probability measure on $[0,1]$ with $m$-th Fourier coefficient $d_m$; see Lemma~\ref{mu-pdf}. \\
$\mu(I)$       & $\displaystyle\int_I \mu(t)\,dt$. \\
$Z(t)$ & The centering of $h$ with respect to $\mu$; see \eqref{Z_M poly gen}. \\
$V_h(\mu)$ & Normalized asymptotic variance of $\sum_{n \in \mathcal A \cap [1,x]} h_n$; see \eqref{Vh-defn}. \\
$\sigma_x$  & $\sqrt{A(x)V_h(\mu)}$; asymptotic standard deviation of $\sum_{n \in \mathcal A \cap [1,x]} h_n$; used in the proof of Theorem~\ref{CLT general smooth variant}.\\
$h_M, Z_M$ & Truncated versions of $h$ and $Z$ respectively; see Lemma \ref{properties of powers of Z_M}.\\
$\gamma_r $ & $r$-th moment of the standard Gaussian distribution.\\
$\mathcal S_x(X)$ & Normalized fluctuation of $\sum_{n \in \mathcal A \cap [1,x]} h_n$; see Theorem \ref{CLT general smooth}.\\
$\displaystyle D(M,r)$ & $r$-th moment of $Z_M$ with respect to $\mu$; see Lemma \ref{properties of powers of Z_M}\ref{DMr-all}.\\
$ E_k(M,n,r)$ & $r$-th moment error term of $Z_M$ at index $n$; see \eqref{r-th moment error term of Z_M at index n}.\\ 
$\mathcal G_{n,k}$ & Thin subsets of $\f_k$ for each $n \in \mathcal{A}$; see \eqref{thin-sets}. \\
$\mathcal H_{n,k}$ & $\f_k \setminus \mathcal G_{n,k}$; $n$-th refined family.\\
$W(X_n)$ & Equals $Z(X_n)$ when $X \in \mathcal{H}_{n,k}$ and $0$ otherwise.\\
$\widetilde{c}_{m,n,k}$ & $m$-th refined family average of $\f_k$ at index $n$; see Theorem \ref{CLT general smooth variant}.\\
$\|Y\|_{\ell}$  & $L^{\ell}$-norm of $Y$ on $\f_k$; see the proof of Theorem~\ref{CLT general smooth variant}.\\ 
$H\colon \R \to \R$ & Continuous, even, two-periodic function with $\sum_{m \in \Z} |\widehat{H}(m)| < \infty$.\\
$C_{m,n,k}, D_m$ & Analogues of $\widetilde{c}_{m,n,k}$ and $d_m$ with $\cos(2\pi mX_n)$ replaced with $\cos(\pi mX_n)$; see Theorem~\ref{2-periodic}.\\ 
$\nu(t)\, dt$ & Probability measure on $[0,1]$ with $m$-th cosine coefficient $D_m$; see Theorem~\ref{2-periodic}.\\
$\nu(I)$ and $V_I$ & $\displaystyle\int_I \nu(t)\, dt$, resp.  $\nu(I) - \nu(I)^2$.\\
$S_M^{\pm}$ & Beurling--Selberg polynomials of degree $M$, approximating $\mathbf{1}_I$; see Lemma~\ref{Selberg Beurling}.\\
$\mathcal N_{I,x}(X)$ & Number of $n \in \mathcal A \cap [1,x]$ with $X_n \in I$, for an interval $I \subseteq (0,1)$; see \eqref{N_(I,x)}. \\
$\mathcal{N}_{I, x}(f, \alpha)$ & Number of primes $p \leq x$ with $\{\alpha\theta_f(p)\} \in I$; see Theorem~\ref{N_(I,x)(f,a)}. \\
$\mathcal{S}_{I,x}$ & Normalized fluctuation of $\mathcal N_{I,x}(X)$; see \eqref{S_(I,x)}.\\
$\mathfrak C_G(X)$ & The centering of $\sum_{n \in \mathcal A \cap [1,x]}G(X_n)$ for a test function $G$; see Prop.~\ref{rem:uniform}. \\
\end{xltabular}

\bibliographystyle{amsalpha}
\bibliography{References}

@unpublished {SWZ,
    AUTHOR = {Sun, Lehan and Wen, Yufei and Zhang, Xudong},
     TITLE = {Remark on the paper ``{F}luctuations in the distribution of
              {H}ecke eigenvalues about the {S}ato--{T}ate measure''},
      NOTE = {Unpublished manuscript},
      YEAR = {2018},
}

@article {KWY2,
    AUTHOR = {Kim, Henry H. and Wakatsuki, Satoshi and Yamauchi, Takuya},
     TITLE = {Central limit theorem for {H}ecke eigenvalues},
   JOURNAL = {Forum Math.},
  FJOURNAL = {Forum Mathematicum},
    VOLUME = {38},
      YEAR = {2026},
    NUMBER = {3},
     PAGES = {793--808},
      ISSN = {0933-7741,1435-5337},
   MRCLASS = {11F66 (11F46 11N75 22E46)},
  MRNUMBER = {5012880},
MRREVIEWER = {Hideshi\ Takayanagi},
       DOI = {10.1515/forum-2025-0044},
       URL = {https://doi.org/10.1515/forum-2025-0044},
}

@article {KWY1,
    AUTHOR = {Kim, Henry H. and Wakatsuki, Satoshi and Yamauchi, Takuya},
     TITLE = {Distribution of {H}ecke eigenvalues for holomorphic {S}iegel
              modular forms},
   JOURNAL = {Acta Arith.},
  FJOURNAL = {Acta Arithmetica},
    VOLUME = {215},
      YEAR = {2024},
    NUMBER = {2},
     PAGES = {161--177},
      ISSN = {0065-1036,1730-6264},
   MRCLASS = {11N75 (11F46)},
  MRNUMBER = {4779685},
MRREVIEWER = {Hideshi\ Takayanagi},
       DOI = {10.4064/aa230831-6-5},
       URL = {https://doi.org/10.4064/aa230831-6-5},
}

@book {Knightly-Li,
    AUTHOR = {Knightly, Andrew and Li, Charles},
     TITLE = {Traces of {H}ecke operators},
    SERIES = {Mathematical Surveys and Monographs},
    VOLUME = {133},
 PUBLISHER = {American Mathematical Society, Providence, RI},
      YEAR = {2006},
     PAGES = {x+378},
      ISBN = {978-0-8218-3739-9; 0-8218-3739-7},
   MRCLASS = {11F72 (11-02 11F25 11F70)},
  MRNUMBER = {2273356},
MRREVIEWER = {Mark\ Rowland\ Budden},
       DOI = {10.1090/surv/133},
       URL = {https://doi.org/10.1090/surv/133},
}

@article {Murty-Prabhu-BAuMS,
    AUTHOR = {Ram Murty, M. and Prabhu, Neha},
     TITLE = {The error term in the {S}ato-{T}ate theorem of {B}irch},
   JOURNAL = {Bull. Aust. Math. Soc.},
  FJOURNAL = {Bulletin of the Australian Mathematical Society},
    VOLUME = {100},
      YEAR = {2019},
    NUMBER = {1},
     PAGES = {27--33},
      ISSN = {0004-9727,1755-1633},
   MRCLASS = {11G05 (11K38)},
  MRNUMBER = {3977298},
MRREVIEWER = {\.{I}lker\ \.{I}nam},
       DOI = {10.1017/s0004972718001648},
       URL = {https://doi.org/10.1017/s0004972718001648},
}

@article {Banks-Shpar,
    AUTHOR = {Banks, William D. and Shparlinski, Igor E.},
     TITLE = {Sato-{T}ate, cyclicity, and divisibility statistics on average
              for elliptic curves of small height},
   JOURNAL = {Israel J. Math.},
  FJOURNAL = {Israel Journal of Mathematics},
    VOLUME = {173},
      YEAR = {2009},
     PAGES = {253--277},
      ISSN = {0021-2172,1565-8511},
   MRCLASS = {11G20 (11G50 11L40 11T24)},
  MRNUMBER = {2570668},
MRREVIEWER = {Steven\ D.\ Galbraith},
       DOI = {10.1007/s11856-009-0091-0},
       URL = {https://doi.org/10.1007/s11856-009-0091-0},
}

@article {Miller-Murty,
    AUTHOR = {Miller, Steven J. and Murty, M. Ram},
     TITLE = {Effective equidistribution and the {S}ato-{T}ate law for
              families of elliptic curves},
   JOURNAL = {J. Number Theory},
  FJOURNAL = {Journal of Number Theory},
    VOLUME = {131},
      YEAR = {2011},
    NUMBER = {1},
     PAGES = {25--44},
      ISSN = {0022-314X,1096-1658},
   MRCLASS = {11G05 (11G20 11K36)},
  MRNUMBER = {2729207},
MRREVIEWER = {D.\ R.\ Heath-Brown},
       DOI = {10.1016/j.jnt.2010.06.013},
       URL = {https://doi.org/10.1016/j.jnt.2010.06.013},
}

@article {Cho-Kim,
    AUTHOR = {Cho, Peter J. and Kim, Henry H.},
     TITLE = {Central limit theorem for {A}rtin {$L$}-functions},
   JOURNAL = {Int. J. Number Theory},
  FJOURNAL = {International Journal of Number Theory},
    VOLUME = {13},
      YEAR = {2017},
    NUMBER = {1},
     PAGES = {1--14},
      ISSN = {1793-0421,1793-7310},
   MRCLASS = {11N60 (11N37)},
  MRNUMBER = {3573409},
MRREVIEWER = {Yuan\ Yi},
       DOI = {10.1142/S1793042117500014},
       URL = {https://doi.org/10.1142/S1793042117500014},
}

@article {Winnie-Li,
    AUTHOR = {Li, Wen Ch'ing Winnie},
     TITLE = {Newforms and functional equations},
   JOURNAL = {Math. Ann.},
  FJOURNAL = {Mathematische Annalen},
    VOLUME = {212},
      YEAR = {1975},
     PAGES = {285--315},
      ISSN = {0025-5831,1432-1807},
   MRCLASS = {10D05},
  MRNUMBER = {369263},
MRREVIEWER = {R.\ A.\ Rankin},
       DOI = {10.1007/BF01344466},
       URL = {https://doi.org/10.1007/BF01344466},
}

@article {AL,
    AUTHOR = {Atkin, A. O. L. and Lehner, J.},
     TITLE = {Hecke operators on {$\Gamma _{0}(m)$}},
   JOURNAL = {Math. Ann.},
  FJOURNAL = {Mathematische Annalen},
    VOLUME = {185},
      YEAR = {1970},
     PAGES = {134--160},
      ISSN = {0025-5831,1432-1807},
   MRCLASS = {10.20},
  MRNUMBER = {268123},
MRREVIEWER = {R.\ A.\ Rankin},
       DOI = {10.1007/BF01359701},
       URL = {https://doi.org/10.1007/BF01359701},
}

@article {Das-Prabhu,
    AUTHOR = {Das, Jishu and Prabhu, Neha},
     TITLE = {A central limit theorem for {H}ilbert modular forms},
   JOURNAL = {J. Ramanujan Math. Soc.},
  FJOURNAL = {Journal of the Ramanujan Mathematical Society},
    VOLUME = {40},
      YEAR = {2025},
    NUMBER = {1},
     PAGES = {29--41},
      ISSN = {0970-1249,2320-3110},
   MRCLASS = {11F41 (11F30 11F72)},
  MRNUMBER = {4875173},
}

@article {Baier-Zhao,
    AUTHOR = {Baier, Stephan and Zhao, Liangyi},
     TITLE = {The {S}ato-{T}ate conjecture on average for small angles},
   JOURNAL = {Trans. Amer. Math. Soc.},
  FJOURNAL = {Transactions of the American Mathematical Society},
    VOLUME = {361},
      YEAR = {2009},
    NUMBER = {4},
     PAGES = {1811--1832},
      ISSN = {0002-9947,1088-6850},
   MRCLASS = {11G05},
  MRNUMBER = {2465818},
MRREVIEWER = {Hui\ Xue},
       DOI = {10.1090/S0002-9947-08-04498-X},
       URL = {https://doi.org/10.1090/S0002-9947-08-04498-X},
}

@article {Johnson,
    AUTHOR = {Johnson, Tobias},
     TITLE = {Exchangeable pairs, switchings, and random regular graphs},
   JOURNAL = {Electron. J. Combin.},
  FJOURNAL = {Electronic Journal of Combinatorics},
    VOLUME = {22},
      YEAR = {2015},
    NUMBER = {1},
     PAGES = {Paper 1.33, 28},
      ISSN = {1077-8926},
   MRCLASS = {05C80 (60B10)},
  MRNUMBER = {3315475},
       DOI = {10.37236/4659},
       URL = {https://doi.org/10.37236/4659},
}

@article {Baier-Prabhu,
    AUTHOR = {Baier, Stephan and Prabhu, Neha},
     TITLE = {Moments of the error term in the {S}ato-{T}ate law for
              elliptic curves},
   JOURNAL = {J. Number Theory},
  FJOURNAL = {Journal of Number Theory},
    VOLUME = {194},
      YEAR = {2019},
     PAGES = {44--82},
      ISSN = {0022-314X,1096-1658},
   MRCLASS = {11G05 (11G40)},
  MRNUMBER = {3860468},
MRREVIEWER = {Daniel\ Delbourgo},
       DOI = {10.1016/j.jnt.2018.07.005},
       URL = {https://doi.org/10.1016/j.jnt.2018.07.005},
}

@article {LNW1,
    AUTHOR = {Lau, Yuk-Kam and Ng, Ming Ho and Wang, Yingnan},
     TITLE = {Statistics of {H}ecke eigenvalues for {${\rm GL}(n)$}},
   JOURNAL = {Forum Math.},
  FJOURNAL = {Forum Mathematicum},
    VOLUME = {31},
      YEAR = {2019},
    NUMBER = {1},
     PAGES = {167--185},
      ISSN = {0933-7741,1435-5337},
   MRCLASS = {11F12 (11F30 11F72)},
  MRNUMBER = {3898221},
MRREVIEWER = {Eren\ Mehmet\ K\i ral},
       DOI = {10.1515/forum-2018-0166},
       URL = {https://doi.org/10.1515/forum-2018-0166},
}

@article {MMP2,
    AUTHOR = {Murty, M. Ram and Murty, V. Kumar and Pujahari, Sudhir},
     TITLE = {An all-purpose {E}rd\"os-{K}ac theorem},
   JOURNAL = {Math. Z.},
  FJOURNAL = {Mathematische Zeitschrift},
    VOLUME = {305},
      YEAR = {2023},
    NUMBER = {3},
     PAGES = {Paper No. 45, 18},
      ISSN = {0025-5874,1432-1823},
   MRCLASS = {11K65 (11N37 11N75)},
  MRNUMBER = {4658616},
MRREVIEWER = {Amrinder\ Kaur},
       DOI = {10.1007/s00209-023-03370-y},
       URL = {https://doi.org/10.1007/s00209-023-03370-y},
}

@article {Mahajan-Sinha,
    AUTHOR = {Mahajan, Jewel and Sinha, Kaneenika},
     TITLE = {Higher moments of the pair correlation function for
              {S}ato-{T}ate sequences},
   JOURNAL = {J. Number Theory},
  FJOURNAL = {Journal of Number Theory},
    VOLUME = {257},
      YEAR = {2024},
     PAGES = {24--97},
      ISSN = {0022-314X,1096-1658},
   MRCLASS = {11F11 (11F25 11F30 11K06)},
  MRNUMBER = {4671616},
MRREVIEWER = {Guodong\ Hua},
       DOI = {10.1016/j.jnt.2023.10.008},
       URL = {https://doi.org/10.1016/j.jnt.2023.10.008},
}

@article {Xi,
    AUTHOR = {Xi, Ping},
     TITLE = {Gaussian distributions of {K}loosterman sums: vertical and
              horizontal},
   JOURNAL = {Ramanujan J.},
  FJOURNAL = {Ramanujan Journal. An International Journal Devoted to the
              Areas of Mathematics Influenced by Ramanujan},
    VOLUME = {43},
      YEAR = {2017},
    NUMBER = {3},
     PAGES = {493--511},
      ISSN = {1382-4090},
   MRCLASS = {11L05 (60F05)},
  MRNUMBER = {3670279},
MRREVIEWER = {Florent Jouve},
       DOI = {10.1007/s11139-017-9905-6},
       URL = {https://doi.org/10.1007/s11139-017-9905-6},
}

@article {Nied,
    AUTHOR = {Niederreiter, Harald},
     TITLE = {The distribution of values of {K}loosterman sums},
   JOURNAL = {Arch. Math. (Basel)},
  FJOURNAL = {Archiv der Mathematik},
    VOLUME = {56},
      YEAR = {1991},
    NUMBER = {3},
     PAGES = {270--277},
      ISSN = {0003-889X},
   MRCLASS = {11L05 (11K38)},
  MRNUMBER = {1091880},
MRREVIEWER = {Jerzy Kaczorowski},
       DOI = {10.1007/BF01190214},
       URL = {https://doi.org/10.1007/BF01190214},
}

@book {Katz,
    AUTHOR = {Katz, Nicholas M.},
     TITLE = {Gauss sums, {K}loosterman sums, and monodromy groups},
    SERIES = {Annals of Mathematics Studies},
    VOLUME = {116},
 PUBLISHER = {Princeton University Press, Princeton, NJ},
      YEAR = {1988},
     PAGES = {x+246},
      ISBN = {0-691-08432-7; 0-691-08433-5},
   MRCLASS = {11G40 (11L05 11T24 14G10 14G15)},
  MRNUMBER = {955052},
MRREVIEWER = {P. Bayer},
       DOI = {10.1515/9781400882120},
       URL = {https://doi.org/10.1515/9781400882120},
}

@article {DP,
    AUTHOR = {David, Chantal and Pappalardi, Francesco},
     TITLE = {Average {F}robenius distributions of elliptic curves},
   JOURNAL = {Internat. Math. Res. Notices},
  FJOURNAL = {International Mathematics Research Notices},
      YEAR = {1999},
    NUMBER = {4},
     PAGES = {165--183},
      ISSN = {1073-7928,1687-0247},
   MRCLASS = {11G05 (11F30 11N36)},
  MRNUMBER = {1677267},
MRREVIEWER = {M.\ Ram\ Murty},
       DOI = {10.1155/S1073792899000082},
       URL = {https://doi.org/10.1155/S1073792899000082},
}

@article {Birch,
    AUTHOR = {Birch, B. J.},
     TITLE = {How the number of points of an elliptic curve over a fixed
              prime field varies},
   JOURNAL = {J. London Math. Soc.},
  FJOURNAL = {The Journal of the London Mathematical Society},
    VOLUME = {43},
      YEAR = {1968},
     PAGES = {57--60},
      ISSN = {0024-6107,1469-7750},
   MRCLASS = {10.14},
  MRNUMBER = {230682},
MRREVIEWER = {A.\ P.\ Ogg},
       DOI = {10.1112/jlms/s1-43.1.57},
       URL = {https://doi.org/10.1112/jlms/s1-43.1.57},
}

@article {Deuring,
    AUTHOR = {Deuring, Max},
     TITLE = {Die {T}ypen der {M}ultiplikatorenringe elliptischer
              {F}unktionenk\"{o}rper},
   JOURNAL = {Abh. Math. Sem. Hansischen Univ.},
  FJOURNAL = {Abhandlungen aus dem Mathematischen Seminar der Hansischen
              Universit\"{a}t},
    VOLUME = {14},
      YEAR = {1941},
     PAGES = {197--272},
      ISSN = {0025-5858},
   MRCLASS = {09.1X},
  MRNUMBER = {5125},
MRREVIEWER = {Saunders\ Mac Lane},
       DOI = {10.1007/BF02940746},
       URL = {https://doi.org/10.1007/BF02940746},
}

@article {JY,
    AUTHOR = {James, Kevin and Yu, Gang},
     TITLE = {Average {F}robenius distribution of elliptic curves},
   JOURNAL = {Acta Arith.},
  FJOURNAL = {Acta Arithmetica},
    VOLUME = {124},
      YEAR = {2006},
    NUMBER = {1},
     PAGES = {79--100},
      ISSN = {0065-1036,1730-6264},
   MRCLASS = {11G05 (11N05 11P55)},
  MRNUMBER = {2262142},
MRREVIEWER = {Francesco\ Pappalardi},
       DOI = {10.4064/aa124-1-7},
       URL = {https://doi.org/10.4064/aa124-1-7},
}

@article {FM,
    AUTHOR = {Fouvry, Etienne and Murty, M. Ram},
     TITLE = {On the distribution of supersingular primes},
   JOURNAL = {Canad. J. Math.},
  FJOURNAL = {Canadian Journal of Mathematics. Journal Canadien de
              Math\'{e}matiques},
    VOLUME = {48},
      YEAR = {1996},
    NUMBER = {1},
     PAGES = {81--104},
      ISSN = {0008-414X,1496-4279},
   MRCLASS = {11G05 (11N45)},
  MRNUMBER = {1382477},
MRREVIEWER = {Andrea\ Mori},
       DOI = {10.4153/CJM-1996-004-7},
       URL = {https://doi.org/10.4153/CJM-1996-004-7},
}

@article {Prabhu-Sinha,
    AUTHOR = {Prabhu, Neha and Sinha, Kaneenika},
     TITLE = {Fluctuations in the distribution of {H}ecke eigenvalues about
              the {S}ato-{T}ate measure},
   JOURNAL = {Int. Math. Res. Not. IMRN},
  FJOURNAL = {International Mathematics Research Notices. IMRN},
      YEAR = {2019},
    NUMBER = {12},
     PAGES = {3768--3811},
      ISSN = {1073-7928},
   MRCLASS = {11F11 (11F25 11F30)},
  MRNUMBER = {3973109},
MRREVIEWER = {Eric Stade},
       DOI = {10.1093/imrn/rnx238},
       URL = {https://doi.org/10.1093/imrn/rnx238},
}

@article {BPS,
    AUTHOR = {Baier, Stephan and Prabhu, Neha and Sinha, Kaneenika},
     TITLE = {Central limit theorems for elliptic curves and modular forms
              with smooth weight functions},
   JOURNAL = {J. Math. Anal. Appl.},
  FJOURNAL = {Journal of Mathematical Analysis and Applications},
    VOLUME = {485},
      YEAR = {2020},
    NUMBER = {1},
     PAGES = {123709, 29},
      ISSN = {0022-247X,1096-0813},
   MRCLASS = {11F11 (11F25 11F41 11G05 11G40)},
  MRNUMBER = {4042634},
MRREVIEWER = {Daniel\ Delbourgo},
       DOI = {10.1016/j.jmaa.2019.123709},
       URL = {https://doi.org/10.1016/j.jmaa.2019.123709},
}

@article {BGHT,
    AUTHOR = {Barnet-Lamb, Tom and Geraghty, David and Harris, Michael and
              Taylor, Richard},
     TITLE = {A family of {C}alabi-{Y}au varieties and potential automorphy
              {II}},
   JOURNAL = {Publ. Res. Inst. Math. Sci.},
  FJOURNAL = {Publications of the Research Institute for Mathematical
              Sciences},
    VOLUME = {47},
      YEAR = {2011},
    NUMBER = {1},
     PAGES = {29--98},
      ISSN = {0034-5318,1663-4926},
   MRCLASS = {11F80 (11F11 11G18 14J32)},
  MRNUMBER = {2827723},
MRREVIEWER = {Neil\ P.\ Dummigan},
       DOI = {10.2977/PRIMS/31},
       URL = {https://doi.org/10.2977/PRIMS/31},
}

@article {Serre,
    AUTHOR = {Serre, Jean-Pierre},
     TITLE = {R\'epartition asymptotique des valeurs propres de
              l'op\'erateur de {H}ecke {$T_p$}},
   JOURNAL = {J. Amer. Math. Soc.},
  FJOURNAL = {Journal of the American Mathematical Society},
    VOLUME = {10},
      YEAR = {1997},
    NUMBER = {1},
     PAGES = {75--102},
      ISSN = {0894-0347,1088-6834},
   MRCLASS = {11F30 (11F25 11G20 11N37 11R45)},
  MRNUMBER = {1396897},
MRREVIEWER = {Glenn\ Stevens},
       DOI = {10.1090/S0894-0347-97-00220-8},
       URL = {https://doi.org/10.1090/S0894-0347-97-00220-8},
}

@article {MS,
    AUTHOR = {Murty, M. Ram and Sinha, Kaneenika},
     TITLE = {Effective equidistribution of eigenvalues of {H}ecke
              operators},
   JOURNAL = {J. Number Theory},
  FJOURNAL = {Journal of Number Theory},
    VOLUME = {129},
      YEAR = {2009},
    NUMBER = {3},
     PAGES = {681--714},
      ISSN = {0022-314X},
   MRCLASS = {11F25 (11N75)},
   MRNUMBER = {2488597},
  MRREVIEWER = {D. R. Heath-Brown},
       DOI = {10.1016/j.jnt.2008.10.010},
       URL = {https://doi.org/10.1016/j.jnt.2008.10.010},
}

@article {H,
    AUTHOR = {Halberstam, H.},
     TITLE = {On the distribution of additive number-theoretic functions},
   JOURNAL = {J. London Math. Soc.},
  FJOURNAL = {The Journal of the London Mathematical Society},
    VOLUME = {30},
      YEAR = {1955},
     PAGES = {43--53},
      ISSN = {0024-6107},
   MRCLASS = {10.0X},
MRREVIEWER = {P. Erd\H{o}s},
       DOI = {10.1112/jlms/s1-30.1.43},
       URL = {https://doi.org/10.1112/jlms/s1-30.1.43},
}

@article {Nagoshi,
    AUTHOR = {Nagoshi, Hirofumi},
     TITLE = {Distribution of {H}ecke eigenvalues},
   JOURNAL = {Proc. Amer. Math. Soc.},
  FJOURNAL = {Proceedings of the American Mathematical Society},
    VOLUME = {134},
      YEAR = {2006},
    NUMBER = {11},
     PAGES = {3097--3106},
      ISSN = {0002-9939,1088-6826},
   MRCLASS = {11F30 (11F11 11F25)},
  MRNUMBER = {2231890},
MRREVIEWER = {Scott\ Ahlgren},
       DOI = {10.1090/S0002-9939-06-08709-0},
       URL = {https://doi.org/10.1090/S0002-9939-06-08709-0},
}

@book {Billingsley,
    AUTHOR = {Billingsley, Patrick},
     TITLE = {Probability and measure},
    SERIES = {Wiley Series in Probability and Statistics},
   EDITION = {Anniversary},
      NOTE = {With a foreword by Steve Lalley and a brief biography of
              Billingsley by Steve Koppes},
 PUBLISHER = {John Wiley \& Sons, Inc., Hoboken, NJ},
      YEAR = {2012},
     PAGES = {xviii+624},
      ISBN = {978-1-118-12237-2},
   MRCLASS = {60-01 (28-01)},
  MRNUMBER = {2893652},
}

@article {murty-prabhu,
    AUTHOR = {Murty, M. Ram and Prabhu, Neha},
     TITLE = {Central limit theorems for sums of quadratic characters,
              {H}ecke eigenforms, and elliptic curves},
   JOURNAL = {Proc. Amer. Math. Soc.},
  FJOURNAL = {Proceedings of the American Mathematical Society},
    VOLUME = {148},
      YEAR = {2020},
    NUMBER = {3},
     PAGES = {965--977},
      ISSN = {0002-9939,1088-6826},
   MRCLASS = {11F30 (11G05 11N37)},
  MRNUMBER = {4055927},
MRREVIEWER = {Francesc\ Fit\'e},
       DOI = {10.1090/proc/14760},
       URL = {https://doi.org/10.1090/proc/14760},
}

@article {TANG2016295,
    AUTHOR = {Tang, Hengcai and Wang, Yingnan},
     TITLE = {Quantitative versions of the joint distributions of {H}ecke
              eigenvalues},
   JOURNAL = {J. Number Theory},
  FJOURNAL = {Journal of Number Theory},
    VOLUME = {169},
      YEAR = {2016},
     PAGES = {295--314},
      ISSN = {0022-314X,1096-1658},
   MRCLASS = {11F11 (11F25)},
  MRNUMBER = {3531241},
MRREVIEWER = {Thomas\ R.\ Shemanske},
       DOI = {10.1016/j.jnt.2016.05.011},
       URL = {https://doi.org/10.1016/j.jnt.2016.05.011},
}

@article {A_Weil,
    AUTHOR = {Weil, Andr\'e},
     TITLE = {On some exponential sums},
   JOURNAL = {Proc. Nat. Acad. Sci. U.S.A.},
  FJOURNAL = {Proceedings of the National Academy of Sciences of the United
              States of America},
    VOLUME = {34},
      YEAR = {1948},
     PAGES = {204--207},
      ISSN = {0027-8424},
   MRCLASS = {10.0X},
  MRNUMBER = {27006},
MRREVIEWER = {H.\ D.\ Kloosterman},
       DOI = {10.1073/pnas.34.5.204},
       URL = {https://doi.org/10.1073/pnas.34.5.204},
}

@book {KN,
    AUTHOR = {Kuipers, L. and Niederreiter, H.},
     TITLE = {Uniform distribution of sequences},
    SERIES = {Pure and Applied Mathematics},
 PUBLISHER = {Wiley-Interscience [John Wiley \& Sons], New
              York-London-Sydney},
      YEAR = {1974},
     PAGES = {xiv+390},
   MRCLASS = {10K05 (22D99)},
  MRNUMBER = {419394},
MRREVIEWER = {P.\ Gerl},
}

@article {Clozel-Harris-Taylor,
    AUTHOR = {Clozel, Laurent and Harris, Michael and Taylor, Richard},
     TITLE = {Automorphy for some {$l$}-adic lifts of automorphic mod {$l$} {G}alois representations},
      NOTE = {With Appendix A, summarizing unpublished work of Russ Mann,
              and Appendix B by Marie-France Vign\'eras},
   JOURNAL = {Publ. Math. Inst. Hautes \'Etudes Sci.},
  FJOURNAL = {Publications Math\'ematiques. Institut de Hautes \'Etudes
              Scientifiques},
    NUMBER = {108},
      YEAR = {2008},
     PAGES = {1--181},
      ISSN = {0073-8301,1618-1913},
   MRCLASS = {11F80 (11G18 11R34)},
  MRNUMBER = {2470687},
MRREVIEWER = {Mark\ Kisin},
       DOI = {10.1007/s10240-008-0016-1},
       URL = {https://doi-org.libproxy1.nus.edu.sg/10.1007/s10240-008-0016-1},
}

@incollection {Tate,
    AUTHOR = {Tate, John T.},
     TITLE = {Algebraic cycles and poles of zeta functions},
 BOOKTITLE = {Arithmetical {A}lgebraic {G}eometry ({P}roc. {C}onf. {P}urdue
              {U}niv., 1963)},
     PAGES = {93--110},
 PUBLISHER = {Harper \& Row, New York},
      YEAR = {1965},
   MRCLASS = {14.05},
  MRNUMBER = {225778},
MRREVIEWER = {A.\ Mattuck},
}

@book{Mont,
	author = {Montgomery, Hugh L.},
	doi = {10.1090/cbms/084},
	isbn = {0-8218-0737-4},
	mrclass = {11-02 (11Kxx 11L07 11Mxx 11Nxx)},
	mrnumber = {1297543},
	mrreviewer = {John B. Friedlander},
	pages = {xiv+220},
	publisher = {Published for the Conference Board of the Mathematical Sciences, Washington, DC; by the American Mathematical Society, Providence, RI},
	series = {CBMS Regional Conference Series in Mathematics},
	title = {Ten lectures on the interface between analytic number theory and harmonic analysis},
	url = {http://dx.doi.org/10.1090/cbms/084},
	volume = {84},
	year = {1994}}

@article{LW,
	author = {Lau, Yuk-Kam and Wang, Yingnan},
	doi = {10.1016/j.jnt.2011.05.014},
	fjournal = {Journal of Number Theory},
	issn = {0022-314X},
	journal = {J. Number Theory},
	mrclass = {11F25 (11N75)},
	mrnumber = {2832823},
	mrreviewer = {Lei Yang},
	number = {12},
	pages = {2262--2281},
	title = {Quantitative version of the joint distribution of eigenvalues of the {H}ecke operators},
	url = {http://dx.doi.org/10.1016/j.jnt.2011.05.014},
	volume = {131},
	year = {2011}}

@article{FR,
	author = {Faifman, Dmitry and Rudnick, Ze\'ev},
	doi = {10.1112/S0010437X09004308},
	fjournal = {Compositio Mathematica},
	issn = {0010-437X},
	journal = {Compos. Math.},
	mrclass = {11M50 (11G20 11M36)},
	mrnumber = {2581242},
	mrreviewer = {Andrew V. Sutherland},
	number = {1},
	pages = {81--101},
	title = {Statistics of the zeros of zeta functions in families of hyperelliptic curves over a finite field},
	url = {http://dx.doi.org/10.1112/S0010437X09004308},
	volume = {146},
	year = {2010}}

@article{CDF,
	author = {Conrey, J. B. and Duke, W. and Farmer, D. W.},
	fjournal = {Acta Arithmetica},
	issn = {0065-1036},
	journal = {Acta Arith.},
	mrclass = {11F25 (11F11)},
	mrnumber = {1438595},
	mrreviewer = {Tomasz Stefanicki},
	number = {4},
	pages = {405--409},
	title = {The distribution of the eigenvalues of {H}ecke operators},
	volume = {78},
	year = {1997}}

@article{BDFLS,
	author = {Bucur, Alina and David, Chantal and Feigon, Brooke and Lal\'\i n, Matilde and Sinha, Kaneenika},
	doi = {10.4310/MRL.2012.v19.n6.a12},
	fjournal = {Mathematical Research Letters},
	issn = {1073-2780},
	journal = {Math. Res. Lett.},
	mrclass = {11G20 (11M50 14G15)},
	mrnumber = {3091611},
	mrreviewer = {Florin Spinu},
	number = {6},
	pages = {1329--1356},
	title = {Distribution of zeta zeroes of {A}rtin-{S}chreier covers},
	url = {http://dx.doi.org/10.4310/MRL.2012.v19.n6.a12},
	volume = {19},
	year = {2012}}

@book{Serre1,
	author = {Serre, Jean-Pierre},
	edition = {Second},
	isbn = {0-201-09384-7},
	mrclass = {11G05 (11F33 11F80)},
	mrnumber = {1043865},
	note = {With the collaboration of Willem Kuyk and John Labute},
	pages = {xxiv+184},
	publisher = {Addison-Wesley Publishing Company, Advanced Book Program, Redwood City, CA},
	series = {Advanced Book Classics},
	title = {Abelian {$l$}-adic representations and elliptic curves},
	year = {1989}}

@incollection{Sarnak,
	author = {Sarnak, Peter},
	booktitle = {Analytic number theory and {D}iophantine problems ({S}tillwater, {OK}, 1984)},
	mrclass = {11F30 (11F37 11F55 11F72)},
	mrnumber = {1018385},
	mrreviewer = {V. Kumar Murty},
	pages = {321--331},
	publisher = {Birkh\"auser Boston, Boston, MA},
	series = {Progr. Math.},
	title = {Statistical properties of eigenvalues of the {H}ecke operators},
	volume = {70},
	year = {1987}}

@book{Rohatgi-Saleh,
	Author = {Rohatgi, Vijay K. and Ehsanes Saleh, A. K. Md.},
	Edition = {Second},
	Isbn = {0-471-34846-5},
	Mrclass = {60-01 (62-01)},
	Mrnumber = {1789794},
	Pages = {xviii+716},
	Publisher = {Wiley-Interscience, New York},
	Series = {Wiley Series in Probability and Statistics: Texts and References Section},
	Title = {An introduction to probability and statistics},
	Year = {2001}}

@article {RT,
    AUTHOR = {Rouse, Jeremy and Thorner, Jesse},
     TITLE = {The explicit {S}ato-{T}ate conjecture and densities pertaining
              to {L}ehmer-type questions},
   JOURNAL = {Trans. Amer. Math. Soc.},
  FJOURNAL = {Transactions of the American Mathematical Society},
    VOLUME = {369},
      YEAR = {2017},
    NUMBER = {5},
     PAGES = {3575--3604},
      ISSN = {0002-9947,1088-6850},
   MRCLASS = {11F30 (11F33 11M41)},
  MRNUMBER = {3605980},
MRREVIEWER = {Hidenori\ Katsurada},
       DOI = {10.1090/tran/6793},
       URL = {https://doi.org/10.1090/tran/6793},
}

@unpublished{L,
	author = {Lester, S.},
	note = {Sieve Theory and its applications, unpublished lecture notes (Link: http://www.math.tau.ac.il/~rudnick/courses/sieves2015.html)},
	title = {Normal order of $\omega(n)$ and the {E}rd\"{o}s-{K}ac theorem},
	urldate = {http://www.math.tau.ac.il/~rudnick/courses/sieves2015.html},
	year = {2015}}
\end{document}